\documentclass[leqno]{amsart}
\usepackage{latexsym,amssymb,amsthm,amsmath,amscd}

\usepackage{latexsym,amssymb,amsthm,amsmath,amscd}
\theoremstyle{plain}  
\newtheorem{theorem}{Theorem}[section]

\newtheorem{lemma}{Lemma}[section]
\newtheorem{proposition}{Proposition}[section]

\newtheorem{corollary}{Corollary}[section]

\newtheorem{definition}{Definition}[section] 
\newtheorem{example}{Example}[section]
\numberwithin{equation}{section}

\theoremstyle{remark}
\newtheorem{remark}{Remark}[section]
 \newtheorem{note}[theorem]{Note}

 \numberwithin{equation}{section}

\def\<{\left < }
\def\>{\right >}
\def\({\left ( }
\def\){\right )}

\begin{document}

\title[Golden Riemannian manifolds]{A comprehensive review of golden Riemannian manifolds}

\author[B.-Y. Chen, M. A. Choudhary and A Perwee]{Bang-Yen Chen, Majid Ali Choudhary and Afshan Perween}

\begin{abstract} In differential geometry, the concept of golden structure, initially proposed by S. I. Goldberg and K. Yano in 1970, presents a compelling area with wide-ranging applications. The exploration of golden Riemannian manifolds was initiated by C. E. Hretcanu and M. Crasmareanu in 2008, following the principles of the golden structure. Subsequently, numerous researchers have contributed significant insights into golden Riemannian manifolds. The purpose of this paper is to provide a comprehensive survey on golden Riemannian manifold done over the past decade.
\end{abstract}


\maketitle

\section{Introduction}

For nearly two thousand years, the concept of the golden ratio has fascinated scholars from diverse disciplines.  The allure of the golden ratio extends beyond mathematicians; it captivates biologists, artists, musicians, historians, architects, and psychologists alike. The golden ratio, known for its aesthetic harmony and proportionality, is extensively employed in iconic architectural structures and artworks, musical composition frameworks, harmonious frequency ratios, and human body measurements. It is likely fair to say that the golden ratio has inspired more scholars in various fields than any other number throughout the history of mathematics \cite{79}.

Polynomial structures on a manifold, discussed in \cite{1}, were the foundation for the concept of golden structure. C. E.  Hretcanu and M.  Crasmareanu investigated some characteristics of the induced structure on an invariant submanifold within a golden Riemannian manifold in \cite{17}. In \cite{2}, Crasmareanu and Hretcanu used a corresponding almost product structure to study the geometry of the golden structure on a manifold.  Hretcanu and Crasmareanu demonstrated in \cite{78} that a golden structure also induces a golden structure on each invariant submanifold. The issue of integrability for golden Riemannian structures was examined by Gezer et al. in \cite{73}. Ozkan studied the golden semi-Riemannian manifold in \cite{80}, where he defined the golden structure's horizontal lift in the tangent bundle. Some more structures on golden Riemannian manifolds were studied by many authors (see, e.g., \cite{76,75, 77}).

However, in the study of differential geometry, the theory of submanifolds is an intriguing subject. Its roots are in Fermat's work on the geometry of surfaces and plane curves. Since then, it has been developed in different directions of differential geometry and mechanics, specifically. It is still a vibrant area of study that has contributed significantly to the advancement of differential geometry in the present era. Among all the submanifolds of an ambient manifold, there are two well-known types of submanifolds: invariant submanifolds and anti-invariant ones.
The differential geometry of submanifolds in golden Riemannian manifolds was initially studied  by Crasmareanu and Hretcanu. Certain characteristics of invariant submanifolds in a Riemannian manifold with golden structure were investigated by Hretcanu and Crasmareanu in \cite{17}, which had been advanced greatly since then. 
Various classifications of submanifolds within golden Riemannian manifolds have been established based on how their tangent bundles react to the golden structure of the ambient manifold, and have been explored by numerous geometers.  Erdoğan and  Yıldırım presented the idea of semi-invariant submanifolds within golden Riemannian manifolds as a generalization of both invariant and anti-invariant types, followed by an analysis of the geometry of their defining distributions \cite{21}. The properties and distributions associated with semi-invariant submanifolds in golden Riemannian manifolds were further explored by  Gök,  Keleş, and  Kılıç \cite{18}. The notion of pointwise slant submanifolds, pointwise bi-slant submanifolds in golden Riemannian manifolds were introduced by Hretcanu and Blaga in \cite{13}. 

R. L. Bishop and B. O'Neill \cite{12} proposed the idea of a warped product. Warped product CR-submanifolds in a K\"ahler manifold, which consist of warped products of holomorphic and totally real submanifolds, were first studied by B.-Y. Chen in \cite{82,83,81,c03,c04,book17}. Later, Riemannian manifolds with golden warped products were studied by Blaga and Hretcanu in \cite{11}, who also investigated submanifolds with pointwise semi-slant and hemi-slant warped products within locally golden Riemannian manifolds \cite{13}.

The similarities between the geometry of semi-Riemannian submanifolds and their Riemannian case are well established, yet the study of lightlike submanifolds presents unique challenges due to the intersection of their normal vector bundle with the tangent bundle. This complexity adds to the intrigue of researching lightlike geometry, which finds practical applications inside mathematical physics, notably in general relativity and electromagnetism \cite{33}. Duggal and Bejancu were pioneers in the study of lightlike submanifolds within semi-Riemannian manifolds \cite{33}. In 2017, Poyraz and Yasar explored lightlike hypersurfaces in golden semi-Riemannian manifolds \cite{5}, and further extended their research to define lightlike submanifolds of golden semi-Riemannian manifolds in 2019 \cite{36}. Subsequently, numerous researchers have investigated various types of lightlike submanifolds in golden semi-Riemannian manifolds, as evidenced by references \cite{44, 40, 41, 55, 37, 39} among others. Additionally, the concept of a lightlike hypersurface in Meta-golden Riemannian manifolds was introduced by Erdoğan et al. \cite{15}.

Finding the optimal inequality between the intrinsic and extrinsic invariants of a Riemannian submanifold is a key problem in submanifold geometry. Chen \cite{58, 84} developed the $\delta$-invariants in this context, which are known today as Chen invariants. Utilizing these invariants along with the mean curvature, which is the key extrinsic invariant of Riemannian submanifolds, he formulated sharp inequality relationships, well-known as Chen inequalities. The study of Chen invariants and Chen inequalities across various submanifolds in diverse ambient spaces has been thoroughly pursued (see, e.g.,  \cite{85,84,87}). Research on Chen inequalities within golden Riemannian manifolds and golden-like statistical manifolds was conducted by Choudhary and Uddin \cite{59} and Bahadir et al. \cite{8}, respectively.

Thanks to F. Casorati \cite{60}, the widely accepted Gaussian curvature can be substituted with the Casorati curvature. Within the framework of Casorati curvatures for submanifolds across diverse ambient spaces, geometric inequalities have been formulated. The rationale behind introducing this curvature by F. Casorati is that it disappears precisely when both principal curvatures of a surface in $\mathbb E^3$ are zero, aligning more closely with the typical understanding of curvature. Numerous scholars have explored Casorati curvatures to derive sharp inequalities for specific submanifolds in varied golden ambient spaces (refer to \cite{66,Chen21,64,63}, etc.). Moreover, P. Wintgen \cite{67} introduced a significant geometric inequality in 1979, involving Gauss curvature, normal curvature, and square mean curvature,  known as Wintgen's inequality \cite{Chen2}. For slant, invariant, C-totally real, and Lagrangian submanifolds in golden Riemannian spaces, generalized Wintgen type inequalities were introduced by Choudhary et al. \cite{70}.

Motivated by the previously mentioned advancements in the subject, the purpose of this paper is to provide a comprehensive survey of the latest progress on golden Riemannian manifold achieved over the past decade.

\section{Preliminaries}\label{S2}

\subsection{Golden Riemannian manifolds}

Take a $(n+m)$-dimensional Riemannian manifold $(\bar{M}, \bar{g})$, and $F$ as a $(1, 1)$-tensor field on $\bar{M}$. Assume that $F$ satisfies the equation.

$$L(X) = X^n + a_nX^{n-1} + \ldots + a_2X + a_1I = 0,$$
where the identity transformation is denoted by $I$ and (for $X = F$) $F^{n-1}(p), F^{n-2}(p)$, $\ldots, F(p), I$ are linearly independent at each point $p$ on $\bar{M}$. Consequently, the polynomial $L(X)$ is referred to as the structure polynomial \cite{25,2,1}. Choosing the structure polynomial as $$L(X) = X^2 +I, \; \text {structure is almost complex;}$$ 
$$L(X) = X^2 -I,\;  \text {structure is almost product;}$$ 
$$L(X) = X^2 , \; \text {structure is almost  tangent.} $$

\begin{definition} 
Suppose a $(n+m)$-dimensional Riemannian manifold
$(\bar{M}, \bar{g})$ with $(1, 1)$-tensor field $\phi$ on $\bar{M}$. If the given equation holds $$\phi^2 - \phi - I = 0,$$ 
The tensor field $\phi$ on $\bar{M}$ is then referred to as a golden structure. When the Riemannian metric $\bar{g}$ is compatible with $\phi$, the manifold $(\bar{M}, \bar{g}, \phi)$ is referred to as a {\rm golden Riemannian manifold} \cite{25,1, 4}. Compatibility of the metric $\bar{g}$ with $\phi$ is defined by the following condition
 $$\bar{g}(\phi X, Y) = \bar{g}(X, \phi Y)$$ $\forall X, Y \in \Gamma(T\bar{M}),$ where $\Gamma(T\bar{M})$ is the set of all vector fields on $\bar{M}$. If we replace $X$in the foregoing equation with $\phi X$, we obtain
 $$\bar{g}(\phi X, \phi Y) = \bar{g}(\phi^2 X, Y) = \bar{g}(\phi X, Y) + \bar{g}(X, Y).$$ Consider a manifold $\bar{M}$ of dimension $(n + m)$ endowed with a type $(1,1)$ tensor field $F$, where $F^2 = I$ but $F$ is not equal to $\pm I$. This tensor field $F$ is known as an almost product structure. If this structure $F$ admits a Riemannian metric $\bar{g}$ satisfying $$\bar{g}(F X, Y) = \bar{g}(X, F Y), \quad \forall X, Y \in \Gamma(T \bar{M}),$$ 
 hence, $(\bar{M}, \bar g )$ is called an almost product Riemannian manifold. Moreover, a golden structure can be produced by an almost product structure $F$ in the manner described below.
 $$\phi = \frac{1}{2}(I + \sqrt{5}F).$$ On the contrary, $\phi$ is a golden structure, then
 $$F = \frac{1}{\sqrt{5}}(2\phi - I)$$ is an almost product structure \cite{25, 2}.
\end{definition}

\begin{example} {\rm  \cite{78}
 Take the Euclidean 6-space $\mathbb R^6$ with conventional coordinates $\left(\zeta_1, \zeta_2, \zeta_3, \zeta_4, \zeta_5, \zeta_6\right)$ and suppose $\phi:\mathbb  R^6 \rightarrow \mathbb R^6$ symbolizes $(1,1)$-tensor field defined by
$$
\phi\left(\zeta_1, \zeta_2, \zeta_3, \zeta_4, \zeta_5, \zeta_6\right)=\left(\psi \zeta_1, \psi \zeta_2, \psi \zeta_3,(1-\psi) \zeta_4,(1-\psi) \zeta_5,(1-\psi) \zeta_6\right)
$$
for any vector field $\left(\zeta_1, \zeta_2, \zeta_3, \zeta_4, \zeta_5, \zeta_6\right) \in \mathbb R^6$, where $\psi=\frac{1+\sqrt{5}}{2}$ and $1-\psi=\frac{1-\sqrt{5}}{2}$ are the roots of the equation $\zeta^2=\zeta+1$. Then one can get
$$
\begin{aligned}
\phi^2\left(\zeta_1, \zeta_2, \zeta_3, \zeta_4, \zeta_5, \zeta_6\right)= & \left(\psi^2 \zeta_1, \psi^2 \zeta_2, \psi^2 \zeta_3,(1-\psi)^2 \zeta_4,(1-\psi)^2 \zeta_5,(1-\psi)^2 \zeta_6\right) \\
= & \left(\psi \zeta_1, \psi \zeta_2, \psi \zeta_3,(1-\psi) \zeta_4,(1-\psi) \zeta_5,(1-\psi) \zeta_6\right) \\
& +\left(\zeta_1, \zeta_2, \zeta_3, \zeta_4, \zeta_5, \zeta_6\right) .
\end{aligned}
$$

Thus, $\phi^2-\phi-I=0$. Moreover,
$$
\langle\phi\left(\zeta_1, \ldots, \zeta_6\right),\left(\Omega_1, \ldots, \Omega_6\right) \rangle
=\langle \left(\zeta_1, \ldots, \zeta_6\right), \phi\left(\Omega_1, \ldots, \Omega_6\right) \rangle,
$$
for every vector fields $\left(\zeta_1, \zeta_2, \zeta_3, \zeta_4, \zeta_5, \zeta_6\right),\left(\Omega_1, \Omega_2, \Omega_3, \Omega_4, \Omega_5, \Omega_6\right) \in {\mathbb R^6}$,  where the standard metric on $\mathbb R^6$ is denoted by $\langle\,, \rangle$. As a result, the manifold $\left(\mathbb R^6,\langle\,, \rangle, \phi\right)$ is a golden riemannian manifold.}
\end{example}

\subsection{Golden-like statistical manifolds}\label{S2.2} 

 An interesting characteristic of golden structures is that they always occur in pairs; for instance, if $\phi$ represents a golden structure, then $\phi^* = I - \phi$ also qualifies as one. This phenomenon is similarly observed in almost tangent structures $(J$ and $-J)$ and almost complex structures $(F$ and $-F)$. Consequently, exploring the relationship between golden and product structures becomes a pertinent question. Consider $\bar{M}$ as a Riemannian manifold and let $\nabla$ represent a torsion-free affine connection. The configuration $(\bar{M}, \bar{g}, \nabla)$ is termed a statistical manifold when $\nabla \bar{g}$ is symmetric. Another affine connection $\nabla^*$ is defined by
$$ 
X\bar{g}(Y,Z) = \bar{g}(\nabla_X Y, Z) + \bar{g} (\nabla^*_{X} Z, Y)
$$
for vector fields $X, Y,$ and  $Z$ on $\bar{M}$. The affine connection $\nabla^*$ is referred to as the conjugate (or dual) to $\nabla$ relative to $\bar{g}$. This connection $\nabla^*$ is torsion-free, ensures that $\nabla^* \bar{g}$ is symmetric, and complies with the equation $\nabla^* =\frac{\nabla +\nabla^*}{2}$. It is evident that the set $(\bar M,\nabla^*, \bar g,)$ forms a statistical manifold. The curvature tensors $R$ and $R^*$ on $\bar M$ correspond to the affine connections $\nabla$ and $\nabla^*$, respectively. Additionally, the curvature tensor field $R^0$ linked with $\nabla^0$ is known as the Riemannian curvature tensor. Then \cite{2, 4, 5} $$\bar g(R(X,Y)Z,W)= -\bar g(Z,R^* (X,Y)W)$$ for vector fields $X, Y, Z,$ and $W$ on $\bar M$, where 
$$R(X,Y)Z = [\nabla_X ,\nabla_Y] Z - \nabla_{[X,Y]} Z.$$ 

Generally, since the dual connections are not metric, the sectional curvature is unable to be defined in a statistical setting as it is in semi-Riemannian geometry. Hence, Opozda introduced two types of sectional curvature on statistical manifolds (See \cite{6,7}).

Considering that $\pi$ is a plane section in $T\bar M$ with an orthonormal basis $\{X,Y\}$, where $\bar M$ is a statistical manifold. The definition of the sectional $K$-curvature is \cite{6}
 $$K(\pi)=\frac{1}{2} [\bar g (R(X,Y)Y,X)+\bar g (R^*(X,Y)Y,X)-\bar g(R^0(X,Y)Y,X).$$
 
\begin{definition} {\rm \cite{8}}
Let ($\bar M,\bar g , \phi$) be a golden semi-Riemannian manifold endowed with a tensor field $\phi^*$ of type
$(1,1)$ satisfying
\begin{equation} \label{2.1}
    \bar g (\phi X,Y)=\bar g (X,\phi^* Y)
\end{equation}
for vector fields $X$ and $Y$. In view of (1), we easily derive
\begin{equation} \label{2.2}
(\phi^*)^2 X=\phi^* X+X,
\end{equation}
\begin{equation} \label{2.3}
\bar g (\phi X,\phi^* Y)=\bar g (\phi X,Y)+\bar g (X,Y).
\end{equation}
Then ($\bar M,\bar g , \phi$) is called golden-like statistical manifold.

According to \eqref{2.2} and \eqref{2.3}, the tensor fields $\phi+\phi^*$ and $\phi-\phi^*$ are symmetric and skew symmetric with respect to $\bar g$, respectively. The equations \eqref{2.1}, \eqref{2.2}, and \eqref{2.3} imply the following proposition.
\end{definition}
\begin{proposition} {\rm \cite{8}}
($\bar M,\bar g , \phi$) is a golden-like statistical manifold if and only if it is ($\bar M,\bar g , \phi^*$).
\end{proposition}

\begin{note}{\rm
    Remark that if one chooses $\phi=\phi^*$ in a golden-like statistical manifold, then we have a golden semi-Riemannian manifold.}
\end{note}

\subsection{Golden Lorentzian manifolds}\label{S2.3} 

For the locally golden space form ${\bar{M}}=\bar M_p (c_p)\times \bar M_q (c_q)$, where $c_p$ and $c_q$ are constant sectional curvatures of Riemannian manifolds $\bar M_p$ and $\bar M_q$, the Riemannian curvature tensor $R$ is expressed as follows in \cite{9}:
 
\begin{align*}
R (X,Y)Z =& \frac{(\mp \sqrt{5}+3)c_p + (\pm \sqrt{5}+3)c_q}{10} [\bar g (Y,Z)X-\bar g (X,Z)Y]\\
& +\frac{(\mp \sqrt{5}-1)c_p + (\pm \sqrt{5}-1)c_q}{10} [\bar g (\phi Y,Z)X-\bar g (\phi X,Z)Y \\ 
&+\bar g(Y,Z)\phi X-\bar g(X,Z)\phi Y]+\frac{c_p+c_q}{5}[\bar g (\phi Y,Z)\phi X-\bar g (\phi X,Z)\phi Y]
\end{align*}
where $X,Y$ and $Z \in \Gamma(TN).$

The golden Lorentzian manifold is defined as

\begin{definition} {\rm \cite{10}
Let us consider a semi-Riemannian manifold $(\bar M^n, \bar g)$ where $\bar g$ has the signature $(-,+,+,...,+(n-1$  times$))$. Then $\bar M$ stands for golden Lorentzian manifold if it is endowed with a golden structure $\phi$ and $\bar g$ is $\phi$-compatible.}
\end{definition}

\begin{example}{\rm \cite{10}
Let $\mathbb{R}_1^3$ represent the semi-Euclidean space and consider the signature of $\bar g$ as $(-,+,+)$. If $\phi$ stands for a $(1,1)$ tensor field, then it is easy to show that if
$$
\phi\left(\zeta_1, \zeta_2, \zeta_3\right)=\frac{1}{2}\left(\zeta_1+\sqrt{5} \zeta_2, \zeta_2+\sqrt{5} \zeta_1, 2 \psi \zeta_3\right),
$$
for any vector field $\left(\zeta_1, \zeta_2, \zeta_3\right) \in \mathbb{R}_1^3$, where $\psi=\frac{1+\sqrt{5}}{2}$ is the golden mean, then
$$
\phi^2=\phi+I,
$$
and hence, $\phi$ is a golden structure on $\mathbb{R}_1^3$. Moreover, $\bar g$ may also be verified to be $\phi$-compatible. Thus, $\left(\mathbb{R}_1^3,\bar g, \phi\right)$ becomes a golden Lorentzian manifold.}

\end{example}

\subsection{Meta-golden semi-Riemannian manifolds}\label{S2.4} 

The following structure is comparable to the golden Ratio (see \cite{14}): In \cite{15}, authors have obtained $\chi =\frac{1}{\psi} +\frac{1}{{\chi'}}$ where meta-golden Chi ratio $\chi =\frac{1+\sqrt{4\psi +5}}{2\psi}$  and $\psi=\frac{1+\sqrt{5}}{2} $, which suggests that $\chi^2-\frac{1}{\psi} \chi-1=0$. Thus, the roots are found as $\frac{\frac{1}{\psi} \pm \sqrt{4+\frac{1}{\psi^2}}}{2}$.
The correlation between the meta-golden Chi ratio $\chi$ and continued fractions was found in \cite{14}. Denoting the positive and negative roots by $\chi = \frac{\frac{1}{\psi} +\sqrt{4+\frac{1}{\psi^2}}}{2}$ and $\chi^. = \frac{\frac{1}{\psi} -\sqrt{4+\frac{1}{\psi^2}}}{2}$, respectively we have \cite{14}
\begin{equation}
    \chi^.=\frac{1}{\psi} \chi,
\end{equation}
\begin{equation}
    \psi \chi^2 = \psi+\chi,
\end{equation}
\begin{equation}
    \psi \chi^{.2} =\psi +\chi^. .
\end{equation}
In \cite{2} it was stated that an endomorphism $\phi$ on a manifold $\bar M$ is an almost golden structure , if
\begin{equation}
    \phi^2 X = \phi X + X,
\end{equation}
$X \in \Gamma (T \bar M)$. Hence, given a semi-Riemannian metric $\bar g$ on $\bar M$, $(\bar g, \phi)$ is referred to as an almost golden semi-Riemannian structure if

\begin{equation}
    \bar g (\phi X , Y)=\bar g ( X ,\phi Y),
\end{equation}
where for $X, Y \in \Gamma(T \bar M)$. Therefore, $(\bar M, \bar g, \phi)$ is called an almost golden semi-Riemannian manifold. In view of (8), we obtain \cite{2}
\begin{equation}
    \bar g = (\phi X , \phi Y)= \bar g ( X ,\phi Y)+ \bar g (X , Y).
\end{equation}

\begin{definition}{\rm  \cite{16}}
    Let $F$ is a $(1,1)$-tensor field on an almost golden manifold $(\bar M, \phi)$  which satisfies,
    $$\phi F^2 X =\phi X + FX$$
    for every $X \in \Gamma(T\bar M)$ Then, $F$ is called an almost meta-golden structure and $(\bar M,\phi, F)$ is called an almost meta-golden manifold.
\end{definition}

\begin{theorem} {\rm \cite{16}}
    A $(1,1)$-tensor field $F$ on an almost golden manifold $(\bar M, \phi)$ is an almost meta-golden structure if
    $$F^2 = \phi F-F+1.$$
\end{theorem}

\begin{definition} {\rm \cite{15}}
Let $F$ be an almost meta-golden structure on $(\bar M, \phi, \bar g)$ If $F$ is compatible with
semi-Riemannian metric $\bar g$ on $\bar M$, namely $$\bar g (FX , Y)= \bar g (X , FY),$$
or $$\bar g (FX , FY)= \bar g (\phi X , FY)-\bar g (X , FY)+ \bar g (X , Y),$$ then $(\bar M, \phi, F, \bar g)$  is known as an almost meta-golden semi-Riemannian manifold where for $X , Y \in \Gamma(T \bar M).$
\end{definition}

\begin{note}{\rm
An almost meta-golden semi-Riemannian manifold is called a meta-golden semi-Riemannian manifold if $\bar \nabla F=0$ here $\bar \nabla$ is Levi-Civita connection of $\bar M$. In this instance, we also have $\bar \nabla \phi=0$
From here throughout the paper, an almost meta-golden semi-Riemannian manifold (resp., meta-golden semi-Riemannian manifold) will be denoted as AMGsR manifold (resp., MGsR manifold).}
\end{note}

\subsection{Norden golden manifolds}\label{S2.5} 

The notion of almost Norden golden manifold can be recalled from \cite{2,30}. Let $\bar M$ be a manifold and $\phi$ an endomorphism on $\bar M$ such that $$\phi^2 =\phi - \frac{3}{2}I.$$ Then $(\bar M,\phi)$ is called an almost complex golden manifold. Let $\bar g$ be a semi-Riemannian metric on $\bar M$ such that 
\begin{equation} \label{ma}
    \bar g(\phi X ,Y)=\bar g(X,\phi Y),
\end{equation}then, $(\bar M,\phi, \bar g)$  is called an almost Norden golden manifold. Note that \eqref{ma} is comparable to $$\bar g(\phi X , \phi Y)= \bar g (\phi X , Y)-\frac{3}{2} \bar g (X ,Y).$$ Moreover, if $\phi$ is parallel with respect to a vector field $X$ on $\bar M$, $(\nabla_X \phi =0),$ then $(\bar M,\phi, \bar g)$ is called a locally decomposable almost Norden golden semi-Riemannian manifold (in short Norden golden semi-Riemannian manifold).

\subsection{Golden warped product Riemannian manifolds}\label{S2.6} 

Let $n$ and $m$ be the dimensions of two Riemannian manifolds $(\bar M_1,\bar g_1)$ and $(\bar M_2,\bar g_2)$, respectively. Indicate projection maps $P$ and $Q$ from the product manifold $\bar M_1 \times \bar M_2$ to $\bar M_1$ and $\bar M_2$, and the lift to $\bar M_1 \times \bar M_2$ of a smooth function $\varphi$ on $\bar M_1$ by $\bar{\varphi}_e := \varphi \circ P$.

In this context, we shall call $\bar M_1$ the base and $\bar M_2$ the fiber of $\bar M_1 \times \bar M_2$, the unique element $\bar{X}$ of $T(\bar M_1 \times \bar M_2)$ that is $P$-related to $X \in T(\bar M_1)$ and to the zero vector field on $\bar M_2$, the horizontal lift of $X$ and the unique element $\bar{Y}$ of $T(\bar M_1 \times \bar M_2)$ that is $Q$-related to $Y \in T(\bar  M_2)$ and to the zero vector field on $\bar M_1$, the vertical lift of $Y$ . Also denote by $\mathcal L(\bar {M}_1)$ the set of all horizontal lifts of vector fields on $\bar M_1$ and by $\mathcal L(\bar {M}_2)$ the set of all vertical lifts of vector fields on $\bar M_2$.

Let$f > 0$ be a smooth function on $\bar M_1$ and 

\begin{equation}\label{2.11}
\bar{g} := P^* \bar{g_1} + {(f \circ P)^2}{Q^*} \bar{g_2}
\end{equation}
be the Riemannian metric on  $\bar M_1 \times \bar M_2$ \cite{11}.
By \cite{12} the product manifold of $\bar M_1$ and $\bar M_2$ together with the Riemannian metric $\bar{g}$ defined by \eqref{2.11} is called a {\it warped product} of $\bar M_1$ and $\bar M_2$ by the warping function $f$ and it is denoted by $\bar{M}:=(\bar M_1 \times_f \bar M_2, \bar{g})$.

The golden warped product Riemannian manifold is defined by Blaga and Hretcanu in \cite{13} as:

\begin{theorem} {\rm \cite{11}} \label{wp}
Let $\left(\bar{M}=\bar M_1 \times_f \bar M_2, \bar{g}, \phi \right)$ be the warped product of two locally golden Riemannian manifolds $\left(\bar M_1,\bar g_1,\phi_1 \right)$ and $\left(\bar M_2,\bar g_2,\phi_1\right)$. Then $\bar{M}$ is locally golden if and only if
$$
\left\{\begin{array}{l}
\left(d f^2 \circ \phi_1\right) \otimes I=d f^2 \otimes \phi_2 \\
\bar g_2\left(\phi_1 \cdot, \cdot\right) \cdot \operatorname{grad}\left(f^2\right)=\bar g_2(\cdot, \cdot) \cdot \phi_1\left(\operatorname{grad}\left(f^2\right)\right).
\end{array} \right.
$$
\end{theorem}

\section{Submanifolds Immersed in Riemannian Manifolds with Golden Structure}\label{S3}

\subsection{Invariant submanifolds with golden structure}\label{S3.1}

Among all the submanifolds of an ambient manifold, invariant submanifolds are one of the common classes. It is commonly known that practically every property of the ambient manifold is inherited by an invariant submanifold.As a result, invariant submanifolds are a dynamic and productive area of study that has greatly influenced the advancement of modern differential geometry. Numerous articles concerning invariant submanifolds in golden Riemannian manifolds have been published.

In \cite{17}, Hretcanu and Crasmareanu studied invariant submanifolds in golden Riemannian manifolds and proved the following three propositions: 

\begin{proposition}\label{P13} Let $N$ be an $n$-dimensional submanifold of  a golden Riemannian manifold $(\bar{M}, \bar{g}, \phi)$ of codimension $r$ and let $(\phi ' , g, u_\alpha, \varepsilon \xi_\alpha,(a_{\alpha \beta})_r)$, $\alpha, \beta \in \{1,2,...,r\}$,  is induced structure on $N$ by structure $(\bar{g}, \phi)$ where $u_\alpha$ is $1$-form, $\xi_\alpha$ tangent vector fields and $\left (a_{\alpha \beta}\right)_r$ is $r \times r$ matrix of real function on $N$. Then an essential and sufficient requirement for $N$ to be invariant is that the induced structure $(\phi ', g)$ on $N$ is a golden Riemannian structure, whenever $\phi '$ is non-trivial.
\end{proposition}

\begin{proposition} \label{P14} 
    If $(N,g,\phi')$ is an invariant submanifold of codimension $r$ in a golden Riemannian manifold $(\bar{M}, \bar{g}, \phi)$ with $\bar{\nabla} \phi=0$ and $(\phi ', g, u_\alpha, \xi_\alpha,(a_{\alpha \beta})_r)$ is the induced structure on $N$ by $(\phi, \bar{g})$ (here $\bar\nabla$ is the Levi-Civita connection defined on $N$ with respect to $g$) then the Nijenhuis torsion tensor field of $\phi '$ vanishes identically on $N$.
\end{proposition}

\begin{proposition} \label{P15} 
 Let $N$ is an invariant submanifold of codimension $r$ in a golden Riemannian manifold $(\bar{M}, \bar{g}, \phi)$ with $\bar{\nabla} \phi=0$ and $\left(\phi ', g, u_\alpha, \xi_\alpha,\left(a_{\alpha \beta}\right)_r\right)$ be the induced structure on $N$. If normal connection $\nabla^{\perp}$ on  normal bundle $TN^{\perp}$ vanishes identically $\left(l_{\alpha \beta}=0\right)$, then the components $\mathcal{N}^{(1)}, \mathcal{N}^{(2)}, \mathcal{N}^{(3)}$ and $\mathcal{N}^{(4)}$ of the Nijenhuis torsion tensor field of $\phi '$ for the structure $\left(\phi ', g, \xi_\alpha, u_\alpha,\left(a_{\alpha \beta}\right)_r\right)$ induced on $N$ have the forms:
 \vskip.05in
    
$\left\{\begin{array}{l}(\text {\rm i})\; \mathcal{N}^{(1)}(X, Y)=\mathcal{N}_{\alpha \beta}^{(4)}(X)=0, \\ (\text {\rm ii})\;\mathcal{N}_\alpha^{(2)}(X, Y)=-\sum_\beta a_{\alpha \beta} g\left(\left(\phi ' A_\beta-A_\beta \phi '\right)(X), Y\right) 
\\ (\text {\rm iii})\; \mathcal{N}_\alpha^{(3)}(X)=\sum_\beta a_{\alpha \beta}\left(\phi ' A_\beta-A_\beta \phi ' \right)(X)-\phi ' \left(\phi ' A_\alpha-A_\alpha \phi '\right)(X),\end{array}\right.$ 

 \vskip.05in
 \noindent for any $X, Y \in \Gamma(TN)$.
\end{proposition}

\begin{remark} {\rm \cite{17}
Under the conditions of Proposition \ref{P15}, if $\phi ' A_\alpha=A_\alpha \phi '$ where $A$ is the  shape operator for $\alpha \in\{1,2, \ldots, r\}$, then the components $\mathcal{N}^{(1)}, \mathcal{N}^{(2)}, \mathcal{N}^{(3)}$ and $\mathcal{N}^{(4)}$ vanishes identically on $N$.}
\end{remark}

Inspired from \cite{17}, G\"ok et al. in \cite{18} have demonstrated the local decomposability of any invariant submanifold of a golden Riemannian manifold and came up with a definition of invariance for  submanifolds in a golden Riemannian manifold. They also determined the prerequisites that must be met for any invariant submanifold to be totally geodesic.

\begin{theorem}\label{T17} {\rm \cite{18}} Let $N$ be an invariant submanifold of a locally decomposable golden Riemannian manifold $(\bar M,\bar g,\phi)$.  Then $N$ is a locally decomposable golden Riemannian
manifold whenever the induced structure $\phi '$ on $N$ is non-trivial.
\end{theorem}

\begin{theorem}\label{T18} {\rm  \cite{18}} Let $N$ be an invariant submanifold of a golden Riemannian manifold $(\bar M,\bar g,\phi)$. Then $N$ is an invariant submanifold if and only if there exists a local orthonormal frame of the normal bundle $TN^{\perp}$ such that it consists of eigenvectors of the golden structure $\phi.$
\end{theorem}

\begin{remark}{\rm 
The result for $N$ to be totally geodesic invariant submanifold was also obtained  in \cite{18} by G\"ok et al.}     
\end{remark}

\subsection{Anti-invariant submanifolds in golden Riemannian manifolds}\label{S3.2}

Some properties of anti-invariant submanifold of a golden Riemannian manifold have been studied in  \cite{22} and some necessary prerequisite for any submanifold in a locally decomposable golden Riemannian manifold to be anti-invariant are obtained. Any anti-invariant submanifold $N$ of a golden Riemannian manifold $(\bar M,\bar g,\phi)$
is a submanifold such that the golden structure $\phi$ of the ambient manifold $\bar M$ carries each tangent vector of the submanifold $N$ into its corresponding normal space in the ambient manifold $\bar M$, that is, $$\phi(T_x N) \subseteq T_x N^\perp$$ for each point $x \in N$ \cite{22}.

\begin{theorem} \label{T20} {\rm  \cite{22}}
   Let $N$ be an n-dimensional submanifold of  a $2n$-dimensional locally decomposable golden Riemannian manifold $(\bar M,\bar g,\phi)$. Then, for any $\alpha, \beta \in \{1,2,...n\}$, $N$ is an anti-invariant submanifold whenever $a_{\alpha \beta} =\delta_{\alpha \beta}$. In addition, the submanifold $N$ is totally geodesic.
\end{theorem}

\begin{theorem}\label{T21}  {\rm \cite{21}}
   Let $N$ be an $n$-dimensional submanifold of a $2n$-dimensional locally decomposable golden Riemannian manifold $(\bar M,\bar g,\phi)$. If $\phi^-1 (\eta_\alpha) \in \Gamma(TN)$ for any $\alpha \in \{1,2,...n\}$, then $N$ is an anti-invariant submanifold. Furthermore, the submanifold $N$ is totally geodesic. 
\end{theorem}

\begin{remark}{\rm
  M. Gök et al. have also obtained results in \cite{22} for the existence of orthonormal frame of  anti-invariant submanifold of  locally decomposable golden Riemannian manifold.}
\end{remark}

Study of non-invariant submanifold of a locally decomposable golden Riemannian manifold in the case that the rank of the set of tangent vector fields of the induced structure on the submanifold by the golden structure of the ambient manifold is less than or equal to the co-dimension of the
submanifold was done in \cite{23} by Gök and Kılıç.

\begin{theorem} \label{T23} {\rm \cite{23}}
    Let $N$ be a submanifold of codimension $r$ in a locally decomposable golden Riemannian manifold $(\bar M,\bar g,\phi)$. If the tangent vector fields $\xi_1 ,...,\xi_r$ are linearly
independent and $\nabla \phi ' = 0$, then $N$ is a totally geodesic submanifold.
\end{theorem}

\begin{theorem} \label{T24} {\rm \cite{23}}
Let $N$ be a submanifold of codimension $r$ in a locally decomposable golden Riemannian manifold $(\bar M,\bar g,\phi)$ and where $\lambda_\alpha$'s are the eigenvalues of the matrix $(a_{\alpha \beta})_{r \times r}$. If $a_{\alpha \beta} =\lambda_\alpha \delta_{\alpha \beta}, \lambda_\alpha \in (1-\psi, \psi)$ for any $\alpha,\beta \in \{1,...,r\},$ and $\nabla \phi '=0,$ then $N$ is totally geodesic.
\end{theorem}

\begin{theorem} \label{T25} {\rm \cite{23}}
Let $N$ be a submanifold of codimension $r$ in a locally decomposable golden Riemannian manifold $(\bar M,\bar g,\phi)$.  If $a_{\alpha \beta} =\lambda_\alpha \delta_{\alpha \beta}, \lambda_\alpha \in (1-\psi, \psi)$ for any $\alpha,\beta \in \{1,...,r\},\, {\rm Trace}\, \phi '$ is constant, and $N$ is totally umbilical, then $N$ is totally geodesic.    
\end{theorem}

\begin{remark} {\rm
{\rm (i)} Gök and Kılıç \cite{23} also obtained some results on the non-invariant submanifold in the event that the tangent vector fields of the induced structure are linearly dependent.
\smallskip

\noindent {\rm (ii)} Stability problem of certain anti-invariant submanifolds in golden Riemannian manifolds have been discussed by the same authors in \cite{24}.

\smallskip
\noindent {\rm (iii)} Effective relations for certain induced structures on submanifold  of codimension 2 in golden Riemannian manifolds were obtained in \cite{28}.}
\end{remark}

\subsection{Slant submanifolds of golden Riemannian manifolds}\label{S3.3}

Given a golden Riemannian manifold $(\bar M,\bar g,\phi)$, let $(N,g)$ be a submanifold of it, and $g$ be the induced metric on $N$. Then, we can write $$\phi X=PX+QX$$ for any $X \in \Gamma(TN)$, where $PX$ and $QX$ are the tangent and transversal components of $\phi X$.

A submanifold $(N,g)$ of a golden Riemannian manifold $(\bar M,\bar g,\phi)$ is referred to as a slant submanifold if, at any point $x$ each nonzero vector $X$ tangent to $N$, the angle between $TN$ and $\phi (X)$, represented by $\theta (X)$, and $\theta (X)$, is independent of the selection of $x \in N$ and $X \in T_x N$. It can also be seen that $N$ is a $\phi$-invariant (resp. $\phi$-anti-invariant) submanifold if the slant angle $\theta = 0$ (resp. $\theta = \frac{\pi}{2}$). The term proper slant (or $\theta$-slant proper) submanifold refers to a slant submanifold that is neither anti-invariant nor invariant.

Using a golden Riemannian manifold, Uddin and Bahadir \cite{25} establish the following characterization of slant submanifolds.

\begin{theorem} \label{T27} {\rm \cite{25}}
   Assume that a golden Riemannian manifold $(\bar{M},\bar g, \phi)$ has a submanifold $(N, g)$. Consequently, $N$ is a slant submanifold if and only if $c \in[0,1]$ is a constant such that
$$
P^2=c (\phi+I),
$$
Additionally, if $\theta$ represents the slant angle of $N$, then $c=\cos ^2 \theta$.
\end{theorem}

\begin{corollary} {\rm \cite{25}}
   Take a golden Riemannian manifold $(\bar{M},\bar g, \phi)$ and let $(N, g)$ be a submanifold of it. After that, $N$ is a slant submanifold if and only if $c \in [0,1]$ exists and ensures
$$
\phi^2=\frac{1}{c} P^2,
$$
where $c=\cos ^2 \theta$ and $\theta$ are slant angles of $N$.
\end{corollary}

\begin{remark} {\rm
    In \cite{25} Uddin and Bahadir have also derived some results of $\phi$-invariant and $\phi$-anti-invariant submanifolds  of a golden Riemannian manifold and provided some examples of such submanifolds. }
\end{remark}

\subsection{Semi-invariant submanifolds of golden Riemannian manifolds}\label{S3.4}

\begin{definition} {\rm \cite{19}}
    Given a golden Riemannian manifold $(\bar M,\bar g,\phi)$, consider $N$ to be a real submanifold of $\bar M$. If $N$ is equipped with a pair of orthogonal distributions $(D,D^{\perp})$ that meet the given conditions, it can be deemed as a semi-invariant submanifold of $\bar M$.    
\vskip.05in
{\rm (i)} $TN= D \oplus D^\perp,$

\vskip.05in
{\rm (ii)} the distribution $D$ is invariant, i.e.,
$\phi D_x =D_x$ for each $x \in N$

\vskip.05in
{\rm (iii)}  the distribution $D^\perp$ is anti-invariant, i.e.,
$\phi D^\perp \subset T_x N^\perp$  for each $x \in N.$

For any $x \in N$, a semi-invariant submanifold $N$ is considered invariant and anti-invariant if $D^{\perp}_{x} = 0$ and $D_x = 0$, respectively.
\end{definition}

Following results have been obtained by Erdogan et al. for semi-invariant submanifolds of the golden Riemannian manifold in \cite{19}.

\begin{theorem} {\rm \cite{19}}
    Assume that $N$ is a semi-invariant submanifold of $(\bar M,\bar g,\phi)$, the golden Riemannian manifold. Consequently, $D$ distribution is integrable if and only if
    $$h(X, \phi ' Y) = h(Y, \phi ' X)$$ here $h$ is the second fundamental form and $X,Y \in \Gamma(D)$ and $Z \in \Gamma(D^\perp)$.
\end{theorem}

\begin{theorem} {\rm \cite{19}}
     Let $N$ be a semi-invariant submanifold of the golden Riemannian manifold $(\bar M,\bar g,\phi)$. Then,  distribution $D$ is integrable if and only if $$\phi ' A_{\phi ' X} Y=A_{\phi ' X} Y$$ has no components in $D$, where $A$ is shape operator and for every $X,Y \in \Gamma(D^\perp)$ and $Z \in \Gamma(D)$.
\end{theorem}

\begin{remark} {\rm
   The conditions for distributions $D$ and $D^\perp$, respectively, of semi-invariant submanifolds of golden Riemannian manifolds to be a totally geodesic foliation have been examined in \cite{19}. The condition for a semi-invariant submanifold $N$ to be totally geodesic is also covered. Also, one can go to \cite{25} for conclusions of a similar kind for semi-invariant submanifolds of the golden Riemannian manifold.}
\end{remark}

\begin{remark} {\rm
   M. Gök et al. in\cite{20} proposed specific characterizations for every submanifold of a golden Riemannian manifold to be semi-invariant in terms of canonical structures on the submanifold, as a consequence of the ambient manifold's golden structure.}
\end{remark}

Totally umbilical semi-invariant 
submanifolds of golden Riemannian manifolds have been studied in \cite{21}.
\begin{theorem} {\rm \cite{21}}
    Let $N$ be a totally umbilical submanifold of a golden Riemannian manifold $\bar{M}$. Then, the distribution $D$ is always integrable.
\end{theorem}

\begin{theorem} {\rm \cite{21}}
Let $N$ be a totally umbilical submanifold of a golden Riemannian manifold $(\bar{M}, \phi)$. Then $D^{\perp}$ is integrable. 
\end{theorem}

\begin{remark} {\rm
    Moreover the properties of semi-invariant submanifolds and  totally umbilical semi-invariant submanifold of golden Riemannian Manifolds having constant sectional curvatures are studied by Sahin et al. in \cite{27}.}
    \end{remark}

\subsection{Skew semi-invariant submanifolds}\label{S3.5}

In \cite{32} Ahmad and  Qayyoom have studied skew semi-invariant submanifolds in golden Riemannian manifold and in locally golden Riemannian manifold:

\begin{definition} {\rm \cite{32}}
    A submanifold $N$ of a golden Riemannian manifold $\bar{M}$ is defined as skew semi-invariant submanifold if there exist an integer $k$ and constant functions $\alpha_i, 1 \leq i \leq k$, defined on $N$ with values in $(0,1)$ such that
    
{\rm (i)} each $\alpha_i, 1 \leq i \leq k$, is a distinct eigenvalue of $\phi^2$ with
$$
T_x N=D_x^0 \oplus D_x^1 \oplus D_x^{\alpha_1} \oplus \cdots \oplus D_x^{\alpha_k}
$$
for $x \in N$, and

{\rm (ii)} the dimensions of $D_x^0, D_x^1$ and $D_k^1, 1 \leq i \leq k$, are independent of $x \in N$.
\end{definition}

\begin{remark} {\rm
   The tangent bundle of $N$ has the following decomposition $$TN=D^0 \oplus D^1 \oplus  D^{\alpha_1} \oplus \cdots \oplus D^{\alpha_k}.$$ If $k=0,$  then $N$ is a semi-invariant submanifold. Also, if $k=0$ and $D_x^0$ (resp., $D_x^1$) is trivial, then $N$ is an invariant (resp., anti-invariant) submanifold of $\bar M.$}
\end{remark}

\begin{definition} {\rm  \cite{32}}
    A submanifold $N$ of a locally golden Riemannian manifold $\bar{M}$ is defined as skew semi-invariant submanifold of order 1 if $N$ is a skew semi-invariant submanifold with $k = 1$.
    In this case, we have
$$
T N=D^{\perp} \oplus D^T \oplus D^\theta,
$$
where $D^\theta=D^{\alpha_1}$ and $\alpha_1$ are constant. A skew semi-invariant submanifold of order 1 is proper if $D^{\perp} \neq 0$ and $D^T \neq 0$.
\end{definition}

\begin{remark} {\rm
    Some lemma's for  proper skew semi-invariant submanifolds of a  locally golden Riemannian manifold have been also discussed in \cite{32}.}
\end{remark}

\subsection{Pointwise slant submanifolds in golden Riemannian manifolds}\label{S3.6}

The notion of slant submanifolds in almost Hermitian manifolds was first introduced by the first authors in \cite{c90,c90.1}. Later, the first author and Garay \cite{CG} extended the notion of slant submanifolds to pointwise slant submanifolds in almost Hermitian manifolds.  Hretcanu  and Blaga \cite{13} defined the notion of pointwise slant submanifolds of golden Riemannian manifolds as follows:

A submanifold $N$ of a golden Riemannian manifold $(\bar M,\bar g, \phi)$  is referred to as {\it  pointwise slant}  \cite{13} if, at every point $x \in N$, the angle $\theta_x (X)$ between $\phi_x$ and $T_x N$  (called the Wirtinger angle) is consistent regardless of nonzero tangent vector $X \in T_x N \setminus \{0\}$, but it depends on $x \in N$. The Wirtinger angle is a real-valued function $\theta$ (called Wirtinger function), verifying $$\cos \theta_x=\frac{\bar{g}(\phi X, T X)}{\|\phi X\| \cdot\|T X\|}=\frac{\|T X\|}{\|\phi X\|},$$
for any $x \in N$ and $X \in T_x N \backslash\{0\}$.
If the Wirtinger function $\theta$ of a pointwise slant submanifold of a golden Riemannian manifold is globally constant, it is referred to as {\it slant submanifold}.

\begin{proposition} {\rm \cite{13}}
   In the golden Riemannian manifold $(\bar{M}, \bar{g}, \phi)$, if $N$ is an isometrically immersed submanifold and $T$ is the map, then $N$ is a pointwise slant submanifold if and only if
$$
T^2=\left(\cos ^2 \theta_x\right)(T+I),
$$
for some real-valued function $x \mapsto \theta_x$, for $x \in N$.
\end{proposition}

\begin{proposition} {\rm \cite{13}}    Let $N$ be a submanifold of the golden Riemannian manifold $(\bar{M}, \bar{g}, \phi)$ that is isometrically immersed. Given $N$ as a pointwise slant submanifold and $\theta_x$ as its Wirtinger angle, then
$$
\left(\nabla_X T^2\right) Y=\left(\cos ^2 \theta_x\right)\left(\nabla_X T\right) Y-\sin \left(2 \theta_x\right) X\left(\theta_x\right)(T Y+Y),
$$
for any $X, Y \in T_x N \backslash\{0\}$ and any $x \in N$.
\end{proposition}

\subsection{Pointwise bi-slant submanifolds in golden Riemannian manifolds}\label{S3.7}

Consider $N$ as an immersed submanifold within the golden Riemannian manifold $(\bar{M}, \bar{g}, \phi)$. We define $N$ as a {\it pointwise bi-slant submanifold} of $\bar{M}$ if there exist two orthogonal distributions $D$ and $D^\perp$ on $N$ such that:

(i) $T N = D \oplus D^\perp$; 

(ii) $\phi(D) \perp D^\perp$ and $\phi(D^\perp) \perp D$; 

(iii) The distributions $D$ and $D^\perp$ are pointwise slant, with slant functions $\theta_{1 x}$ and $\theta_{2 x}$, for $x \in N$. The pair $\{\theta_1, \theta_2\}$ of slant functions is referred to as the bi-slant function.

\vskip.05in
A pointwise bi-slant submanifold $N$ is called {\it proper} if its bi-slant functions $\theta_1, \theta_2 \neq 0 ; \frac{\pi}{2}$ and both $\theta_1$ and $\theta_2$ are not constant on $N$. Specifically, if $\theta_1 = 0$ and $\theta_2 \neq 0 ; \frac{\pi}{2}$, then $N$ is called a pointwise semi-slant submanifold; if $\theta_1 = \frac{\pi}{2}$ and $\theta_2 \neq 0 ; \frac{\pi}{2}$, then $N$ is called a pointwise hemi-slant submanifold.
$T\left(D\right) \subseteq D$ and $T\left(D^\perp \right) \subseteq D^\perp$ are verified by the distributions $D$ and $D^\perp$ on $N$, if $N$ is a pointwise bi-slant submanifold of $\bar{M}$ \cite{13}.

\begin{remark} {\rm
Some examples of pointwise bi-slant submanifolds in golden Riemannian manifolds were given in \cite{13}, where Blaga and Hretcanu gave fundamental lemmas for pointwise bi-slant, pointwise semi-slant, and pointwise hemi-slant submanifolds in locally golden Riemannian manifolds.}  
\end{remark}

\subsection{CR-submanifolds of a golden Riemannian manifold}\label{S3.8}

 A submanifold $N$ within a golden Riemannian manifold $\bar{M}$ is called a CR-submanifold if there exists a differentiable distribution $D: X \to D_x \subseteq T_x N$ on $N$ that meets the following criteria:

(i) $D$ is holomorphic, meaning $\phi D_x = D_x$ for every $x \in N$, and 

(ii) The orthogonal complementary distribution $D^\perp : x \to D^\perp \subseteq T_x N$ is completely real, i.e., $\phi D^\perp \subset T_x N^\perp$ for each $x \in N$. If $\dim D^\perp_x = 0$ (or, $\dim D_x = 0$), then the CR-submanifold $N$ is a holomorphic submanifold (or, a totally real submanifold, respectively). If $\dim D^\perp_x = \dim T_x N^\perp$, then the CR-submanifold is an anti-holomorphic submanifold (or a generic submanifold). A submanifold is considered as proper CR-submanifold if it is neither holomorphic nor totally real \cite{29}.

 In \cite{29} authors have defined and studied CR-submanifolds of a golden Riemannian manifold.

\begin{proposition} {\rm \cite{29}} Let $N$ be a CR-submanifold of a locally golden Riemannian manifold $\bar M$. Then
 $$\bar g(\phi A_{\phi Y} X,Z)+\bar g(\nabla_X Y,\phi Z)+\bar g(\nabla_X Y,Z)=0,$$  $$A_{\phi \xi '} Z=-A-\xi ' \phi Z,$$  $$A_{\phi Y} W=A_{\phi W}Y$$ 
 for $X \in TN, Z \in D, Y,W \in D^\perp$ and $\xi ' \in \mathcal{V}$, where $\mathcal{V}$ the complementary orthogonal subbundle of $\phi (D^\perp)$ in $TN^\perp$,  $\xi '$ is the unit normal vector field, and $X, Y, Z, W$ are vector fields
\end{proposition}

\begin{lemma} {\rm \cite{29}}
   Consider $N$ as a CR-submanifold of $\bar M$, a locally golden Riemannian manifold. Given any $Y,W \in D^\perp$, then
 $$(\nabla^{\perp} _{W} \phi Y- \nabla^{\perp} _{Y} \phi W) \in \phi D^\perp.$$
\end{lemma}

\begin{remark} {\rm
    The  integral condition of $D$ of CR-submanifolds of golden Riemannian manifold is also discussed in \cite{29}.}
\end{remark}

The following outcomes were attained by Ahmad and Qayyoom in \cite{29} from their study of totally umbilical CR-submanifolds of golden Riemannian manifolds:

\begin{lemma} {\rm \cite{29}}
  Assume that $N$ is a totally umbilical CR-submanifold of $\bar M$, a locally golden Riemannian manifold. Then, either $H$, the mean curvature vector, is perpendicular to $\phi (D^\perp)$, or the totally real distribution $D^\perp$ is one-dimensional.
\end{lemma}

\begin{theorem} {\rm \cite{29}}
For a locally golden Riemannian manifold $\bar M$, let $N$ be a totally umbilical CR-submanifold. Following, $\bar K (\pi)=0$ 0 for every CR-section $\pi$, i.e., the CR-sectional curvature of $\bar M$ vanish.
\end{theorem}

\section{Warped Product Manifolds in Golden Riemannian Manifolds}\label{S4}

The golden warped product Riemannian manifold were defined by Blaga et al. in \cite{11} as mentioned in Theorem \eqref{wp} . In \cite{11} authors have also studied its properties with a special view towards its curvature and have attained the following outcomes:

\begin{theorem} {\rm \cite{11}}
    Let $\left(\bar{M}=\bar M_1 \times_f \bar M_2, \bar{g}, \phi \right)$ $($with $\bar{g}$ given by Equation \eqref{2.11}$)$ be the warped product of the golden Riemannian manifolds $\left(\bar M_1,\bar g_1, \phi_1\right)$ and $\left(\bar M_2,\bar g_2, \phi_2\right)$. If $\bar M_1$ and $\bar M_2$ have $\phi_1$ and $\phi_2$-invariant Ricci tensors, respectively (i.e. $Q_{\bar M_i} \circ \phi_i=\phi_i \circ Q_{\bar M_i}, i \in\{1,2\}$ ), then $\bar{M}$ has $\phi$-invariant Ricci tensor if and only if we have
$$
\operatorname{Hess}(f)\left(\phi_1 \cdot, \cdot\right)-\operatorname{Hess}(f)\left(\cdot, \phi_1 \cdot\right) \in\{0\} \times T\left(\bar M_2\right) .
$$
\end{theorem}

\begin{remark} {\rm
    In \cite{11} authors have also provided  example of golden warped product Riemannian manifold.}
\end{remark}

In \cite{13} Blaga and Hretcanu have studied warped product pointwise bi-slant submanifolds and warped product pointwise semi-slant or hemi-slant submanifolds in golden Riemannian manifolds and have derived the following results:

\begin{definition} {\rm \cite{13}}
The warped product $N_1 \times_{f} N_2$ of two pointwise slant submanifolds $N_1$ and $N_2$ within a golden Riemannian manifold $(\bar M,\bar g, \phi)$ is referred to as a warped product pointwise bi-slant submanifold. Furthermore, the pointwise bi-slant submanifold $N_1 \times_{f} N_2$ is termed proper if both submanifolds $N_1$ and $N_2$ are proper pointwise slant in $(\bar M,\bar g, \phi)$.
\end{definition}

\begin{definition}\hskip-.03in  {\rm \cite{13}}
Consider $N := N_1  \times_{f}  N_2$ as a warped product bi-slant submanifold within a golden Riemannian manifold $(\bar M,\bar g, \phi)$, where one of the components $N_i$ $ (i \in \{1, 2\})$ is either an invariant submanifold or an anti-invariant submanifold in $\bar M$, and the other component is a pointwise slant submanifold in $\bar M$ with the Wirtinger angle $\theta_x \in [0, \frac{\pi}{2}]$. In this context, the submanifold $N$ is referred to as a warped product pointwise semi-slant submanifold or a warped product pointwise hemi-slant submanifold in the golden Riemannian manifold $(\bar M,\bar g,\phi)$.
\end{definition}

\begin{theorem} {\rm \cite{13}}
    Consider $N:=N_T \times N_\theta$ as warped product pointwise semi-slant submanifold within a locally golden Riemannian manifold $(\bar{M}, \bar{g}, \phi)$ with the pointwise slant angle $\theta_x \in\left(0, \frac{\pi}{2}\right)$, for $x \in N_\theta$, then the warping function $f$ is constant on the connected components of $N_T$.
\end{theorem}

\begin{remark} {\rm
 Blaga and Hretcanu \cite{HB20,13} gave examples of warped product pointwise bi-slant submanifolds, warped product semi-slant submanifolds, and warped product hemi-slant submanifolds within a golden Riemannian manifold. Additionally, they discussed various results related to warped product pointwise semi-slant submanifolds and warped product pointwise hemi-slant submanifolds.}
\end{remark}

{\rm In \cite{32}, Ahmad and Qayyoom introduced and examined warped product skew semi-invariant submanifolds within a locally golden Riemannian manifold. They also explored the conditions required and sufficient for a skew semi-invariant submanifold in such a manifold to be classified as a locally warped product.}

\begin{proposition} {\rm \cite{32}}
For a locally golden Riemannian manifold $\bar{M}$, let $N=N_1 \times_f N_T$ be a $\left(D^\theta, D^T\right)$-mixed totally geodesic proper skew semi-invariant submanifold with integrable distribution $D^{\perp}$. Then $N$ is a locally warped product submanifold if
$$
A_{\phi V} \phi X=-V(\ln f)(T X-X),
$$
$$\bar g\left(A_{\eta T Z} X, Y\right)+\bar g\left(A_{\eta Z} X, \phi Z\right)=\sin ^2 \theta X(\ln f)[\bar g(Y, Z)+\bar g(Y, T Z)].$$ for any $\eta , V \in TN^\perp .$
\end{proposition} 

\begin{lemma} {\rm \cite{32}}  Assume that $N=N_1 \times_f N_T$ is a warped product proper skew semi-invariant submanifold of a locally golden Riemannian manifold. Then we have
$$
\bar g(h(X, V), \phi W)=0
$$
$$
\bar g(h(X, V), N Z)=0.
$$
\end{lemma}

\begin{lemma} {\rm \cite{32}}
    Let $N=N_1 \times_f N_T$ be a warped product proper skew semi-invariant submanifold $N$ of a locally golden Riemannian manifold. Then
$$
\bar g(h(X, \phi Y), \phi V)=-V(\ln f)[\bar g(\phi Y, X)+\bar g(Y, X)] .
$$
\end{lemma}

\begin{theorem} {\rm \cite{32}}
    Consider $N=N_1 \times_f N_T$ as a $(p+q+r)$-dimensional warped product proper skew semi-invariant submanifold within a $(2p+2q+r)$-dimensional locally golden Riemannian manifold $\bar{M}$. The following statements hold true:
    
{\rm (i)} The squared norm of the second fundamental form of $N$ meets the condition
$$\|h\|^2 \geq r\left\{2\left\|\nabla^{\perp}(\ln f)\right\|^2+2 \cos ^2 \theta\left\|\nabla^\theta(\ln f)\right\|^2\right\},$$
where $r=\operatorname{dim}\left(N_T\right)$, and $\nabla^{\perp}(\ln f)$ and $\nabla^\theta(\ln f)$ are gradients of $(\ln f)$ on $D^{\perp}$ and $D^\theta$, respectively.

{\rm (ii)} Assume that the equality sign remains unchanged. It follows that $N$ is a mixed entirely geodesic and $N_1$ a totally geodesic submanifold of $\bar{M}$. Furthermore, $N_T$ is never going to be $\bar{M}$'s minimal submanifold.
\end{theorem}

\section{Lightlike Submanifolds of Golden Semi-Riemannian Manifolds}\label{S5}

Research in the geometry of degenerate submanifolds shows significant differences compared to non-degenerate submanifolds. This variation stems from the fact that the tangent bundle of non-degenerate submanifolds intersects trivially with the normal vector bundle, whereas this intersection is non-trivial in degenerate submanifolds. The work of Duggal, Bejancu and Kupeli on lightlike submanifolds within semi-Riemannian manifolds is documented in \cite{33} and \cite{34}, respectively. Furthermore, the study of lightlike submanifolds in semi-Riemannian and specifically golden semi-Riemannian manifolds is a critical field in differential geometry, attracting numerous researchers.

For a lightlike submanifold $N$ of a semi-Riemannian manifold $(\bar M,\bar g)$,  Duggal and  Bejancu \cite{33} defined   the notion of
{\it radical distribution} $Rad(TN)$ and the notion of {\it normal bundle} $TN^{\perp}$ such that
\begin{align*} &Rad(TN) = TN \cap TN^{\perp},\end{align*} where
$T^{\perp}N = \cup_{x \in N} \big\{ X \in T_x\overline{M} \mid \bar{g}(X,Y) = 0, \forall Y\in  T_x\,N\big\}.$

\subsection{Lightlike Hypersurfaces}\label{S5.1}

{\rm Poyraz and Yaşar have introduced lightlike hypersurfaces of a golden semi-Riemannian manifold in \cite{5}.
Let $N$ be a lightlike hypersurface of a golden semi-Riemannian manifold $(\bar M, \bar g, \phi)$. For any $X \in \Gamma(TN)$ and $\mathcal{T} \in {\rm Trace}(TN)$, can be decompose $$\phi X=\phi ' X+u_1(X)\mathcal{T},\quad \phi \mathcal{T}=U_1 +u_2(\xi ') \mathcal{T},$$ where $\phi ' X \in \Gamma(TN)$ and $u_1$ and $u_2$ are 1-form defined by $$u_1 (X)=g(X, \phi \xi '),\quad u_2(X)=g(X,\phi \mathcal{T}).$$

\begin{definition} {\rm \cite{5}}  Let $N$ be a lightlike hypersurface of a golden semi-Riemannian manifold $(\bar{M}, \bar{g}, \phi)$. 
Then
    
{\rm (i)} $N$ is  called  a screen semi-invariant lightlike hypersurface if $\phi ({Rad}(T N)) \subset S(T N)$. Here, a screen distribution on $N$, denoted by $S(TN)$, is defined as a non-degenerate complementary vector bundle of $TN^\perp$ in $TN$. Additionally, $\phi (\operatorname{ltr}(T N)) \subset S(T N)$ where $\operatorname{ltr}(T N)$ denotes the lightlike transversal bundle associated with the hypersurface N;

{\rm (ii)} $N$ is known as radical anti-invariant lightlike hypersurface if $\phi ({Rad}(T N)) \subset \operatorname{ltr}(T N)$.
\end{definition}

\begin{theorem} {\rm  \cite{5}} Let $N$ be a lightlike hypersurface of a golden semi-Riemannian manifold $(\bar{M}, \bar{g}, \phi)$ and consider the induced structure $(g, \phi ')$ on $TN$. Then the next three statements are equivalent.

   {\rm (i)} $N$ is invariant.

{\rm (ii)} $u_1$ vanishes on $N$.

{\rm (iii)} $\phi '$ is golden structure on $N$.
\end{theorem}

\begin{theorem} {\rm \cite{5}}
An anti-invariant lightlike hypersurface of a golden semi-Riemann manifold is not radical.
\end{theorem}

In \cite{5}, Poyraz and Yaşar have derived certain findings regarding screen semi-invariant lightlike hypersurfaces in a golden semi-Riemannian manifold.

\begin{theorem} {\rm \cite{5}}
    A golden semi-Riemannian manifold $(\bar M,\bar g,\phi)$ has a screen semi-invariant lightlike hypersurface denoted by $(N, g, S(TN))$. When $N$ is totally geodesic in $\bar M$ and $u = 0$, where $u$ is 1-form, then only the vector field $U =\phi \xi'$ on $N$ is parallel to $\nabla$.
\end{theorem}

Both mixed geodesic lightlike hypersurfaces and totally geodesic lightlike hypersurfaces are introduced in \cite{5}.

\begin{theorem} {\rm \cite{5}}
    Let $(N, g, S(TN))$ be a screen semi-invariant lightlike hypersurface of a golden semi-Riemannian manifold $(\bar M,\bar g,\phi)$. Then, the following claims are equivalent:

   {\rm (i)} N is mixed geodesic.

    {\rm (ii)} There is no $D_2$-component of $A_\mathcal{T}$.

   {\rm  (iii)} There is no $D_1$-component of $A^*_{\xi '}$.
\end{theorem}

\begin{theorem} {\rm \cite{5}}
     Let $(N, g, S(TN))$ be a screen semi-invariant lightlike hypersurface of a golden semi-Riemannian manifold $(\bar M,\bar g,\phi)$. It follows that $N$ is totally geodesic if and only if, for any $X \in \Gamma (TN)$ and $Y \in \Gamma (D)$, we have $$(\nabla_{X} \phi ')Y=0,\quad (\nabla_{X} \phi ')U_1=A_\mathcal{T} X.$$
\end{theorem}

\begin{theorem} {\rm \cite{5}}
   Assume that the locally golden product space form is $\bar M = M_p (c_p) \times M_q (c_q)$, and that the screen semi-invariant lightlike hypersurface of $\bar M$ is $N$. If $N$ is totally umbilical, then $c_p = - (\psi + 1)c_q$.
\end{theorem}

\begin{remark} {\rm
 Poyraz and Yaşar have presented additional findings on the screen semi-invariant lightlike hypersurface in the locally golden product space form as described in \cite{5}.}
\end{remark}

In \cite{5} Poyraz and Yaşar have also studied  screen conformal screen semi-invariant lightlike hypersurfaces and have obtained the following outcomes:

\begin{theorem} {\rm \cite{5}}
    Assume that the golden semi-Riemannian manifold $(\bar M,\bar g,\phi)$ has a screen conformal screen semi-invariant lightlike hypersurface $(N, g, S(TN))$. The leaf $N^b$ of $S(TN)$ is totally geodesic in both $N$ and $\bar M$, and $N$ is totally geodesic in $\bar M$ if $N$ or $S(TN)$ is totally umbilical.
\end{theorem}

\begin{theorem} {\rm \cite{5}}
    Consider $(N, g, S(TN))$ be a screen conformal screen semi-invariant lightlike hypersurface within a locally golden product space form $\bar M = M_p (c_p) \times M_q (c_q)$  Then, we have $c_p = c_q = 0$.
\end{theorem}

\begin{corollary}{\rm \cite{5}}
There exists no screen conformal screen semi-invariant lightlike hypersurface within a locally golden product space form $\bar M = M_p (c_p) \times M_q (c_q)$ with $c_p ,c_q \neq 0.$
\end{corollary}

\subsection{Invariant lightlike submanifolds}\label{S5.2}

In \cite{38}, researchers explored invariant lightlike submanifolds within golden semi-Riemannian manifolds and identified certain criteria for such submanifolds to qualify as locally product manifolds in the context of golden semi-Riemannian manifolds.

\begin{theorem} {\rm \cite{38}}
   Assume that $N$ is a lightlike submanifold of $(\bar{M}, \bar{g}, \phi)$, a golden semi-Riemannian manifold. For $N$ to be an invariant, it is required and sufficient that the induced structure $(\phi ', g)$ on $N$ is a golden semi-Riemannian structure.
\end{theorem}
\begin{theorem} {\rm \cite{38}}
   For a golden semi-Riemannian manifold $(\bar{M}, \bar{g}, \phi)$, let $N$ be an invariant lightlike submanifold. The induced structure $\phi '$ on the submanifold $N$ is parallel to the induced connection $\nabla$ if the golden structure $\phi$ is parallel to the Levi-Civita connection $\bar{\nabla}$ of $\bar{M}$.
\end{theorem}

\begin{theorem}{\rm \cite{38}}
    Assume that $N$ is an invariant lightlike submanifold of $(\bar{M}, \bar{g}, \phi)$, a golden semi-Riemannian manifold. Then, for the induced golden structure $\phi$ on $N$, the Nijenhuis tensor expression is given by
$$
\begin{gathered}
\mathcal{N}{_\phi '} (X, Y)=\left(\nabla_{\phi ' X} \phi ' \right) Y-\left(\nabla_{\phi ' Y} \phi ' \right) X+\left(\nabla_X \phi ' \right) \phi ' Y-\left(\nabla_Y \phi ' \right) \phi ' X \\
+\left(\nabla_Y \phi ' \right) X-\left(\nabla_X \phi ' \right) Y,
\end{gathered}
$$
for any $X, Y \in \Gamma(T N)$.
\end{theorem}

\begin{theorem} {\rm \cite{38}}
    For a golden semi-Riemannian manifold $(\bar{M}, \bar{g}, \phi)$, let $N$ be an invariant lightlike submanifold. As a result, $N$ is a totally geodesic, lightlike submanifold of $\bar{M}$.
\end{theorem}

 In \cite{36}, Poyraz and Yaşar explored various characteristics of semi-invariant lightlike submanifolds within golden semi-Riemannian manifolds, leading to the subsequent findings.
 
\begin{definition} {\rm \cite{36}}
 Consider a lightlike submanifold $(N, g, S(T N))$ of a golden semi-Riemannian manifold $(\bar{M}, \bar{g}, \phi)$.
 
{\rm (i)}  $N$ is an invariant lightlike submanifold if $\phi({Rad}(T N))={Rad}(T N)$ and $\phi(S(T N))=S(T N)$,

{\rm (ii)} $N$ is a semi-invariant lightlike submanifold if $$\phi({Rad}(TN)) \subset S(TN),\;\; \phi(\operatorname{ltr}(TN)) \subset S(TN)\;\; {\rm and}\;\:\; {S(TN^\perp)} \subset S(TN).$$

{\rm (iii)} $N$ is a radical anti-invariant lightlike submanifold if $\phi({Rad}(T N))=\operatorname{ltr}(T N)$.
\end{definition}

Let $(N, g, S(T N))$  be a semi-invariant lightlike submanifold of a golden semi-Riemannian manifold $(\bar{M}, \bar{g}, \phi)$. If we set $D_1=\phi({Rad}(T N)), D_2=\phi(I t r(T N))$ and $D_3=\phi\left(S\left(T N^{\perp}\right)\right)$ then we have
$$
S(T N)=D_0 \perp\left\{D_1 \oplus D_2\right\} \perp D_3 .
$$
Therefore,
$$
\begin{aligned}
T N & =D_0 \perp\left\{D_1 \oplus D_2\right\} \perp D_3 \perp {Rad}(T N) \quad {\rm and}\\
T \bar{M} & =D_0 \perp\left\{D_1 \oplus D_2\right\} \perp D_3 \perp\{{Rad}(T N) \oplus \operatorname{ltr}(T N)\} \perp S\left(T N^{\perp}\right) .
\end{aligned}
$$
According to this definition, one can write
$$D=D_0 \perp D_1 \perp {Rad}(T N)\;\; {\rm and} \;\;D^{\perp}=D_2 \perp D_3 .$$
Thus, we have
$T N=D \oplus D^{\perp} .$

\begin{proposition} {\rm \cite{36}}
Regarding $\phi$, we know that the distributions $D_0$ and $D$ are invariant distributions.
\end{proposition}

\begin{theorem} {\rm \cite{36}}
A golden semi-Riemannian manifold $(\bar{M}, \bar{g}, \phi)$ has a lightlike submanifold denoted by $N$. Hence, the ensuing assertions are equivalent.

{\rm (i)} $N$ is invariant

{\rm (ii)} The 1-forms $u_i$ and $w_\alpha$ vanish on TN $\forall$$i$ and $\alpha$.

{\rm (iii)} $\phi '$ is golden structure on $N$.
\end{theorem}

\begin{theorem} {\rm \cite{36}}
Let $N$ be a totally umbilical semi-invariant lightlike submanifold of a golden semi-Riemannian manifold ($\bar M, \bar g, \phi$).
Then we have $c_p = c_q = 0$.
\end{theorem}

\begin{remark} {\rm
In \cite{37}, Poyraz and Doğan identified several criteria for the integrability of distributions on semi-invariant lightlike submanifolds within golden semi-Riemannian manifolds, and explored both totally geodesic and mixed geodesic distributions of such submanifolds.}
\end{remark}

The authors explore the geometry of screen semi-invariant lightlike submanifolds within a golden semi-Riemannian manifold in \cite{39}.

\begin{definition} {\rm  \cite{39}}
    Let $N$ be a lightlike submanifold and $(\bar M, \bar g, \phi)$ be a golden semi-Riemannian manifold. Then, it is possible to define $N$ to be a screen semi-invariant lightlike submanifold of $\bar M$ if the following conditions are met:
$$
\begin{gathered}
\phi ({Rad}(T N)) \subseteq S(T N), \;\;
\phi (l t r(T N)) \subseteq S(T N) .
\end{gathered}
$$
Using the description above, one may also define a non-degenerate distribution $D_0$ for a screen semi-invariant lightlike submanifold of a golden semi-Riemannian manifold, such that $S(T N)$ is decomposed as follows:
$$
S(T N)=D_0 \perp D_1 \oplus D_2,
$$
where $D_1=\phi ({Rad}(T N))$ and $D_2=\phi (\operatorname{ltr}(T N))$.
\end{definition}

\begin{theorem} {\rm cite{39}}
   With a screen semi-invariant lightlike submanifold, let $\bar M$ be a golden semi-Riemannian manifold. The invariant distribution $D$ is integrable for any vector fields $U, V \in \Gamma(D)$, if and only if we have $$h^l\left(\phi V, \phi U\right)=h^l\left(U, \phi V\right)+h^l(U, V),$$ where $h^l$ is second fundamental form on $\Gamma(ltr(TN))$.
\end{theorem}

\begin{theorem} {\rm \cite{39}}
    Let $\bar M$ be a golden semi-Riemannian manifold with a screen semi-invariant lightlike submanifold. The radical distribution ${Rad}(T N)$ is integrable for any vector fields $U, V \in \Gamma({Rad}(T N))$ if and only if
$$
\begin{gathered}
\nabla^*_U \phi V-\nabla^*_V \phi U=\left(A^*_U-A^*_V\right)\quad
 {\rm or } \quad
\phi \nabla^*_U \phi V-\phi \nabla^*_V \phi U=\left(\nabla^*_U \phi V-\nabla^*_V \phi U\right),
\end{gathered}
$$
where $\nabla^*$ is linear connection on $(S(TN))$, $A^*$ shape  operator of distributions $(S(TN))$ and $Rad(TN)$.
\end{theorem}

\begin{theorem} {\rm \cite{39}}
    Let $N$ be a screen semi-invariant lightlike submanifold of a golden semi-Riemannian manifold $\bar M$. Then, for any $U, V \in \Gamma(S(T N))$, the screen distribution $S(T N)$ is integrable if and only if
$$
\begin{gathered}
\nabla^*_U \phi V-\nabla^*_V \phi U=\left(\nabla^*_U V-\nabla^*_V U\right) \;\; {\rm or}\;\;
\nabla^*_U \phi V=\nabla^*_V \phi U .
\end{gathered}
$$
\end{theorem}
\begin{theorem} {\rm \cite{39}}
    Let $N$ be a screen semi-invariant lightlike submanifold of a golden semi-Riemannian manifold $\bar M$. Then, for any $U \in \Gamma(T N)$ and $\xi \in \Gamma({Rad}(TN))$, induced connection $\nabla$ on $N$ is a metric connection if and only if one of the followings is satisfied:
$$
\begin{gathered}
\nabla^*_U \phi \xi ' = -A^*_{\xi '} U \;\; {\rm or} \;\;
A^*_{\xi '} U = 0
\end{gathered}
$$
\end{theorem}

\begin{remark} {\rm
    The essential criteria for these distributions to form complete geodesic foliations are also established in \cite{39}.}
\end{remark}

\begin{remark} {\rm
Poyraz has also conducted research on screen semi-invariant lightlike submanifolds within golden semi-Riemannian manifolds, as documented in \cite{43}. This study includes a discussion on various characteristics of these submanifolds and the verification of certain properties specific to totally umbilical screen semi-invariant lightlike submanifolds of golden semi-Riemannian manifolds. }
\end{remark}

\subsection{Transversal lightlike submanifolds}\label{S5.3}

Research on the geometry of screen transversal lightlike submanifolds and their anti-invariant counterparts in golden semi-Riemannian manifolds was conducted in \cite{40}. The authors in \cite{40} explored the geometry of distributions and established the essential and adequate conditions for the induced connections in these manifolds to qualify as metric connections. Furthermore, they provided a characterization of screen transversal anti-invariant lightlike submanifolds within golden semi-Riemannian manifolds.

\begin{definition} {\rm \cite{40}}
For a golden semi-Riemannian manifold $\bar M$, let $N$ be a lightlike submanifold. If
$$\phi (Rad (TN)) \subset S(TN)^\perp,$$
where $S(TN)^\perp$ refers to the screen distribution of the normal bundle $TN^\perp$. Then $N$ is a screen transversal lightlike submanifold of $\bar M$ golden semi-Riemannian manifold.
\end{definition}

\begin{definition} {\rm \cite{40}}
Let $N$ be a screen transversal lightlike submanifold of golden semi-Riemannian manifold $\bar M$. Then, if $\phi \left(S\left(T N\right)\right) \subset S\left(T N\right)^{\perp}$, then, $N$ is a screen transversal anti-invariant submanifold of $\bar M$.

If $N$ is a screen transversal anti-invariant submanifold of $\bar M$, then, 
$$
S\left(T N\right)^{\perp}=\phi  ({Rad} (T N)) \oplus \phi ({\rm l t r} (T N)) \oplus \phi (S(T N)) \perp D^\perp,
$$
where $D^\perp$ is the orthogonal complement non-degenerate distribution to $$\phi (Rad(TN)) \oplus \phi({\rm ltr}(TN)) \oplus \phi (S(TN)).$$
\end{definition}

\begin{proposition} {\rm \cite{40}}
    Let $N$ be a screen transversal anti-invariant lightlike submanifold of golden semi-Riemannian manifold $\bar M$. In such case, $D^\perp$ is an invariant distribution with regard to $\phi$.
\end{proposition}

\begin{theorem} {\rm \cite{40}}
    Let $N$ be a screen transversal anti-invariant lightlike submanifold within golden semi-Riemannian manifold $\bar M$. Radical distribution is integrable if and only if
$$
\nabla_U^s \phi V=\nabla_V^s \phi U
$$
for $U, V \in \Gamma({Rad} TN)$ and $\nabla^s $ refers to the screen connection on the screen bundle $S(TN).$
\end{theorem}

\begin{theorem} {\rm \cite{40}}
 Consider $N$ to be a screen transversal anti-invariant lightlike submanifold of golden semi-Riemannian manifold $\bar M$. In this instance, the screen distribution is integrable if and only if
$$
\nabla_U^s \phi V-\nabla_V^s \phi  U=h^s(U, V)-h^s(V, U),
$$
for $U, V \in \Gamma(S(T N))$ and $h^s$ is the screen second fundamental form.   
\end{theorem}

\begin{remark} {\rm
Erdogan also explores the geometry of screen transversal lightlike submanifolds, radical screen transversal lightlike submanifolds, and screen transversal anti-invariant lightlike submanifolds within golden semi-Riemannian manifolds in \cite{42}.}
\end{remark}

The study of the geometry of radical screen transversal lightlike submanifolds is carried out in \cite{40} by Erdoğan.

\begin{definition} {\rm \cite{40}}
Assume $N$ is a screen transversal lightlike submanifold of golden semi-Riemannian manifold $\bar M$. If $\phi \left(S\left(T N\right)\right)=S\left(T N \right)$, then $N$ is called a radical transversal lightlike screen submanifold of $\bar M$.
\end{definition}

\begin{theorem} {\rm \cite{40}}
    Let $N$ be a radical screen transversal lightlike submanifold of golden semi-Riemannian manifold $\bar M$. In this case, the screen distribution is integrable if and only if
$$
U, V \in \Gamma\left(S\left(T N \right)\right) . \quad h^s(U, \phi V)=h^s(V, \phi U),
$$
\end{theorem}

\begin{theorem} {\rm \cite{40}}
    If $N$ be a radical screen transversal lightlike submanifold of golden semi-Riemannian manifold $\bar  M$. The radical distribution is integrable if and only if
$$
A_{\phi V} U-A_{\phi U} V=A^*_U V-A^*_V U
$$
for $U, V \in \Gamma({Rad} (T N))$.
\end{theorem}

\begin{theorem} {\rm \cite{40}}
Given a golden semi-Riemannian manifold $\bar M$, let $N$ be a radical screen transversal lightlike submanifold of it. The screen distribution define totally geodesic foliation if and only if there is no component of $h^s(U, \phi V)-h^s(U, V)$ in $\phi l t r (T N)$ for $U, V \in \Gamma(S(T N))$.
\end{theorem}

\begin{theorem} {\rm \cite{40}}
    Let $N$ be a radical screen transversal lightlike submanifold of golden semi-Riemannian manifold $\bar M$. The radical distribution define totally geodesic foliation if and only if there is no component of $A_{\phi V} U$ in $S(T N)$ and $A^*_V U=0$, for $U, V \in \Gamma({Rad}TN))$.
\end{theorem}

\begin{theorem} {\rm \cite{40}}
    Let $N$ be a radical screen transversal lightlike submanifold of golden semi-Riemannian manifold $\bar M$. The connection induced on $N$ is a metric connection if and only if there is no component of $h^s(V, U)$ in $\phi (Rad (T N))$ or of $A_{\phi \xi '} U$ in $S(T N)$ for $U, V \in S(T N)$.
\end{theorem}

The study of the geometry of radical transversal lightlike submanifolds of golden Semi-Riemannian manifolds is carried out in \cite{41}. The authors have also investigated the geometry of the distributions and obtained the necessary and sufficient conditions for the induced connection on the manifolds to be a metric connection in \cite{41}.

\begin{definition} {\rm \cite{41}}
Consider a lightlike submanifold of a golden semi-Riemannian manifold as $\big(N, g, S\left(T N\right), S\left(T N\right)^{\perp}\big)$. If $\phi ({Rad(TN)})=\operatorname{ltr}\left(T N\right)$ and $\phi (S(T N))=S\left(T N\right)$ are satisfied, then the lightlike submanifold is called a radical transversal light-like submanifold.
\end{definition}

\begin{proposition} {\rm \cite{41}}
Let $\bar M$ be a manifold that is golden semi-Riemannian. In this instance, the manifold $\bar M$ does not have a 1-lightlike radical transversal lightlike submanifold.
\end{proposition}

\begin{theorem}  {\rm \cite{41}}
For a locally golden semi-Riemannian manifold $\bar M$, let $N$ be a radical transversal lightlike submanifold. The induced connection $\nabla$ on $D$ is metric connection if and only if the following conditions are met: $U \in \Gamma (TN)$ and $\xi ' \in \Gamma(Rad(TN))$ $$A_{\phi \xi '} U \in \Gamma (Rad(TN)).$$
\end{theorem}

\begin{theorem} {\rm \cite{41}}
For a locally golden semi-Riemannian manifold $\bar M$, let $N$ be a radical transversal lightlike submanifold. Here, the distribution $S(T N)$ must meet the following requirements in order to be integrable: for all  $U, V \in \Gamma(S(T N))$ we have
$$
h^l(V, S U)=h^l(U, S V) .
$$
\end{theorem}

\begin{note}
The necessary and sufficient condition for radical distribution definition of totally geodesic foliation on $N$ is also discussed in \cite{41}.
\end{note}

 The geometry of transversal lightlike submanifolds of golden semi-Riemannian manifolds is studied in \cite{41}. The authors have also looked into the geometry of the distributions and have determined what is required in order for the induced connection on the manifold to be a metric connection in \cite{41}.

\begin{definition} {\rm  \cite{41}}
Let $\big(N, g, S\left(T N\right), S\left(T N\right)^{\perp}\big)$ be lightlike submanifold of a golden semi-Riemannian manifold. If $\phi ({Rad(TN)})=\operatorname{ltr}\left(T N\right)$ and $\phi (S(T N)) \subseteq S\left(T N\right)$ are satisfied, then the lightlike submanifold is referred as transversal light-like submanifold.
\end{definition}

\begin{proposition} {\rm  \cite{41}}
 Let $N$ be its transversal lightlike submanifold of a golden semi-Riemannian manifold. In this case  the subbundle $\mu$ is orthogonal complement to the
$\phi (S(TN))$ in the $S(TN)$, the distribution is invariant under $\phi$.
\end{proposition}

\begin{proposition} {\rm \cite{41}}
Let $\bar M$ be golden semi-Riemannian manifold. In this case, there exists no 1-lightlike transversal submanifold of  $\bar M$.
\end{proposition}

\begin{note}
    It is also established in \cite{41} that the screen distribution constitutes both a necessary and sufficient condition for defining totally geodesic foliation on $N$. 
\end{note}

\subsection{CR-lightlike submanifolds of golden semi-Riemannian manifolds}\label{S5.4}

CR-lightlike submanifolds within a golden semi-Riemannian manifold were explored and characterized by Ahmad and Qayyoom in \cite{44}. They further examined various aspects of geodesic CR-submanifolds in a golden semi-Riemannian manifold. 

Additionally, they documented numerous significant findings regarding totally geodesic and totally umbilical CR-submanifolds in a golden Riemannian manifold.

\begin{definition} {\rm \cite{44}}
If both of the following two requirements are met, a submanifold $N$ of a golden semi-Riemannian manifold $\bar M$ is referred to as a CR-lightlike submanifold:

{\rm (i)} $\phi (Rad(TN))$ is a distribution on N such that $$Rad(TN) \cap \phi (Rad(TN) )= \{0\},$$

{\rm (ii)} There exist vector bundles $S(TN)$, $S(TN^\perp)$, ${\rm ltr}(TN)$, $D_0$, and $D'$ over $N$ such that
\[
S(TN) = \{\phi(Rad(TN)) \oplus D'\} \perp D_0, \quad \phi D_0 = D_0, \quad \phi D' = Z_1 \perp Z_2,
\]
where $D_0$ is a non-degenerate distribution on $N$ and $Z_1, Z_2,$ are vector bundles of
$ltr(TN)$ and $S(TN^\perp)$ respectively.
\end{definition}

\begin{lemma} {\rm \cite{44}}
Let the screen distribution be totally geodesic and $N$ a CR-lightlike submanifold of a golden semi-Riemannian manifold. Then $\nabla_X Y \in  \Gamma (S(TN))$, where $X, Y \in \Gamma (S(TN))$.
\end{lemma}

\begin{lemma} {\rm \cite{44}}
Let $N$ be a CR-lightlike submanifold a locally golden semi-Riemannian manifold $\bar M$. Then $\nabla_X \phi X = \phi \nabla_X X$ for any $X \in \Gamma (D_0)$.
\end{lemma}

\subsubsection{Geodesic CR-lightlike submanifolds}\label{S5.4.1}

\begin{definition} {\rm \cite{44}}
When the second fundamental form $h$ of a CR-lightlike submanifold of a golden semi-Riemannian manifold satisfies certain conditions, it is referred to as a mixed geodesic CR-lightlike submanifold:
$$ h(X, U) = 0, \quad \text{where } X \in \Gamma(D) \text{ and } U \in \Gamma(D'). $$
\end{definition}

\begin{definition} {\rm \cite{44}}
A CR-lightlike submanifold in golden semi-Riemannian manifold is referred as $D$-geodesic CR-lightlike submanifold if its second fundamental form $h$ satisfies
$h(X, Y ) = 0$ for $X, Y \in \Gamma(D)$.
\end{definition}

\begin{definition} {\rm \cite{44}}
A CR-lightlike submanifold in a golden semi-Riemannian manifold is known as $D'$-geodesic CR-lightlike submanifold if its second fundamental form $h$ satisfies 
$h(U, V) = 0$ for $U, V \in \Gamma(D')$.  
\end{definition}

\begin{theorem} {\rm \cite{44}}
Given a golden semi-Riemannian manifold $\bar M$, the submanifold $N$ of $\bar M$ should be CR-lightlike. In that case, $N$ is totally geodesic if and only if
$(Z_{\xi'} g)(X, Y ) = 0$
and
$(Z_W g)(X, Y ) = 0$
for any $X, Y \in \Gamma(TN), \xi ' \in \Gamma(Rad(TN))$ and $W \in \Gamma(S(TN^\perp))$.
\end{theorem}

\begin{theorem} {\rm \cite{44}}
Let $\bar M$ be a golden semi-Riemannian manifold and $N$ be a CR-lightlike submanifold of it. Then $N$ is mixed geodesic if and only if we have
$$A^*_{\xi'} X \in \Gamma(D_0\perp \phi Z_1)$$
and
$$A_W X \in \Gamma(D_0 \perp Rad(TN) \perp \phi Z_1),$$ for any $X \in \Gamma(D), \xi' \in \Gamma(Rad(TN))$ and $W \in \Gamma (S(TN^\perp)).$ 
\end{theorem}

\subsubsection{Totally umbilical CR-lightlike submanifolds}\label{S5.4.2}

\begin{definition}{\rm \cite{44}}
A lightlike submanifold $N$ within a semi-Riemannian manifold $\bar M$ is termed totally umbilical within $\bar M$ if it possesses a smooth transversal vector field $\mathcal{T}$, belonging to $\Gamma({\rm Trace}(TN))$ on $N$, called the transversal curvature vector field of $N$, such that $h(X, Y ) = \mathcal{T} g(X, Y )$
for $X, Y \in \Gamma(TN).$
\end{definition}

\begin{theorem} {\rm \cite{44}}
Consider $N$ as a totally umbilical CR-lightlike submanifold of $\bar M$, a golden manifold. The CR-lightlike sectional curvature of $N$ then disappears, that is
$K(\pi) = 0$ for any CR-lightlike section $\pi$.
\end{theorem}

\begin{remark} {\rm
    Analogous findings for CR-lightlike, totally geodesic lightlike, and totally umbilical lightlike submanifolds of a golden semi-Riemannian manifold are discussed in \cite{qa}.}
\end{remark}

\subsection{Half lightlike submanifolds}\label{S5.5}

Submanifolds of lightlike type with codimension 2 are termed either half lightlike or coisotropic, depending on the rank of their radical distribution. These are further divided into two subclasses \cite{46}. A half lightlike submanifold represents a particular instance of the broader r-lightlike submanifolds where $r = 1$. Its geometric structure is more comprehensive compared to that of a coisotropic submanifold or a lightlike hypersurface \cite{33}.

Poyraz et al. examined half lightlike submanifolds within a golden Semi-Riemannian manifold in their work referenced as \cite{45}.They demonstrated the absence of radical anti-invariant half lightlike submanifolds within such a context and provided findings regarding screen semi-invariant half lightlike submanifolds. Additionally, their study encompassed screen conformal half lightlike submanifolds within a golden Semi-Riemannian manifold.

\begin{definition} {\rm \cite{45}}
Suppose $N$ is a half lightlike submanifold of a golden Semi-Riemannian manifold $\bar M$.

{\rm (i)} $N$ is an invariant half lightlike submanifold if $\phi (TN) = TN$.

{\rm (ii)} $N$ is a screen semi-invariant half lightlike submanifold if $\phi (Rad(TN)) \subset S(TN)$ and $\phi ({\rm ltr}(TN)) \subset S(TN)$.

{\rm (iii)} $N$ is a radical anti-invariant half lightlike submanifold if $\phi (Rad(TN)) = {\rm ltr}(TN).$
\end{definition}

\begin{theorem} {\rm \cite{45}}
If $N$ is a half lightlike submanifold of a golden semi-Riemannian manifold $\bar M$, then the following three statements are equivalent.

{\rm (i)} $N$ is invariant.

{\rm (ii)} $u_1$ and $u_2$ vanish on $N$, where $u_1$ and $u_2$ are 1-forms on $N$.

{\rm (iii)} $\phi '$ is a golden structure on $N$.
\end{theorem}

\begin{theorem} {\rm \cite{45}}
No radical anti-invariant half lightlike submanifold exists within a golden semi-Riemannian manifold.
\end{theorem}

\begin{theorem}{\rm \cite{45}}
Suppose $N$ is a screen conformal totally umbilical screen semi-invariant half lightlike submanifold of a golden semi-Riemannian manifold $\bar M$. Then, both $N$ and the screen distribution $S(TN)$ are totally geodesic.
\end{theorem}

\begin{corollary} {\rm \cite{45}}
For a screen semi-invariant half lightlike submanifold $N$ of a golden semi-Riemannian manifold $\bar M$, the condition $$h_1(X,Y) = 0,$$ holds, indicating that the vector field $V$ results in the degeneration of the local second fundamental form of $N$.     
\end{corollary}

\begin{theorem} {\rm \cite{45}}
If $N$ is a screen semi-invariant half lightlike submanifold of a golden semi-Riemannian manifold $\bar M$, then the distribution $D$ is integrable if and only if the following condition is satisfied:
$$h_1 (\phi X,\phi Y) = h_1(\phi X, Y) + h_1(X, Y )\quad {and}$$
$$h_2(\phi X,\phi Y) = h_2(\phi X, Y) + h_2(X, Y ),$$
for any $X, Y \in \Gamma(D),$ where $h_1$ and $h_2$ are the local second fundamental forms of $N.$
\end{theorem}

\begin{definition} {\rm \cite{45}}
If $N$ is a screen semi-invariant half lightlike submanifold of a golden semi-Riemannian manifold $\bar M$, then $N$ is mixed totally geodesic if and only if
$$h_1(X, Y ) = h_2(X, Y ) = 0,$$
for any $X \in \Gamma(D)$ and $Y \in \Gamma(D^\perp).$
\end{definition}

\begin{remark} {\rm
Results regarding a totally umbilical screen semi-invariant half lightlike submanifold of a golden semi-Riemannian manifold  is also discussed in \cite{45}.}
\end{remark}

\subsection{Generic lightlike submanifolds}\label{S5.6}

If $N$ is a real $r$-lightlike submanifold of a semi-Riemannian manifold $\bar M$, it's called a {\it golden generic lightlike submanifold} if the screen distribution $S(TN)$ of $N$ is characterized by

\begin{align}
S(TN) &= \phi (S(TN)^\perp) \perp D_0 \notag \\
& = \phi (Rad(TN)) \oplus \phi ({\rm ltr}(TN)) \perp \phi (S(TN)^\perp) \perp D_0
\end{align}
where $D_0$ is a non-degenerate almost complex distribution on $N$ with respect to $\phi$ \cite{48}.

\begin{theorem} {\rm \cite{48}}
For a golden generic lightlike submanifold $N$ of a golden semi-Riemannian manifold $\bar M$, the Nijenhuis tensor field, concerning the golden structure $\phi$, is null.
\end{theorem}

\begin{theorem}{\rm \cite{48}}
If $N$ is a golden generic lightlike submanifold of a golden semi-Riemannian manifold $\bar M$, then $g$ serves as a golden structure on $D$.
\end{theorem}

\begin{definition} {\rm \cite{48}}
A golden generic lightlike submanifold $N$ is termed mixed geodesic if its second fundamental form $h$ fulfills the condition:
$$ h(Y, Z) = 0, \quad \text{for } Y \in \Gamma(D) \text{ and } Z \in \Gamma(D'). $$
\end{definition}

\begin{theorem} {\rm \cite{48}}
If $N$ is a totally umbilical golden generic lightlike submanifold of a golden semi-Riemannian manifold $\bar M$, then the distribution $D$ is inherently integrable.
\end{theorem}

\begin{remark} {\rm
    Results of  Minimal golden generic lightlike submanifolds is also discussed in \cite{48}.}
\end{remark}

Yadav and Kumar have conducted research on screen generic lightlike submanifolds as documented in \cite{49}, yielding the following outcomes:

\begin{definition} {\rm \cite{49}}
If $N$ is a real $r$-lightlike submanifold of a golden semi-Riemannian manifold $\bar M$, it's termed a screen generic lightlike submanifold of $\bar M$ if the following conditions hold:

{\rm (i)}  $Rad(TN)$ is invariant with respect to $\phi$, that is
$$\phi (Rad (TN)) = Rad(TN).$$

{\rm (ii)} There exist a subbundles $D_0$ of $S(TN)$ such that
$$D_0 = \phi (S(TN)) \cap S(TN),$$
where $D_0$ is a non-degenerate distribution on $N$.
\end{definition}

\begin{proposition} {\rm \cite{49}}
A generic $r$-lightlike submanifold is a screen generic lightlike submanifold with $\mu = 0$, where $\mu$ represents a non-degenerate invariant distribution.
\end{proposition}

\begin{proposition} {\rm \cite{49}}
No coisotropic, isotropic, or totally lightlike proper screen generic lightlike submanifold exists within a golden semi-Riemannian manifold $\bar M$. Any screen generic isotropic, coisotropic, or totally lightlike submanifold in $\bar M$ is an invariant submanifold.
\end{proposition}

\begin{definition} {\rm \cite{49}}
A screen generic lightlike submanifold of a golden semi-Riemannian manifold is termed a $D$-geodesic screen generic lightlike submanifold if its second fundamental form $h$ satisfies the condition: $h(X, Y ) = 0$ for any $X, Y \in \Gamma(D).$
\end{definition}

\begin{definition} {\rm \cite{49}}
A screen generic lightlike submanifold $N$ of a golden semi-Riemannian manifold is called a mixed geodesic screen generic lightlike submanifold if its second fundamental form $h$ satisfies the condition: $h(X, Y ) = 0,$ for any $X \in \Gamma(D)$ and $ Y \in \Gamma(D').$
Thus $N$ is a mixed geodesic screen generic lightlike submanifold if we have
$h^l (X, Y )=0 \quad \text{and} \quad h^s(X, Y )=0$
for any $X \in \Gamma(D)$ and $Y \in \Gamma(D').$
\end{definition}

\begin{remark} {\rm
In \cite{49}, the conditions for the induced connection to meet the criteria of a metric connection are discussed, along with the classification of totally umbilical screen generic lightlike submanifolds of golden semi-Riemannian manifolds as totally geodesic. Moreover, the paper delves into the study of minimal screen generic lightlike submanifolds of golden semi-Riemannian manifolds.}
\end{remark}

\subsection{Golden GCR-lightlike submanifolds}\label{S5.7}

\begin{definition} {\rm \cite{50}}
A real lightlike submanifold $(N, g, S(TN))$ of a golden semi-Riemannian manifold $(\bar M,\bar g,\phi)$ is termed a golden generalized Cauchy–Riemann (GCR)-lightlike submanifold if the following conditions hold:

{\rm (i)} There exist two subbundles $D_1$ and $D_2$ of $Rad(TN)$ such that
$$Rad(TN) = D_1 \oplus D_2, \;\; \phi (D_1) = D_1,\;\;  \phi (D_2) \subset S(TN).$$

{\rm (ii)} There exist two subbundles $D_0$ and $D'$ of $S(TN)$ such that
$$S(TN) = \{\phi (D_2) \oplus D'\} \perp D_0, \;\; \phi (D_0) = D_0, \;\; \phi (Z_1 \perp Z_2) = D',$$
where $D_0$ is a non-degenerate distribution on $N$, $Z_1$ and $Z_2$ are vector subbundles of $ltr (TN)$ and $S(TN^\perp)$, respectively.

Let $\phi (Z_1) = N_1$ and $\phi (Z_2) = N_2$. Then we have
$$D'= \phi (Z_1) \perp \phi (Z_2) = N_1 \perp N_2.$$
Thus, the following decomposition:
$$TN = D \oplus D', \;\; D = Rad(TN) \perp D_0 \perp \phi (D_2).$$
Thus, $N$ is a proper golden GCR-lightlike submanifold of a golden
semi-Riemannian manifold if $D_0 \neq {0}, D_1 \neq {0}, D_2 \neq {0}$ and $Z_2 \neq {0}.$
\end{definition}

\begin{theorem} {\rm \cite{50}}
If $N$ is a golden GCR-lightlike submanifold of a golden semi-Riemannian manifold, then $\phi'$ serves as a golden structure on $D$.
\end{theorem}

\begin{theorem} {\rm \cite{50}}
If $N$ is a golden GCR-lightlike submanifold of a golden semi-Riemannian manifold, then the distribution $D$ is integrable if and only if we satisfy the condition:

{\rm (i)} $\bar g (h^l (X,\phi Y), \xi')=\bar g (h^l (Y, \phi X), \xi')$ and

{\rm (ii)} $\bar g (h^s (X,\phi Y), W)=\bar g (h^s (Y, \phi X), W)$

\noindent for any $X, Y \in \Gamma(D), \xi' \in \Gamma(D_2)$
and $W \in \Gamma(Z_2)$.
\end{theorem}

\begin{definition} {\rm \cite{50}}
 A golden GCR-lightlike submanifold of a golden semi-Riemannian manifold is termed a D-geodesic golden GCR-lightlike submanifold if its second fundamental form $h$ satisfies: $h(X, Y ) = 0$ for any $X, Y \in \Gamma(D)$.   
\end{definition}

Discussions on minimal golden GCR-lightlike submanifolds can also be found in \cite{50}.

\begin{definition} {\rm \cite{50}}
A lightlike submanifold isometrically immersed in a semi-Riemannian manifold is considered minimal if it meets the following two conditions:

{\rm (i)} $ h^s = 0$ on $Rad(TN)$ and

{\rm (ii)} $ \operatorname{Trace} h = 0,$ where Trace is written with respect to $g$ restricted to $S(TN)$.
\end{definition}

\begin{remark} {\rm
In \cite{50} Poyraz have also provided examples and results of minimal golden GCR-lightlike submanifolds of golden semi-Riemannian manifolds.}
\end{remark}

\begin{remark} {\rm
Poyraz discussed comparable findings on the geometry of golden GCR-lightlike submanifolds in golden semi-Riemannian manifolds in \cite{51}.}
\end{remark}

\subsection{Slant lightlike submanifolds}\label{S5.8}

Acet explored screen pseudo slant lightlike submanifolds in golden semi-Riemannian manifolds in \cite{52}.

\begin{definition} {\rm \cite{52}}
If $N$ is a lightlike submanifold of a golden semi-Riemannian manifold $\bar M$, it is termed a screen pseudo-slant submanifold of $\bar M$ if the following conditions are satisfied:

{\rm (i)} The radical distribution $Rad(TN)$ is an invariant distribution with respect to $\phi$, i.e., $\phi (Rad(TN)) = Rad(TN),$

{\rm (ii)} There exist  non-degenerate orthogonal distributions $D_0$ and $D^\perp$ on $N$ such that $S(TN) = D_0 \perp D^\perp,$

{\rm (iii)} The distribution $D_0$ is anti-invariant, i.e., $\phi (D_0) \subset S(TN^\perp).$

{\rm (iv)} The distribution $D^\perp$ is slant with angle $\theta (\neq \frac{\pi}{2})$ i.e., for each $x \in N$ and each non-zero vector $X \in (D^\perp)_x,$ the angle $\theta$ between $\phi X$ and the vector subspace $(D^\perp)_x$ is a constant $(\neq \frac{\pi}{2})$ which is independent of the choice of $x \in N$ and $X \in (D^\perp)_x.$
\end{definition}

\begin{remark} \hskip-.01in {\rm
Some non-trivial examples of  screen pseudo-slant lightlike submanifold in a golden semi-Riemannian manifold are studied and the conditions for integrability of distributions of a screen pseudo-slant lightlike submanifold of a golden semi-Riemannian manifold are also provided in \cite{52}.}
\end{remark}

In \cite{53}, Yadav and Kumar explored characteristics of screen generic lightlike submanifolds in golden semi-Riemannian manifolds.

\begin{definition} {\rm \cite{53}}
If $N$ is a $2q$-lightlike submanifold of a golden semi-Riemannian manifold $\bar M$ with index $2q$ (where $q$ is an integer indicating the dimension of the radical bundle of the lightlike submanifold $N$ within the golden semi-Riemannian manifold $\bar M)$, and if $2q < \text{dim}(N)$, then $N$ is termed a screen semi-slant lightlike submanifold of $\bar M$ if the following conditions are satisfied:

{\rm (i)} $Rad(TN)$ is invariant with respect to $\phi$, i.e., $\phi (Rad(TN)) = Rad(TN),$

{\rm (ii)} there exist non-degenerate orthogonal distributions $D_0$ and $D^\perp$ on $N$ such that $S(TN) = D_0 \oplus_{orth}  D^\perp,$

{\rm (iii)} the distribution $D_0$ is an invariant distribution, i.e., $\phi D_0 = D_0,$

{\rm (iv)} the distribution $D^\perp$ is slant with angle $\theta (\neq 0),$ i.e. for each $x \in N$ and each non-zero vector $X \in (D^\perp)_x,$ the angle $\theta$ between $\phi X$ and the vector subspace $(D^\perp)_x$ is a non-zero constant, which is independent of the choice of $x \in N$ and $X \in (D^\perp)_x.$
\end{definition}

\begin{theorem} {\rm \cite{53}}
$N$ is a $2q$-lightlike submanifold of a golden semi-Riemannian manifold $\bar M$, it qualifies as a screen semi-slant lightlike submanifold of $\bar M$ if and only if

{\rm (i)} $ltr(TN)$ and $D_0$ are invariant with respect to $\phi$,

{\rm (ii)} there exists a constant $c \in [0, 1)$ such that $\phi'^2 X = c (\phi X + X),$ for any $X \in \Gamma (D^\perp),$ where $D_0$ and $D^\perp$ are non-degenerate orthogonal distribution and $\phi'$ is $(1,1)$-tensor field on $N$ such that $S(TN) = D_0 \oplus_{orth} D^\perp .$ Moreover, in this case $c = \cos^2 {\theta}$ and $\theta$ is slant angle of $D^\perp$.
\end{theorem}

\begin{corollary} {\rm \cite{53}}
For $N$ as a screen semi-slant lightlike submanifold of a golden semi-Riemannian manifold $\bar M$ with a slant angle $\theta$, for any $X, Y \in \Gamma(D^\perp)$, we have:

{\rm (i)} $g(\phi' X,\phi' Y)= \cos^2 {\theta} (g(X,Y)+g(X,\phi' Y)),$

{\rm (ii)} $g(\mathcal{T} X,\mathcal{T} Y)=\sin^2 {\theta}(g(X,Y)+g(\phi' X,Y)).$
\end{corollary}

\begin{remark} {\rm
Yadav and Kumar investigate the essential and adequate
 conditions for the integrability and totally geodesic foliation of the distributions $Rad(TN)$, $D_0$, and $D^\perp$ of screen semi-slant lightlike submanifolds of golden semi-Riemannian manifolds in \cite{53}.}
\end{remark}

\begin{definition} {\rm \cite{54}}
If $N$ is a $2q$-lightlike submanifold of a golden semi-Riemannian manifold $\bar M$ with index $2q$, and if $2q < \text{dim}(N)$, then $N$ is termed a semi-slant lightlike submanifold of $\bar M$ if the following conditions are satisfied:

{\rm (i)} $\phi Rad(TN)$ is a distribution on $N$ such that $Rad (TN) \cap \phi (Rad (TN))=\{0\};$

{\rm (ii)} there exist non-degenerate orthogonal complementary distributions $D_0$ and $D^\perp$ on $N$ such that $S(TN)= (\phi (Rad (TN)) \oplus \phi ({\rm ltr} (TN)) \oplus_{orth} D_0 \oplus_{orth} D^\perp ;$

{\rm (iii)} the distribution $D_0$ is an invariant distribution, i.e., $\phi D_0 =D_0 ;$

{\rm (iv)} the distribution $D^\perp$ is slant with angle $\theta (\neq 0),$ i.e. for each $x \in N$ and each non-zero vector $X \in (D^\perp)_x,$ the angle $\theta$ between $\phi X$ and the vector subspace $(D^\perp)_x$ is a non-zero constant, which is independent of the choice of $x \in N$ and $X \in (D^\perp)_x.$
\end{definition}

\begin{proposition} {\rm \cite{54}}
A golden semi-Riemannian manifold $(\bar M, \bar g, \phi)$ does not have any isotropic or totally lightlike proper semi-slant lightlike submanifolds.
\end{proposition}

\begin{theorem} {\rm \cite{54}}
A golden semi-Riemannian manifold $\bar{M}$ of index $2q$ has a q-lightlike submanifold, denoted by $N$. Then $N$ is a semi-slant lightlike submanifold of $\bar{M}$ if and only if

{\rm (i)} $\phi {Rad}(T N)$ is a distribution on $N$ such that ${Rad}(T N) \cap \phi ({Rad}(T N))=0$;

{\rm (ii)} the screen distribution $S(T N)$ splits as $$S(T N)=\phi ({Rad}(T N)) \oplus \phi (\operatorname{ltr}(T N))\oplus_{\text {orth}} D_0 \oplus_{\text {orth}} D^\perp;$$

{\rm (iii)} there exists a constant $c \in[0,1)$ such that $\phi'^2 X=c(\phi' X+X)$ for any $X \in \Gamma\left(D^\perp\right)$. Moreover, in this case we have $c=\cos ^2 \theta$, where $\theta$ is the slant angle of $D^\perp$.
\end{theorem}

\begin{remark} {\rm
Kumar and Yadav discussed the necessary and sufficient conditions for the integrability of distributions and the geometry of the leaves of the foliation determined by the distributions in \cite{54}.}
\end{remark}

\begin{remark} {\rm
Kumar and Yadav in \cite{55} explore the concept of screen slant lightlike submanifolds within golden semi-Riemannian manifolds, along with the essential conditions for the integrability of distributions and the structural details of the foliation's leaves governed by these distributions.}
\end{remark}

\begin{definition} {\rm \cite{56}}
A golden semi-Riemannian manifold $\bar M$ of index $2 q$, such that $2 q<\operatorname{dim}(N)$, has a $q$-lightlike submanifold called $N$. Then, $N$ is a bi-slant lightlike submanifold of $\bar{M}$ if the following conditions are satisfied:

{\rm (i)} $\phi {Rad}(T N)$ is a distribution on $N$ such that ${Rad}(T N) \cap \phi ({Rad} (T N))=\{0\}$;

{\rm (ii)} there exists non-degenerate orthogonal distributions $D, D_0$ and $D^\perp$ on $N$ such that
$$
S(T N)=\phi ({Rad}(T N)) \oplus \phi (\operatorname{ltr}(T N)) \perp D \perp D_0 \perp D^\perp ;
$$

{\rm (iii)} the distribution $D$ is an invariant distribution, i.e., $\phi D=D$;

{\rm (iv)} the distribution $D_0$ is slant with angle $\theta_1(\neq 0)$, i.e., for each $x \in N$ and each non-zero vector $X \in$ $\left(D_0\right)_x$, the angle $\theta_1$ between $\phi X$ and the vector space $(D_0)_x$ is a non-zero constant, which is independent of the choice of $x \in N$ and $X \in (D_0)_x$;

{\rm (v)} the distribution $D^\perp$ is slant with angle $\theta_2(\neq 0)$, i.e., for each $x \in N$ and each non-zero vector $X \in$ $(D^\perp)_x$, the angle $\theta_2$ between $\phi X$ and the vector space $\left(D^\perp\right)_x$ is a non-zero constant, which is independent of the choice of $x \in N$ and $X \in (D^\perp)_x$.

The $\theta_1$ and $\theta_2$ constant angles are referred to as the slant angles of the distributions $D_0$ and $D^\perp$, respectively. A bi-slant lightlike submanifold is said to be proper if $D_0 \neq\{0\}, D^\perp \neq\{0\}$ and $\theta_1 \neq \frac{\pi}{2}, \theta_2 \neq \frac{\pi}{2}$.
\end{definition}

\begin{theorem}  {\rm \cite{56}}
If $N$ is a $q$-lightlike submanifold of a golden semi-Riemannian manifold $\bar{M}$ with index $2q$, then $N$ is a bi-slant lightlike submanifold if and only if

{\rm (i)} there exist a distribution $\phi ({Rad}(T N))$ on $N$ such that $${Rad}(T N) \cap \phi ({Rad}(T N))=\{0\};$$

{\rm (ii)} there exist a screen distribution $S(T N)$ which can be spitted as
$$
S(T N)=\phi ({Rad}(T N)) \oplus \phi (\operatorname{ltr}(T N)) \perp D \perp D_0 \perp D^\perp$$
such that $D$ is an invariant distribution on $N$, i.e., $\phi (D)=D$;

{\rm (iii)} there exist a constant $c_1 \in[0,1)$ such that $\phi'^2 X=c_1(\phi' +I) X$  for any $ X \in \Gamma\left(D_0\right)$;

{\rm (iv)} there exist a constant $c_2 \in[0,1)$ such that $\phi'^2 X=c_2(\phi' +I) X$ for any $ X \in \Gamma\left(D^\perp\right)$.

In that case, $c_1=\cos ^2 \theta_1$ and $c_2=\cos ^2 \theta_2$, where $\theta_1$ and $\theta_2$ represents the slant angles of $D_0$ and $D^\perp$ respectively.
\end{theorem}

\begin{remark} {\rm
In \cite{56}, the integrability conditions of distributions on bi-slant lightlike submanifolds and the necessary and sufficient conditions for foliations determined by distributions on bi-slant lightlike submanifolds of golden semi-Riemannian manifolds to be geodesic are obtained.}
\end{remark}

\begin{definition} {\rm \cite{57} }
Given a golden semi-Riemannian manifold $\bar M$ of index 2q, let $N$ be a q-lightlike submanifold such that $2q < dim (N)$. If the following criteria are met, then $N$ is a slant lightlike submanifold of $\bar M$:

{\rm (i)} $\phi Rad(TN)$ is a distribution on $N$ such that $Rad(TN) \cap \phi (Rad(TN)) = \{0\}$,

{\rm (ii)} there exists a non-degenerate orthogonal complementary distribution $D$ on $N$
such that $$S(TN) = \phi (Rad(TN)) \oplus \phi ({\rm ltr}(TN)) \oplus_{orth} D,$$

{\rm (iii)} the distribution $D$ is slant with angle $\theta (\neq 0)$, i.e. for each $x \in N$ and each non-zero vector $X \in (D)_x,$ the angle $\theta$ between $\phi X$ and the vector subspace $(D)_x$ is a non-zero constant, which is independent of the choice of $x \in N$ and $X \in (D)_x$.
\end{definition}

\begin{theorem} {\rm \cite{57}}
If $N$ is a $q$-lightlike submanifold of a golden semi-Riemannian manifold $\bar{M}$ with index $2q$, then $N$ is a slant lightlike submanifold of $\bar{M}$ if and only if we have:

{\rm (i)} $\phi (Rad(TN))$ is a distribution on $N$ such that $Rad(TN) \cap \phi (Rad(TN))=0,$ 

{\rm (ii)} the screen distribution $S(T N)$ split as $$S(T N)=\phi ({Rad} (T N)) \oplus \phi ({\rm l t r}(T N)) \oplus_{\text {orth}} D, $$

{\rm (iii)} there exists a constant $c \in[0,1)$ such that $\phi'^2 X=c(\phi' X+X)$ for any $X \in \Gamma(D)$. Moreover, in this case $c=\cos ^2 \theta$ and $\theta$ is the slant angle of $D$.
\end{theorem}

\begin{remark} {\rm
The criteria for integrability of distributions and the curvature characteristics of slant lightlike submanifolds in golden semi-Riemannian manifolds are further explored in \cite{57}.}
\end{remark}

\section{Lightlike Submanifolds of Meta-Golden Manifolds}\label{S6}

F.E. Erdoğan et. al {\rm \cite{15}} proved the following results.

\begin{theorem} {\rm \cite{15}}
Given an AMGsR manifold $(\bar M,\phi, F, \bar g)$, let $N$ be a lightlike hypersurface of $\bar M$. Here, the almost golden structure $\bar g$ induces a structure $(\phi, g,u, X)$ on $N$ that satisfies the following equalities: $$\phi'^2 X=\phi' X + X,\;\; u(\phi' X)=0,\;\; \phi' U=0,$$ $$u(W^2)-u(W)-1=0,\;\; g(\phi' X,\phi' Y)= g(\phi' X, Y)+ g(X, Y),$$ where for $X, Y$ and $U \in \Gamma(TN), W \in \Gamma(ltr (TN)).$
\end{theorem}

\begin{proposition} {\rm \cite{15}}
Assume that the manifold $(\bar M,\phi, F, \bar g)$ is an AMGsR. Then $\nabla \phi F=0$.
\end{proposition}

\begin{definition}{\rm \cite{15}}
Consider the AMGsR manifold $(\bar M,\phi, F, \bar g,)$ and the lightlike hypersurface $N$ of $\bar M$. Then

{\rm (i)} if $\phi F (TN) \subset TN$ $N$ is called invariant,

{\rm (ii)} if $\phi F (Rad(TN)) \subset S(TN)$ and $\phi F (ltr(TN)) \subset S(TN)$, then $N$ is called  screen semi-invariant,

{\rm (iii)} if $\phi F (Rad(TN)) \subset (ltr(TN))$, then $N$ is called  a radical anti-invariant,
lightlike hypersurface.
\end{definition}

\begin{theorem}  {\rm \cite{15}}
Every AMGsR manifold does not admit a radical anti-invariant lightlike hypersurface.
\end{theorem}

The screen semi-invariant lightlike hypersurface of almost meta-golden semi-Riemannian manifolds is also studied in  {\rm \cite{15}}.

\begin{corollary}  {\rm \cite{15}}
Given an AMGsR manifold $(\bar M,\phi, F, \bar g,)$ and a screen semi-invariant lightlike hypersurface $N$ of $\bar M$, we have $h(X, Z)=0$ for any $X, Z \in \Gamma(TN)$.
\end{corollary}

\begin{proposition}  {\rm \cite{15}}
 Let $N$ be a screen semi-invariant lightlike hypersurface of an AMGsR manifold $(\bar M,\phi, F, \bar g,)$.  Then $F D_0\subset S(TN)$ for the distribution $D_0.$
\end{proposition}

\begin{corollary}  {\rm \cite{15}}
Let $N$  be a screen semi-invariant lightlike hypersurface of an AMGsR manifold $(\bar M,\phi, F, \bar g)$. Then the distribution $D_0$ is $F$-invariant.
\end{corollary}

\begin{theorem}  {\rm \cite{15}} Let $N$  be a screen semi-invariant lightlike hypersurface of an AMGsR manifold $(\bar M,\phi, F, \bar g)$. 
Then the distribution $D$ can be integrated if and only if we have $$h (F Y, F X )= h (X, F \phi Y )-h (X, F Y )+ h (X , Y )$$ for any $X , Y \in \Gamma(D).$
\end{theorem}

\begin{theorem}  {\rm \cite{15}}
Consider an AMGsR manifold $(\bar M,\phi, F, \bar g,)$ and $N$ is a totally umbilici screen semi-invariant lightlike hypersurface of $\bar M$. Then $N$ is totally geodesic.  
\end{theorem}

\begin{theorem}  {\rm \cite{15}}
In an AMGsR manifold $(\bar M,\phi , F, \bar g)$, suppose $N$ is a screen semi-invariant lightlike hypersurface. If the screen distribution $(S(TN))$ is totally umbilical, then it is also totally geodesic.
\end{theorem}

\section{Lightlike Submanifolds of An  Almost Norden Golden Semi-Riemannian Manifold}\label{S7}

Investigations into certain classifications of lightlike hypersurfaces within almost Norden golden semi-Riemannian manifolds, including invariant and screen semi-invariant types, are discussed in \cite{er}.

\begin{theorem} \hskip-.03in {\rm \cite{er}}
In an almost Norden golden semi-Riemannian manifold $(\bar M,\bar g,\phi)$ with $(N, g)$ as a lightlike hypersurface of $\bar M$, the following three statements are equivalent:

{\rm (i)} $N$ is $\phi$-invariant.

{\rm (ii)} The 1-form $u$ vanishes on $N$.

{\rm (iii)} $\phi$ is an almost Norden golden structure on $N$.
\end{theorem}

\begin{theorem}  {\rm \cite{er}}
The almost Norden golden semi-Riemannian has no radical anti-invariant lightlike hypersurface.
\end{theorem}

\begin{proposition} {\rm \cite{er}}
Let $N$ be a screen semi-invariant lightlike hypersurface and $(\bar M,\bar g,\phi)$ be an almost Norden golden semi-Riemannian manifold. A $\phi$-invariant distribution is then $D_0$.   
\end{proposition}

\begin{theorem}  {\rm \cite{er}}
 Let $(\bar M,\bar g,\phi)$ be an almost Norden golden semi-Riemannian manifold and $N$ be a
screen semi-invariant lightlike hypersurface. Then the following three statements are equivalent:

{\rm (i)} $D$ is a parallel distribution.

{\rm (ii)} $D$ is totally geodesic.

{\rm (iii)} $(\nabla_X \phi)Y = 0,$ where $X, Y \in \Gamma(D).$ 
\end{theorem}

\begin{theorem}  \hskip-.03in {\rm \cite{er}}
In an almost Norden golden semi-Riemannian manifold $(\bar M,\bar g,\phi)$ if $N$ is a screen semi-invariant lightlike hypersurface that is totally umbilical, then $N$ is also totally geodesic in $\bar M$.
\end{theorem}

\begin{remark} {\rm
In \cite{er}, one can find illustrations of invariant and screen semi-invariant lightlike hypersurfaces in almost Norden golden semi-Riemannian manifolds.}
\end{remark}

\section{Warped Product of Screen-Real Lightlike Submanifolds of Golden Semi-Riemannian Manifolds}\label{S8}

The study on warped product of screen real lightlike submanifolds in a golden semi-Riemannian manifold is discussed in \cite{g}.

\begin{theorem} {\rm \cite{g}}
 Let $N=N_1 \times_f N_T$ be a warped product lightlike submanifold, then, for any $X \in Rad(TN),\, Y \in \Gamma(S(T\bar M)),$ we have $\nabla_X Y \in \Gamma(S(TN))$.
\end{theorem}

\begin{theorem} {\rm \cite{g}}
    If $(N, g, S(TN))$ is an irrotational screen-real m-lightlike submanifold of a golden semi-Riemannian manifold, then the induced connection is metric.
\end{theorem}

\begin{theorem} {\rm \cite{g}}
    There is no concept of an irrotational screen-real $m$-lightlike submanifold that can be expressed as warped product lightlike submanifolds.
\end{theorem}

\section{Chen Invariants and Inequalities}\label{S9}

Let $\bar M^{n}$ be a Riemannian $n$-manifold.   Let us choose a local field of orthonormal frame
$e_1,\ldots,e_n$ on $\bar M^n$.  Denote by $K(e_i\wedge e_j)$  the sectional curvatures of $\bar M^{n}$ of the plane section spanned by $e_i$ and $e_j$. 

The {\it scalar curvature} $\tau$ of $\bar M^{n}$ at $p$ is defined by
\begin{align}\label{9.1} \tau(p)=\sum_{1\leq i<j\leq n} K(e_i\wedge e_j).\end{align}
Similarly, if $\mathcal{L}$ is an $\ell$-dimensional linear subspace of $T_{p}\bar M^{n}$ with $2<\ell<n$, then the scalar curvature $\rho(\mathcal{L})$ of $\mathcal{L}$ is defined by:

\begin{align}\label{9.2}  \tau(L)=\sum_{1\leq i<j\leq \ell} K(e_i\wedge e_j),\end{align}
where $e_{1},\ldots,e_{\ell}$ is an orthonormal basis of $\mathcal{L}$.

\subsection{Chen invariants}\label{S9.1}

Let $n$ be a positive integer $\geq 3$. For a positive integer $k\leq \frac{n}{2}$, let $\mathcal S(n,k)$ denote the set  consisting of $k$-tuples $(n_1,\ldots,n_k)$ of integers $\geq 2$ such that  $n_1< n\;\; and \;\; n_1+\cdots+n_k\leq n.$  
Put ${\mathcal S}(n)=\cup_{k\geq 1} \mathcal S(n,k)$.

 For a given point $p$ in a Riemannian $n$-manifold $\bar M^{n}$ and each $(n_1,\ldots,n_k)\in\mathcal S(n)$, the first author introduced in \cite{58,c98,c00} the following invariants:
\begin{align}\notag&\delta{(n_1,\ldots,n_k)}(p)=\tau(x)-\inf\{\tau(\mathcal{L}_1)+\cdots+ \tau(\mathcal{L}_k)\},
\end{align}
where $\mathcal{L}_1,\ldots,\mathcal{L}_k$ run over all $k$ mutually orthogonal subspaces of $T_p \bar M^{n}$ such that $\dim \mathcal{L}_j=n_j$ and $j=1,\ldots,k$. 
In particular, we have
\vskip.05in 

(a) $\delta(\emptyset)=\tau$,
\vskip.05in

(b) $\delta{(2)}=\tau-\inf K$, where $K$ is the sectional curvature,

\vskip.05in
(c) $\delta(n-1)(p)=\max Ric(p)$.

\begin{remark} {\rm $\delta(2)$ is known today as the first Chen invariant among all of the invariants $\delta{(n_1,\ldots,n_k)}$.}
\end{remark}

\subsection{Chen inequalities}\label{S9.2}

 \noindent For each $(n_1,\ldots,n_k)\in \mathcal S(n,k)$, we put 
\begin{align*}& a(n_1,\ldots,n_k)={1\over2} {{n(n-1)}}-{1\over2}\sum_{j=1}^k {{n_j(n_j-1)}},\\&
 b(n_1,\ldots,n_k)= {{n^2(n+k-1-\sum_j n_j)}\over{2(n+k-\sum_j n_j)}}.\end{align*}
The first author proved the following optimal universal inequalities (see \cite{c98,c00,c05,84}).

\begin{theorem} \label{T:9.2}   Let $N$ be an $n$-dimensional submanifold of a Riemannian manifold $\bar M^{m}$. Then, for each point $p\in N$ and each $k$-tuple $(n_1,\ldots,n_k)\in \mathcal S(n)$,  we have
\begin{align}\label{9.3}\delta{(n_1,\ldots,n_k)}(p) \leq  b(n_1,\ldots,n_k)\|H\|^2(p)+a(n_1,\ldots,n_k)\max\bar K(p),\end{align}  where $\|H\|^2$ is the squared mean curvature of $N$ and $\max\bar K(p)$ is the maximum of the sectional curvature function of $\bar M^m$ restricted to $2$-plane sections of the tangent space $T_pN$ at $p$.

The equality case of inequality \eqref{9.3} holds at  $p\in N$ if and only if the following  conditions hold:

{\rm (a)} There is an  orthonormal basis  $e_1,\ldots,e_n,\xi_{n+1},\ldots,\xi_{m}$ at $p$ such that  the shape operators of $N$ in $\bar M^m$ at $p$ take the following  form:
 \begin{align}\label{9.4}\font\b=cmr8 scaled \magstep2 \def\bigzerol{\smash{\hbox{ 0}}} \def\bigzerou{\smash{\lower.0ex\hbox{\b 0}}} A_{e_r}=\left( \begin{matrix} A^r_{1} & \hdots & 0 \\ \vdots  & \ddots& \vdots &\bigzerou \\ 0 &\hdots &A^r_k\\  \\&\bigzerou & &\mu_rI \end{matrix} \right),\quad  r=n+1,\ldots,m, \end{align}
where $I$ is an identity matrix and  $A^r_j$ is a symmetric $n_j\times n_j$  submatrix such  that 
\begin{align}\label{9.5}\hbox{\rm trace}\,(A^r_1)=\cdots=\hbox{\rm trace}\,(A^r_k)=\mu_r.\end{align}

{\rm (b)} For mutual orthogonal subspaces $\mathcal{L}_1,\ldots,\mathcal{L}_k\subset T_p  N$ satisfying
$\delta(n_1,\ldots,n_k)=\tau-\sum_{j=1}^k \tau (\mathcal{L}_j)$ at $p$, we have
$ \bar K(e_{\alpha_i},e_{\alpha_j})=\max \bar K(p)$
 for $\alpha_i\in \Gamma_i,\alpha_j\in \Gamma_j$ and $0\leq i\ne j\leq k$, where \begin{equation}\begin{aligned}\notag &\Gamma_0=\{1,\ldots, n_1\},\; \ldots,\;\Gamma_{k-1}=\{n_1+\cdots+n_{k-1}+1,\ldots, n_1+\cdots+n_k\},\\&\hskip1.2in  \Gamma_{k}=\{n_1+\cdots+n_{k}+1,\ldots, n\}.\end{aligned}\end{equation}
\end{theorem}

An important case of  Theorem \ref{T:9.2} is the following.

 \begin{theorem} \label{T:9.3} {\rm \cite{c98,c00}}
For an $n$-dimensional submanifold $N$ of a real space form $R^m(c)$ of constant curvature $c$, we have 
\begin{align}\label{9.6} \delta(n_1,\ldots,n_k)\leq b(n_1,\ldots, n_k)\|H\|^2+a(n_1,\ldots,n_k)c.\end{align}

The equality case of inequality \eqref{9.6} holds at a point $p\in N$ if and only if there is an orthonormal basis 
$e_1,\ldots,e_n,\xi_{n+1},\ldots,\xi_{m}$  such that  the shape operators at $p$ take the forms  \eqref{9.4} and \eqref{9.5}.
\end{theorem}

\section{Inequalities in Golden Riemannian Manifolds}\label{S10}

 Following Chen's inequalities, many researchers have studied Chen-type inequalities within golden Riemannian manifolds.

\subsection{Chen type inequality in golden Riemannian manifolds}\label{S10.1}

The following findings about Chen-type inequalities for slant submanifolds in golden Riemannian manifolds were discovered by Uddin and Choudhary in \cite{59}.

\begin{theorem} {\rm \cite{59}} \label{T:10.1}
The following inequality holds for any proper $\theta$-slant submanifold $N^n$ that is isometrically immersed in a locally golden product manifold $\bar{M}^m$.
\begin{equation}\label{10.1}
\begin{aligned}
\delta_N(p) & \leq \frac{(n-2)}{2}\left[\frac{n^2}{(n-1)}\|H\|^2+\frac{1}{10}\left(c_p+c_q\right)\{3(n+1)-2 \operatorname{Trace}(\phi)\}\right] \\
& +\frac{1}{10}\left(c_p+c_q\right)\left[(\operatorname{Trace}(T)+(4-n)) \cos ^2 \theta-\operatorname{Trace}^2(\phi)\right] \\
& +\frac{1}{4 \sqrt{5}}\left(c_p-c_q\right)(n-2)[2 \operatorname{Trace}(\phi)-(\mathrm{n}+1)],\;\;  p\in N^{n} .
\end{aligned}
\end{equation}
\end{theorem}

For the equality case,

\begin{theorem}\label{T:10.2} {\rm \cite{59}}
When all conditions of the above Theorem \ref{T:10.1} are met, equality in Equation \eqref{10.1} is achieved at $p \in N $ if and only if $\{e_1, \ldots, e_n, e_{n+1}, \ldots, e_m\}$, the shape operator $A$ has the following form:
\begin{equation}\label{10.2}
    A_{n+1}=\left(\begin{array}{ccccc}
c & 0 & 0 & \ldots & 0 \\
0 & d & 0 & \ldots & 0 \\
0 & 0 & c+d & \ldots & 0 \\
\vdots & \vdots & \vdots & \ddots & \vdots \\
0 & 0 & 0 & \ldots & 0 \\
0 & 0 & 0 & \ldots & c+d
\end{array}\right),\;\;  A_s=\left(\begin{array}{ccccc}
c_s & d_s & 0 & \ldots & 0 \\
d_s & -c_s & 0 & \ldots & 0 \\
0 & 0 & 0 & \ldots & 0 \\
\vdots & \vdots & \vdots & \ddots & \vdots \\
0 & 0 & 0 & \ldots & 0 \\
0 & 0 & 0 & \ldots & 0
\end{array}\right), 
\end{equation} for $n+2 \leq s \leq m$.
\end{theorem}

Then they calculated an inequality involving $\delta(n_1, \ldots, n_k)$.

\begin{theorem} \label{T:10.3}  {\rm \cite{59}}
In each proper $\theta$-slant submanifold $N^n$ immersed in $\bar{M}^m$, the following inequality is true: 
\begin{equation}\label{10.3}
\begin{aligned}
\delta\left(n_1, \ldots, n_k\right) & \leq T_3-\frac{1}{10}\left(c_p+c_q\right)\left\{\cos ^2 \theta+\operatorname{Trace}(\phi)\right\}\left(n-\sum_{j=1}^k n_j\right) \\
&\hskip-.4in  -\frac{1}{4 \sqrt{5}}\left(c_p-c_q\right)\left\{\left(n+\sum_{j=1}^k n_j\right)-2 \operatorname{Trace}(\phi)-1\right\}\left(n-\sum_{j=1}^k n_j\right) ,
\end{aligned}
\end{equation}
where $$T_3=d\left(n_1, \ldots, n_\mu\right)\|H\|^2+\frac{3}{10}\left(c_p+c_q\right) b\left(n_1, \ldots, n_k\right).$$

Additionally, the equality sign in \eqref{10.3} holds at a point $p \in N$  if and only if there exists an orthonormal basis $\left\{e_1, \ldots, e_n, e_{n+1}, \ldots, e_m\right\}$ and $A$ such that
\begin{equation*}
    A_{n+1}=\left(\begin{array}{ccccc}
a_1 & 0 & 0 & \ldots & 0 \\
0 & a_2 & 0 & \ldots & 0 \\
0 & 0 & a_3 & \ldots & 0 \\
\vdots & \vdots & \vdots & \ddots & \vdots \\
0 & 0 & 0 & \ldots & 0 \\
0 & 0 & 0 & \ldots & a_n
\end{array}\right),\; \; A_s=\left(\begin{array}{cccccc}
B_1^s & \ldots & 0 & 0 & \ldots & 0 \\
\vdots & \ddots & \vdots & \vdots & \ddots & \vdots \\
0 & \ldots & B_k^s & 0 & \ldots & 0 \\
0 & \ldots & 0 & c_s & \ldots & 0 \\
\vdots & \ddots & \vdots & \vdots & \ddots & \vdots \\
0 & \ldots & 0 & 0 & \ldots & c_s
\end{array}\right), 
\end{equation*}
for $s \in\{n+2, \ldots, m\}$, where $a_1, \ldots, a_n$ satisfy
$$
a_1+\cdots+a_{n_1}=\cdots=a_{n_1+\ldots n_{k-1}+1}+\cdots+a_{n_1+\ldots n_k}=a_{n_1+\ldots n_k+1}=\cdots=a_n,
$$
and $B_i^s$ is a symmetric $n_i \times n_i$ submatrix satisfying
$$
\operatorname{Trace}\left({B}_1^{\mathrm{s}}\right)=\cdots=\operatorname{Trace}\left({B}_k^{\mathrm{s}}\right)=\mathrm{c}_{\mathrm{s}} .
$$
\end{theorem}

\begin{remark} {\rm
Additionally, Uddin and Choudhary  have deduced in \cite{59} a special case of Theorem \ref{T:10.1} and Theorem \ref{T:10.3} for $\phi$-invariant submanifolds $N^n$ immersed in locally golden product manifold $\bar{M}^m$ and  inequalities for Ricci curvature tensor.}
\end{remark}

\subsection{$\delta$-Casorati curvature in golden Riemannian manifolds}\label{S10.2}

In 1890, Casorati \cite{60} introduced what is now termed the Casorati curvature for surfaces in a Euclidean $3$-space $E^3$. Casorati favored this curvature over Gauss curvature because the latter may vanish for surfaces that intuitively seem curved, whereas the former only vanishes at planar points. The Casorati curvature $C$ of a submanifold in a Riemannian manifold is defined generally as the normalized squared norm of the second fundamental form.Decu et al. introduced normalized Casorati curvatures $\delta_C(n - 1)$ and $\hat{\delta}_C (n - 1)$ in 2007 (refer to \cite{61}), aligning with the essence of $\delta$-invariants. In 2008, they extended normalized Casorati curvatures to generalized normalized $\delta$-Casorati curvatures $\delta_C(r; n-1)$ and $\hat{\delta}_C(r; n-1)$ in \cite{62}. Concurrently, they were able to ascertain the optimal inequality with respect to the (intrinsic) scalar curvature and the (extrinsic) $\delta$-Casorati curvature.

Certainly! We can recall the Weingarten and Gauss formulas in this context.

For a Riemannian manifold $(\bar{M}, \bar g)$ and a Riemannian submanifold $N$ isometrically immersed in $\bar{M}$, where $\bar{\nabla}$ and $\nabla$ are the Levi-Civita connections on $\bar{M}$ and $N$ respectively, and $h$ represents the second fundamental form of $N$, the Weingarten and Gauss formulas are as follows:
$$
\begin{gathered}
\bar{\nabla}_X Y=\nabla_X Y+h(X, Y), \\
\bar{\nabla}_X \xi'=-A_{\xi '} X+\nabla_X^{\perp} \xi',
\end{gathered}
$$
$\forall$ $X, Y \in \Gamma(T N)$ and $ \xi' \in \Gamma\left(T^{\perp} N\right)$, where $A_{\xi'}$ denotes the shape operator of $N$ associated with $\xi'$ and $\nabla^{\perp}$ represents the connection in the normal bundle.
The relationship between $A_{\xi'}$ and $h$ can then be recalled.
$$
\bar g\left(A_{\xi'} X, Y\right)=\bar g(h(X, Y), \xi') .
$$
The Gauss formula is written as
$$
\bar{R}(X, Y, Z, W)=R(X, Y, Z, W)-\bar g(h(X, W), h(Y, Z))+\bar g(h(X, Z), h(Y, W)),
$$
for any vector fields tangent to $N$, such as $X, Y, Z, W$.
Assume that the local orthonormal tangent frame is $\left\{e_1, \ldots, e_n\right\}$, and the local orthonormal normal frame is $\left\{e_{n+1}, \ldots, e_m\right\}$. The definition of the scalar curvature is
$$
\tau=\sum_{1 \leq i<j \leq n} R\left(e_i, e_j, e_j, e_i\right),
$$
and the normalized scalar curvature $\rho$ by
$$
\rho=\frac{2 \tau}{n(n-1)} .
$$
For $N$, the mean curvature vector $H$ is
$$
H=\frac{1}{n} \sum_{i=1}^n h\left(e_i, e_i\right) .
$$
The components of $h$ are
$$
h_{i j}^r=\bar g\left(h\left(e_i, e_j\right), e_r\right), \;\; \forall i, j \in\{1, \ldots, n\},\; \forall r \in\{n+1, \ldots, m\} .
$$
Then
$$
\|H\|^2=\frac{1}{n^2} \sum_{r=n+1}^m\left(\sum_{i=1}^n h_{i i}^r\right)^2
$$
and
$$
\|h\|^2=\sum_{r=n+1}^m \sum_{i, j=1}^n\left(h_{i j}^r\right)^2 .
$$
The Casorati curvature $C$ of $N$ is defined by
$$
C=\frac{1}{n}\|h\|^2 .
$$
Let $\mathcal{L}$ be a $l$-dimensional subspace of $T_p N, l \geq 2$, and assume that $p \in N$. The scalar curvature of $\mathcal{L}$ for an orthonormal basis $\left\{e_1, \ldots, e_l\right\}$ can be expressed as
$$
\tau(\mathcal{L})=\sum_{1 \leq i<j \leq t} R\left(e_i, e_j, e_j, e_i\right) .
$$
One defines
$$
C (\mathcal{L})=\frac{1}{t} \sum_{r=n+1}^m \sum_{i, j=1}^t\left(h_{i j}^r\right)^2 .
$$
Assume that a hyperplane of $T_p N$ is $\mathcal{L}$. Then, the normalized $\delta$-Casorati curvatures $\delta_c(n-1)$ and $\widehat{\delta}_c(n-1)$ are expressed by
$$
\begin{gathered}
{\left[\delta_C (n-1)\right]_p=\frac{1}{2} C_p+\frac{n+1}{2 n} \inf \{C(\mathcal{L})\},} \\
{\left[\widehat{\delta}_{C}(n-1)\right]_p=2 C_p-\frac{2 n-1}{2 n} \sup C(\mathcal{L})\} .}
\end{gathered}
$$

The generalized normalized $\delta$-Casorati curvatures of $N$ contain the following expression for any real number $r>0$.

$\text { If } 0<r<n(n-1) \text {, }$
$$
\left[\delta_C(r ; n-1)\right]_p=r C_p+\frac{1}{r n} \cdot \mathcal{A}_1 \inf \{C(\mathcal{L})\},
$$

and, if $r>n(n-1)$,

$$
\left[\widehat{\delta}_{C}(r ; n-1)\right]_p=r C_p+\frac{1}{r n} \cdot \mathcal{A}_1 \sup \{C(\mathcal{L})\},
$$
with $\mathcal{A}_1=(n-1)(n+r)\left(n^2-n-r\right)$ \cite{64}.

In modern differential geometry, the study of $\delta$-Casorati curvatures becomes a highly active research subject. Many researchers have obtained intriguing findings on $\delta$ -Casorati curvatures in the golden Riemannian manifold.

Choudhary and Park obtained the following results in \cite{63} regarding $\delta$-Casorati curvatures of slant submanifolds of locally product golden space forms.

\begin{theorem}\label{T188} {\rm  \cite{63}}
Given an $(n+m)$-dimensional locally golden product space form $\big(\bar{M}^{{n+m}}=M_p\left(c_p\right) \times M_q\left(c_q\right), g, \phi\big)$, let $N$ be a $n$-dimensional $\theta$-slant proper submanifold. Then we have:

{\rm (i)} The curvature $\delta_C(r ; n-1)$, which is the generalized normalized $\delta$-Casorati, satisfies
\begin{equation}\label{10.4}
\begin{aligned}
\rho \leq &\, \frac{\delta_C(r ; n-1)}{n(n-1)} -\left(\frac{(1-\psi) c_p-\psi c_q}{2 \sqrt{5}}\right)\times 
\\
&\left\{1+\frac{1}{n(n-1)} \operatorname{Trace}^2 \phi-\cos ^2 \theta\left\{\frac{1}{n-1}+\frac{1}{n(n-1)} \operatorname{Trace} \phi \right\}\right\} \\
& -\left(\frac{(1-\psi) c_p+\psi c_q}{4}\right) \frac{2}{n} \operatorname{Trace} \phi
\end{aligned}
\end{equation}
for any real number $r$ such that $0<r<n(n-1)$.

{\rm (ii)}  The generalized normalized $\delta$-Casorati curvature $\widehat{\delta}_C(r ; n-1)$ satisfies
\begin{equation}\label{10.5}
\begin{aligned}
\rho \leq &\, \frac{\widehat{\delta}_C(r ; n-1)}{n(n-1)}-\left(\frac{(1-\psi) c_p-\psi c_q}{2 \sqrt{5}}\right)\times \\
& \left\{1+\frac{1}{n(n-1)} \operatorname{Trace}^2 \phi-\cos ^2 \theta\left\{\frac{1}{n-1}+\frac{1}{n(n-1)} \operatorname{Trace} \phi \right\}\right\} \\
& -\left(\frac{(1-\psi) c_p+\psi c_q}{4}\right) \frac{2}{n} \operatorname{Trace} \phi
\end{aligned}
\end{equation}
for any real number $r>n(n-1)$.

Furthermore, if and only if $N$ is an invariantly quasi-umbilical submanifold with trivial normal connection in $\bar{M}$, then the equalities in the relations \eqref{10.4} and \eqref{10.5} hold, such that the shape operators $A_r, r \in$ $\{n+1, \ldots, n+m\}$, with respect to some orthonormal tangent frame $\left\{e_1, \ldots, e_n\right\}$ and orthonormal normal frame $\left\{e_{n+1}, \ldots, e_{n+m}\right\}$, have the following forms:
\begin{equation}
A_{n+1}=\left(\begin{array}{cccccc}
d & 0 & 0 & \ldots & 0 & 0 \\
0 & d & 0 & \ldots & 0 & 0 \\
0 & 0 & d & \ldots & 0 & 0 \\
\vdots & \vdots & \vdots & \ddots & \vdots & \vdots \\
0 & 0 & 0 & \ldots & d & 0 \\
0 & 0 & 0 & \ldots & 0 & \frac{n(n-1)}{r} d
\end{array}\right), \quad A_{n+2}=\cdots=A_{n+m}=0.
\end{equation}
\end{theorem}

In golden Riemannian space forms, Choudhary and Park in \cite{7} have also obtained sharp inequalities for $\phi$-invariant and $\phi$-anti-invariant
submanifolds as a consequence of Theorem \ref{T188}.

\begin{theorem}\label{T189}  {\rm  \cite{63}}
Consider $N$ be an n-dimensional invariant submanifold of a locally golden product space form $\big(\bar{M}^{{n+m}}=M_p\left(c_p\right) \times M_q\left(c_q\right), g, \phi\big)$. Then we have:

{\rm (i)}  The generalized normalized $\delta$-Casorati curvature $\delta_C(r ; n-1)$ satisfies
\begin{equation}\label{10.7}
\begin{aligned}
\rho \leq &\, \frac{\delta_C(r ; n-1)}{n(n-1)}-\left(\frac{(1-\psi) c_p-\psi c_q}{2 \sqrt{5}}\right)\left\{1+\frac{1}{n(n-1)} \operatorname{Trace}^2 \phi\right\} \\
& +\left(\frac{(1-\psi) c_p-\psi c_q}{2 \sqrt{5}}\right)\left\{\frac{1}{n-1}+\frac{1}{n(n-1)} \operatorname{Trace} \phi \right\}\\
&-\left(\frac{(1-\psi) c_p+\psi c_q}{4}\right) \frac{2}{n} \operatorname{Trace} \phi
\end{aligned}
\end{equation}
for any real number $r$ such that $0<r<n(n-1)$.

{\rm (ii)}  The generalized normalized $\delta$-Casorati curvature $\widehat{\delta}_C(r ; n-1)$ satisfies
\begin{equation}\label{10.8}
\begin{aligned}
\rho \leq &\, \frac{\widehat{\delta}_C(r ; n-1)}{n(n-1)}-\left(\frac{(1-\psi) c_p-\psi c_q}{2 \sqrt{5}}\right)\left\{1+\frac{1}{n(n-1)} \operatorname{Trace}^2 \phi\right\} \\
& +\left(\frac{(1-\psi) c_p-\psi c_q}{2 \sqrt{5}}\right)\left\{\frac{1}{n-1}+\frac{1}{n(n-1)} \operatorname{Trace} \phi \right\} \\
& -\left(\frac{(1-\psi) c_p+\psi c_q}{4}\right) \frac{2}{n} \operatorname{Trace} \phi
\end{aligned}
\end{equation}
for any real number $r>n(n-1)$.

Furthermore, the equalities in the relations \eqref{10.7} and \eqref{10.8} hold if and only if $N$ is an invariantly quasi-umbilical submanifold with trivial normal connection in $\bar{M}$, such that the shape operators $A_r, r \in$ $\{n+1, \ldots, n+m\}$, take the following forms with respect to some orthonormal tangent frame $\left\{e_1, \ldots, e_n\right\}$ and orthonormal normal frame $\left\{e_{n+1}, \ldots, e_{n+m}\right\}$:

\begin{equation}
A_{n+1}=\left(\begin{array}{cccccc}
d & 0 & 0 & \ldots & 0 & 0 \\
0 & d & 0 & \ldots & 0 & 0 \\
0 & 0 & d & \ldots & 0 & 0 \\
\vdots & \vdots & \vdots & \ddots & \vdots & \vdots \\
0 & 0 & 0 & \ldots & d & 0 \\
0 & 0 & 0 & \ldots & 0 & \frac{n(n-1)}{r} d
\end{array}\right), \quad A_{n+2}=\cdots=A_{n+m}=0 .
\end{equation}
\end{theorem}

\begin{theorem}\label{T190}  {\rm  \cite{63}}
Let $N$ be an n-dimensional anti-invariant submanifold within $(n+m)$-dimensional locally golden product space form $\left(\bar{M}=M_p\left(c_p\right) \times M_q\left(c_q\right), g, \phi\right)$. Thus

{\rm (i)}  The generalized normalized $\delta$-Casorati curvature $\delta_C(r ; n-1)$ satisfies

\begin{equation}\label{10.10}
\rho \leq \frac{\delta_C(r ; n-1)}{n(n-1)}+\left(-\frac{(1-\psi) c_p-\psi c_q}{2 \sqrt{5}}\right)
\end{equation}
for any real number $r$ such that $0<r<n(n-1)$

{\rm (ii)} The generalized normalized $\delta$-Casorati curvature $\widehat{\delta}_C(r ; n-1)$ satisfies
\begin{equation}\label{10.11}
\rho \leq \frac{\widehat{\delta}_C(r ; n-1)}{n(n-1)}+\left(-\frac{(1-\psi) c_p-\psi c_q}{2 \sqrt{5}}\right)
\end{equation}
for any real number $r>n(n-1)$.

Additionally, the equalities hold in the relations \eqref{10.10} and \eqref{10.11} if and only if $N$ is an invariantly quasi-umbilical submanifold with trivial normal connection in $\bar{M}$, such that the shape operators $A_r, r \in$ $\{n+1, \ldots, n+m\}$, take the following forms with respect to some orthonormal tangent frame $\left\{e_1, \ldots, e_n\right\}$ and orthonormal normal frame $\left\{e_{n+1}, \ldots, e_{n+m}\right\}$:
\begin{equation}
A_{n+1}=\left(\begin{array}{cccccc}
d & 0 & 0 & \ldots & 0 & 0 \\
0 & d & 0 & \ldots & 0 & 0 \\
0 & 0 & d & \ldots & 0 & 0 \\
\vdots & \vdots & \vdots & \ddots & \vdots & \vdots \\
0 & 0 & 0 & \ldots & d & 0 \\
0 & 0 & 0 & \ldots & 0 & \frac{n(n-1)}{r} d
\end{array}\right), \quad A_{n+2}=\cdots=A_{n+m}=0.
\end{equation}
\end{theorem}

\begin{theorem}\label{T191}  {\rm  \cite{63}}
Assume that the locally golden product space form of dimension $(n+m)$ is $\left(\bar{M}=M_p\left(c_p\right) \times M_q\left(c_q\right), g, \phi\right)$. For any n-dimensional $\theta$-slant proper submanifold $N$ of $\bar{M}$,

{\rm (i)}  the normalized $\delta$-Casorati curvature $\delta_C(n-1)$ satisfies
\begin{equation} \label{a}
\begin{aligned}
\rho \leq &\, \delta_C(n-1)-\left(\frac{(1-\psi) c_p-\psi c_q}{2 \sqrt{5}}\right)\left\{1+\frac{1}{n(n-1)} \operatorname{Trace}^2 \phi\right\} \\ & +\left(\frac{(1-\psi) c_p-\psi c_q}{2 \sqrt{5}}\right) \cos ^2 \theta\left\{\frac{1}{n-1}+\frac{1}{n(n-1)} \operatorname{Trace} \phi \right\} \\ & -\left(\frac{(1-\psi) c_p+\psi c_q}{4}\right) \frac{2}{n} \operatorname{Trace} \phi
\end{aligned}
\end{equation}

{\rm (ii)}  the normalized $\delta$-Casorati curvature $\widehat{\delta}_C(n-1)$ satisfies

\begin{equation} \label{b}
\begin{aligned}
\rho \leq \,& \widehat{\delta}_C(n-1)-\left(\frac{(1-\psi) c_p-\psi c_q}{2 \sqrt{5}}\right)\left\{1+\frac{1}{n(n-1)} \operatorname{Trace}^2 \phi\right\} \\
& +\left(\frac{(1-\psi) c_p-\psi c_q}{2 \sqrt{5}}\right) \cos ^2 \theta\left\{\frac{1}{n-1}+\frac{1}{n(n-1)} \operatorname{Trace} \phi \right\} \\
& -\left(\frac{(1-\psi) c_p+\psi c_q}{4}\right) \frac{2}{n} \operatorname{Trace} \phi .
\end{aligned}
\end{equation}

\noindent Regarding any  invariant submanifold $N$ of $\bar{M}$, we have

{\rm (i)} the normalized $\delta$-Casorati curvature $\delta_C(n-1)$ satisfies
\begin{equation} \label{c}
\begin{aligned}
\rho \leq & \, \delta_C(n-1) - \left(\frac{(1-\psi) c_p - \psi c_q}{2 \sqrt{5}}\right) \left\{1 + \frac{1}{n(n-1)} \operatorname{Trace}^2 \phi\right\} \\
& + \left(\frac{(1-\psi) c_p - \psi c_q}{2 \sqrt{5}}\right) \left\{\frac{1}{n-1} + \frac{1}{n(n-1)} \operatorname{Trace} \phi \right\} \\
& -\left(\frac{(1-\psi) c_p + \psi c_q}{4}\right) \frac{2}{n} \operatorname{Trace} \phi
\end{aligned}
\end{equation}

{\rm (ii)} the normalized $\delta$-Casorati curvature $\widehat{\delta}_C(n-1)$ satisfies
\begin{equation} \label{d}
\begin{aligned}
\rho \leq & \widehat{\delta}_C(n-1)-\left(\frac{(1-\psi) c_p-\psi c_q}{2 \sqrt{5}}\right)\left\{1+\frac{1}{n(n-1)} t^2 \phi\right\} \\
& +\left(\frac{(1-\psi) c_p-\psi c_q}{2 \sqrt{5}}\right)\left\{\frac{1}{n-1}+\frac{1}{n(n-1)} \operatorname{Trace} \phi \right\} \\
& -\left(\frac{(1-\psi) c_p+\psi c_q}{4}\right) \frac{2}{n} \operatorname{Trace} \phi .
\end{aligned}
\end{equation}

\noindent  For any $n$-dimensional anti-invariant submanifold $N$ of $\bar{M}$,

{\rm (i)} the normalized $\delta$-Casorati curvature $\delta_C(n-1)$ satisfies

\begin{equation} \label{e}
\rho \leq \delta_C(n-1)-\left(\frac{(1-\psi) c_p-\psi c_q}{2 \sqrt{5}}\right)
\end{equation}

{\rm (ii)} the normalized $\delta$-Casorati curvature $\widehat{\delta}_C(n-1)$ satisfies
\begin{equation} \label{f}
\rho \leq \widehat{\delta}_C(n-1)-\left(\frac{(1-\psi) c_p-\psi c_q}{2 \sqrt{5}}\right) .
\end{equation}

Additionally, \eqref{a}, \eqref{c}, and \eqref{e} hold as equalities if and only if the submanifold $N^n$ is invariantly quasi-umbilical with trivial normal connection in $\bar{M}$, In this case, the shape operators $A_r, r \in\{n+1, \ldots, n+m\}$ with respect to some orthonormal tangent frame $\left\{e_1, \ldots, e_n\right\}$ and orthonormal normal frame $\left\{e_{n+1}, \ldots, e_{n+m}\right\}$ fulfill the following requirements:
\begin{equation}
A_{n+1}=\left(\begin{array}{cccccc}
d & 0 & 0 & \ldots & 0 & 0 \\
0 & d & 0 & \ldots & 0 & 0 \\
0 & 0 & d & \ldots & 0 & 0 \\
\vdots & \vdots & \vdots & \ddots & \vdots & \vdots \\
0 & 0 & 0 & \ldots & d & 0 \\
0 & 0 & 0 & \ldots & 0 & 2 d
\end{array}\right), \quad A_{n+2}=\cdots=A_{n+m}=0,
\end{equation}

Furthermore, \eqref{b}, \eqref{d}, and \eqref{f} hold as equalities if and only if the submanifold $N^n$ is invariantly quasi-umbilical with trivial normal connection in $\bar{M}$, such that the shape operators $A_r, r \in\{n+1, \ldots, n+m\}$ with respect to some orthonormal tangent frame $\left\{e_1, \ldots, e_n\right\}$ and orthonormal normal frame $\left\{e_{n+1}, \ldots, e_{n+m}\right\}$
satisfy:
\begin{equation}
A_{n+1}=\left(\begin{array}{cccccc}
2 d & 0 & 0 & \ldots & 0 & 0 \\
0 & 2 d & 0 & \ldots & 0 & 0 \\
0 & 0 & 2 d & \ldots & 0 & 0 \\
\vdots & \vdots & \vdots & \ddots & \vdots & \vdots \\
0 & 0 & 0 & \ldots & 2 d & 0 \\
0 & 0 & 0 & \ldots & 0 & d
\end{array}\right), \quad A_{n+2}=\cdots=A_{n+m}=0.
\end{equation}
\end{theorem}

\begin{remark} {\rm
 Some sharp inequalities for slant submanifolds immersed in golden Riemannian space forms with a semi-symmetric metric connection have been deduced by Lee et al. in \cite{65}. Furthermore, they characterize submanifolds in the equality case. They conclude by talking about these inequalities for a few special submanifolds.} 
\end{remark}

Geometric inequalities for the Casorati curvatures on submanifolds on golden Riemannian manifolds with constant golden sectional curvature were established by Choudhary and Mihai in \cite{64}. Let $(\bar{M}, \phi,\bar g)$ be a locally decomposable golden Riemannian manifold with constant golden sectional curvature. Next, on a submanifold $N$, the following are the optimal inequalities for $\delta_C(r; n-1)$ and $\hat{\delta}_C(r; n-1)$:

\begin{theorem}\label{T193}  {\rm  \cite{64}}
Consider $M$ be an n-dimensional Riemannian manifold isometrically immersed in $(\bar{M}, \phi,\bar g)$ under the condition $\mathcal{A}_2=n(n-1)$, we have:

{\rm (i)} $\delta_C(r ; n-1)$ satisfies
\begin{equation}\label{10.21}
    \rho \leq \frac{\delta_C(r ; n-1)}{\mathcal{A}_2}+\frac{c}{3 \mathcal{A}_2}\left\{n^2-3 n+2\|\phi\|^2-2 n\|\phi\|\right\}
\end{equation}
if $0 < r < \mathcal{A}_2$

{\rm (ii)} $\hat{\delta}_C(r; n-1)$ satisfies
\begin{equation}\label{10.22}
    \rho \leq \frac{\widehat{\delta}_C(r ; n-1)}{\mathcal{A}_2}+\frac{c}{3 \mathcal{A}_2}\left\{n^2-3 n+2\|\phi\|^2-2 n\|\phi\|\right\}
\end{equation}
if $r>\mathcal{A}_2$.

Furthermore, if and only if $N$ is invariantly quasi-umbilical, the equality holds in \eqref{10.21} or \eqref{10.22}, There exists an orthonormal tangent frame $\left\{e_1, \ldots, e_n\right\}$ and an orthonormal normal frame $\left\{e_{n+1}, \ldots, e_m\right\}$ such that $A_r, r \in\{n+1, \ldots, m\}$, have the following forms, and the normal connection of $N$ in $\bar{M}$ is trivial.
\begin{equation}
    A_{n+1}=\left(\begin{array}{cccccc}
d & 0 & 0 & \ldots & 0 & 0 \\
0 & d & 0 & \ldots & 0 & 0 \\
0 & 0 & d & \ldots & 0 & 0 \\
\vdots & \vdots & \vdots & \ddots & \vdots & \vdots \\
0 & 0 & 0 & \ldots & d & 0 \\
0 & 0 & 0 & \ldots & 0 & \frac{1}{r} \cdot \mathcal{A}_2  \cdot d
\end{array}\right), \quad A_{n+2}=\cdots=A_m=0.
\end{equation}
\end{theorem}

\begin{theorem}\label{T:10.11} {\rm  \cite{64}}
     Let $N$ be an isometrically immersed n-dimensional submanifold in $(\bar{M}, \phi,\bar g)$.
Then:

{\rm (i)} $\delta_C(r ; n-1)$ satisfies
\begin{equation}\label{10.24}
    \rho^G \leq \frac{\delta_C(r ; n-1)}{\mathcal{A}_2}+\frac{c}{3 \mathcal{A}_2}\{\|\phi\|(3\|\phi\|-3-n)\},
\end{equation}
for $0 < r < \mathcal{A}_2$.

ii. $\hat{\delta}_C(r; n-1)$ satisfies

\begin{equation}\label{10.25}
   \rho^G \leq \frac{\widehat{\delta}_C(r ; n-1)}{\mathcal{A}_2}+\frac{c}{3 \mathcal{A}_2}\{\|\phi\|(3\|\phi\|-3-n)\},
\end{equation}
for $r>\mathcal{A}_2$.

Furthermore, $A_r$ achieves the following forms, and the equality holds in \eqref{10.24} or \eqref{10.25} if and only if $N$ meets the equality's criteria as stated in Theorem \ref{T:10.10}:
\begin{equation}
A_{n+1}=\left(\begin{array}{cccccc}
d & 0 & 0 & \ldots & 0 & 0 \\
0 & d & 0 & \ldots & 0 & 0 \\
0 & 0 & d & \ldots & 0 & 0 \\
\vdots & \vdots & \vdots & \ddots & \vdots & \vdots \\
0 & 0 & 0 & \ldots & d & 0 \\
0 & 0 & 0 & \ldots & 0 & \frac{1}{r} \mathcal{A}_2 d
\end{array}\right), \quad A_{n+2}=\cdots=A_m=0 .
\end{equation}
\end{theorem}

\begin{remark} {\rm
In \cite{64} Choudhary and Mihai have established the consequences of theorem \ref{T:10.10} and \ref{T:10.11} and obtained inequality cases for Casorati curvature on an anti-invariant submanifold in $(\bar{M}, \phi,\bar g)$. }
\end{remark}

\begin{remark} {\rm
    Regarding sharp inequalities concerning $\delta$-Casorati curvatures for slant submanifolds of golden Riemannian space forms, refer to \cite{66} as well.}
\end{remark}

\subsection{Wintgen type inequality in golden Riemannian manifolds}\label{S10.3}

In the context of four-dimensional Euclidean space, the Wintgen inequality represents a critical geometric inequality that incorporates the Gauss curvature, normal curvature, and squared mean curvature, all of which are intrinsic invariants. In 1979, P. Wintgen \cite{67} formulated this inequality to show that for any surface $M^2$ in $E^4$, the Gauss curvature $K$, the normal curvature $K^{\perp}$, and the squared mean curvature $\|H\|^2$ meet the following condition:
$$
\|H\|^2 \geq K + |K^\perp|
$$
and the equality is valid if and only if the ellipse of curvature of $M^2$ in $E^4$ is a circle. This result was further generalized by I. V. Guadalupe et al. in \cite{68} for arbitrary codimension $m$ in real
space forms $\bar{M}^{m+2} (c)$ as follows:
$$
\|H\|^2 + c \geq K + |K^\perp|.
$$
They also discussed the conditions under which equality is achieved.

De Smet, Dillen, Verstraelen, and Vrancken \cite{DDVV}
proposed an inequality for submanifolds in real space forms, referred to as the generalized Wintgen inequality or DDVV conjecture, which extends the Wintgen inequality. This conjecture was independently proven by Ge and Tang in \cite{6}. Different researchers have obtained DDVV inequality for various classes of submanifolds in various ambient manifolds in the recent years.
For slant, invariant, C-totally real, and Lagrangian submanifolds in golden Riemannian space forms, authors have obtained generalized Wintgen type inequalities in \cite{70} and discussed the equality cases.

For slant submanifolds the Generalized Wintgen inequality is as follows:

\begin{theorem}\label{T:10.14} {\rm \cite{70}} Let $N$ be a $n$-dimensional $\theta$-slant proper submanifold of a locally golden product space form $\big(\bar{M}=M_p\left(c_p\right) \times M_q\left(c_q\right), g, \phi\big)$. Then we have:
\begin{equation}
    \begin{aligned}
\rho_\eta \leq & \, \|\mathcal{H}\|^2-2 \rho-2\left(\frac{(1-\psi) c_p-\psi c_q}{2 \sqrt{5}}\right)\left\{1+\frac{1}{n(n-1)} \operatorname{Trace}^2 \phi\right\} \\
& +2\left(\frac{(1-\psi) c_p-\psi c_q}{2 \sqrt{5}}\right)\left\{\frac{1}{n-1}+\frac{1}{n(n-1)} \operatorname{Trace} P \right\} \\
& -\left(\frac{(1-\psi) c_p+\psi c_q}{4}\right) \frac{4}{n} \operatorname{Trace} \phi .
\end{aligned}
\end{equation} 
For scalar normal curvature, $\rho_\eta$ is used.
\end{theorem}

Next, Choudhary et al. have established the generalized Wintgen inequality for invariant submanifold of golden Riemannian space forms with the aid of Theorem \ref{T:10.14}.

\begin{theorem}\label{T:10.15}  {\rm \cite{70}}
consider $N$ an n-dimensional invariant submanifold within a locally golden product space form $\left(\bar{M}=M_p\left(c_p\right) \times M_q\left(c_q\right), g, \phi\right)$. Then
\begin{equation}
\begin{aligned}
\rho_\eta \leq &\, \|\mathcal{H}\|^2-2 \rho-2\left(\frac{(1-\psi) c_p-\psi c_q}{2 \sqrt{5}}\right)\left\{1+\frac{1}{n(n-1)} \operatorname{Trace}^2 \phi\right\} \\
& +2\left(\frac{(1-\psi) c_p-\psi c_q}{2 \sqrt{5}}\right)\left\{\frac{1}{n-1}+\frac{1}{n(n-1)} \operatorname{Trace} P \right\} \\
& -\left(\frac{(1-\psi) c_p+\psi c_q}{4}\right) \frac{4}{n} \operatorname{Trace} \phi .
\end{aligned}
\end{equation}
\end{theorem}

The generalized Wintgen inequality of a locally golden product space form for a Lagrangian submanifold.

\begin{theorem} \label{T:10.16}  {\rm \cite{70}}
Suppose $N$ as a Lagrangian submanifold in a locally golden product space form $\left(\bar{M}=M_p\left(c_p\right) \times M_q\left(c_q\right), g, \phi\right)$. Then
\begin{equation}
\begin{aligned}
\left(\rho^{\perp}\right)^2 \geq \,& \rho_N^2 -\frac{2}{n(n-1)}\left(\frac{(1-\psi) c_p-\psi c_q}{2 \sqrt{5}}\right)^2 \\
& -\frac{4}{n(n-1)}\left(\frac{(1-\psi) c_p-\psi c_q}{2 \sqrt{5}}\right)\left\{\left(\frac{(1-\psi) c_p-\psi c_q}{2 \sqrt{5}}\right)+\rho\right\}
\end{aligned}
\end{equation}
\end{theorem}

\section{Inequalities in Golden-Like Statistical Manifolds}\label{S11}

\subsection{Chen type inequality in golden-like statistical manifolds}\label{S11.1}

Some fundamental inequalities for the curvature invariants of statistical submanifolds in golden-like statistical manifolds have been obtained by Bahadir et al. in \cite{8}.

\begin{theorem} {\rm \cite{8}}
  A golden-like statistical manifold of dimension $m$ is given by $(\bar M,\bar g, \phi)$ and $N$ be its statistical submanifold of dimension $n$. Then, 
 \begin{equation}
\begin{aligned}
(\tau- K(\pi))-\,&(\tau_0-K_0(\pi))
\geq \\& -\left(\frac{(1-\psi) c_p-\psi c_q}{2 \sqrt{5}}\right)\left[n(n-2)+\operatorname{Trace}^2(\phi)-\operatorname{Trace}(\phi^*)\right]\\
&-\left(\frac{(1-\psi) c_p+\psi c_q}{4}\right) 2(n-1) \operatorname{Trace}(\phi) \\
& +\left(\frac{(1-\psi) c_p-\psi c_q}{2 \sqrt{5}}\right)[1+\Psi(\pi)+\Theta(\pi)]\\ 
&- \frac{n^2(n-2)}{4(n-1)}\left[\|H\|^2+\left\|H^*\right\|^2\right]+2 \hat{K}_0(\pi)-2 \hat{\tau}_0 .
\end{aligned}
 \end{equation}

\end{theorem}

\begin{corollary} {\rm \cite{8}}
Let $N$ be the totally real statistical submanifold of dimension $n$ of a golden-like statistical manifold $(\bar M,\bar g, \phi)$ of dimension $m$.
 Then,
\begin{equation}
\begin{aligned}
(\tau-K(\pi))-\left(\tau_0-K_0(\pi)\right) \geq & -\left(\frac{(1-\psi) c_p-\psi c_q}{2 \sqrt{5}}\right)[n(n-2)-1]\\
&\hskip-.6in -\frac{n^2(n-2)}{4(n-1)}\left[\|H\|^2+\left\|H^*\right\|^2\right]+2 \hat{K}_0(\pi)-2 \hat{\tau}_0 .
\end{aligned}
\end{equation}
\end{corollary}

\subsection{$\delta$-Casorati curvature in golden-like statistical manifolds}\label{S11.2}

 Bahadir et al. have deduced in \cite{8}  optimal relationships for the generalized normalized $\delta$-Casorati curvature of a  statistical submanifold in a golden-like statistical manifold.
 
\begin{theorem} {\rm \cite{8}}
 Let $N^n$ be a statistical submanifold in a golden-like statistical manifold $\bar M^m$. Then we have the following optimal relationships for the generalized normalized $\delta$-Casorati curvature, :
 
{\rm (i)} For any real number $r$, such that $0<r<n(n-1)$,
\begin{equation}
\begin{aligned}
\rho \leq &\; \frac{\delta_C^0(r ; n-1)}{n(n-1)}+\frac{1}{(n-1)} C^0-\frac{n}{(n-1)} g\left(H, H^*\right)-\frac{2 n}{n(n-1)}\left\|H^0\right\|^2 \\
& -\frac{1}{n(n-1)}\left(\frac{(1-\psi) c_p-\psi c_q}{2 \sqrt{5}}\right)\left[n(n-2)+\operatorname{Trace}^2(\phi)-\operatorname{Trace}\left(\phi^*\right)\right]\\
&-\frac{2}{n}\left(\frac{(1-\psi) c_p+\psi c_q}{4}\right) \operatorname{Trace}(\phi),
\end{aligned}
\end{equation}
where $\delta_C^0(r ; n-1)=\frac{1}{2}\left[\delta_C(r ; n-1)+\delta_C^*(r ; n-1)\right]$.

{\rm (ii)} For any real number $r>n(n-1)$,
\begin{equation}
\begin{aligned}
\rho \leq &\; \frac{\widehat{\delta}_C^0(r ; n-1)}{n(n-1)}+\frac{1}{(n-1)} C^0-\frac{n}{(n-1)} g\left(H, H^*\right)-\frac{2 n}{n(n-1)}\left\|H^0\right\|^2 \\
& -\frac{1}{n(n-1)}\left(\frac{(1-\psi) c_p-\psi c_q}{2 \sqrt{5}}\right)\left[n(n-2)+\operatorname{Trace}^2(\phi)-\operatorname{Trace}\left(\phi^*\right)\right]\\
& -\frac{2}{n}\left(\frac{(1-\psi) c_p+\psi c_q}{4}\right) \operatorname{Trace}(\phi),
\end{aligned}
\end{equation}
where $\widehat{\delta}_C^0(r ; n-1)=\frac{1}{2}\left[\widehat{\delta}_C(r ; n-1)+\widehat{\delta}_C^*(r ; n-1)\right]$.
\end{theorem}

\begin{corollary}  {\rm \cite{8}}
Given a golden-like statistical manifold $\bar M^m$, let $N^n$ be a totally real statistical submanifold of it. Then for the generalized normalized $\delta$-Casorati curvature, we get the following optimal relationships:

{\rm (i)} For any real number $r$, such that $0<r<n(n-1)$,
\begin{equation}
\begin{aligned}
\rho \leq \,&\frac{\delta_C^0(r ; n-1)}{n(n-1)}+C^0-\frac{n}{(n-1)} g\left(H, H^*\right)-\frac{2 n}{n(n-1)}\left\|H^0\right\|^2\\
&-\left(\frac{(1-\psi) c_p-\psi c_q}{2 \sqrt{5}}\right)\left(\frac{n-2}{n-1}\right),
\end{aligned}
\end{equation}
where $\delta_C^0(r ; n-1)=\frac{1}{2}\left[\delta_C(r ; n-1)+\delta_C^*(r ; n-1)\right]$.

{\rm (ii)} For any real number $r>n(n-1)$,
\begin{equation}
\begin{aligned}
\rho \leq \,&\frac{\hat{\delta}_C^0(r ; n-1)}{n(n-1)}+\frac{1}{(n-1)} C^0-\frac{n}{(n-1)} g\left(H, H^*\right)-\frac{2 n}{n(n-1)}\left\|H^0\right\|^2\\
&-\left(\frac{(1-\psi) c_p-\psi c_q}{2 \sqrt{5}}\right)\left(\frac{n-2}{n-1}\right),
\end{aligned}
\end{equation}
where $\hat{\delta}_C^0(r ; n-1)=\frac{1}{2}\left[\hat{\delta}_C(r ; n-1)+\widehat{\delta}_C^*(r ; n-1)\right]$.
\end{corollary}

\section{Inequalities in Golden Lorentzian Manifolds}\label{S12}

\subsection{$\delta$-Casorati curvature in golden Lorentzian manifolds}\label{S12.1}

In \cite{10}, Choudhary et al. have deduced sharp geometric inequalities which involve generalized
normalized $\delta$-Casorati curvatures concerning submanifolds of golden Lorentzian manifolds equipped with generalized symmetric metric $U$-connection, and obtained the following results:

\begin{theorem} {\rm \cite{10}} \label{T204}
Below are the inequalities for the submanifold $N^{n}$ of $\bar M^{m}$.

{\rm (i)}  for $\delta_{C}(r ; n-1)$, we have
\begin{equation}\label{12.1}
\begin{aligned}
\rho \leq & \,\frac{1}{\left(n^2-n\right)}\left[\delta_C(r ; n-1)\right] \\
& +\frac{(\mp \sqrt{5}+3) c_p+( \pm \sqrt{5}+3) c_q-10 \alpha^2}{10(n-1)}(n-\varepsilon) \\
& +\frac{( \pm \sqrt{5}-1) c_p+(\mp \sqrt{5}-1) c_q-10 \alpha \beta}{10\left(n^2-n\right)}[(2 n \varepsilon-2) \text{\rm Trace }\phi ] \\
& +\frac{c_p+c_q-5 \beta^2}{5\left(n^2-n\right)}\left[(\text{\rm Trace }\, \phi )^2-\operatorname{Trace}\, \phi-n \varepsilon\right],
\end{aligned}
\end{equation}
here the real number $r$ satisfies $n^{2}-n>r>0$;

{\rm (ii)}  for $\widehat{\delta}_{C}(r ; n-1)$, we have
\begin{equation}\label{12.2}
\begin{aligned}
\rho \leq &\, \frac{1}{\left(n^2-n\right)}\left[\widehat{\delta}_C(r ; n-1)\right] \\
& +\frac{(\mp \sqrt{5}+3) c_p+( \pm \sqrt{5}+3) c_q-10 \alpha^2}{10(n-1)}(n-\varepsilon) \\
& +\frac{(\pm \sqrt{5}-1) c_1+(\mp \sqrt{5}-1) c_2-10 \alpha \beta}{10\left(n^2-n\right)}[(2 n \varepsilon-2) \text{\rm Trace }\phi ] \\
& +\frac{c_p+c_q-5 \beta^2}{5\left(n^2-n\right)}\left[(\operatorname{Trace} \phi)^2-\operatorname{Trace} \phi-n \varepsilon\right],
\end{aligned}
\end{equation}
where $n^2-n<r$.
Furthermore, if the shape operator in an orthonormal frame $\left\{e_1, \ldots\right.$, $\left.e_n, e_{n+1}, \ldots, e_m\right\}$, can be expressed as follows, then the relations in Equations \eqref{12.1} and \eqref{12.2} become equalities:
$$
A_{n+1}=\left(\begin{array}{cccccc}
d & 0 & 0 & \ldots & 0 & 0 \\
0 & d & 0 & \ldots & 0 & 0 \\
0 & 0 & d & \ldots & 0 & 0 \\
\vdots & \vdots & \vdots & \ddots & \vdots & \vdots \\
0 & 0 & 0 & \ldots & d & 0 \\
0 & 0 & 0 & \ldots & 0 & \frac{d n}{r}(n-1)
\end{array}\right),\;\;  A_{n+2}=\cdots=A_m=0.
$$
\end{theorem}

\begin{corollary}\label{C205} {\rm \cite{10}} 
The following inequalities hold for any submanifold $N^{n}$ immersed in a locally golden product Lorentzian manifold $\bar M^{m}$ equipped with a generalized symmetric metric $U$-connection:

{\rm (i)}  for $\delta_{C}(r ; n-1)$, we have
\begin{equation}\label{12.3}
\begin{aligned}
\rho \leq \,& \delta_C(n-1) +\frac{(\mp \sqrt{5}+3) c_p+( \pm \sqrt{5}+3) c_q-10 \alpha^2}{10(n-1)}(n-\varepsilon) \\
& +\frac{( \pm \sqrt{5}-1) c_p+(\mp \sqrt{5}-1) c_q-10 \alpha \beta}{10\left(n^2-n\right)}[(2 n \varepsilon-2) \text{\rm Trace } \phi] \\
& +\frac{c_p+c_q-5 \beta^2}{5\left(n^2-n\right)}\left[(\text{\rm Trace }\, \phi)^2-\text{ \rm Trace }\, \phi-n\varepsilon\right]
\end{aligned}
\end{equation}
where $0 \leq r \leq n^2-n$;

{\rm (ii)} for $\widehat{\delta}_C(n-1)$, we have
\begin{equation}\label{12.4}
\begin{aligned}
\rho \leq\, & \widehat{\delta}_C(n-1) \\
& +\frac{(\mp \sqrt{5}+3) c_p+( \pm \sqrt{5}+3) c_q-10 \alpha^2}{10(n-1)}(n-\varepsilon) \\
& +\frac{( \pm \sqrt{5}-1) c_p+(\mp \sqrt{5}-1) c_q-10 \alpha \beta}{10\left(n^2-n\right)}[(2 n \varepsilon-2) \text{\rm Trace}\, \phi] \\
& +\frac{c_p+c_q-5 \beta^2}{5\left(n^2-n\right)}\left[(\text{\rm Trace}\, \phi)^2-\text{\rm Trace}\, \phi-n \varepsilon\right],
\end{aligned}
\end{equation}
where $n^2-n<r$.

In addition, Equations \eqref{12.3} and \eqref{12.4} hold for equality if for an orthonormal frame $\left\{e_1, \ldots, e_n, e_{n+1}, \ldots, e_m\right\}$, operator $A$ can be expressed as follows
\[
A_{n+1}=\left(\begin{array}{cccccc}
d & 0 & 0 & \ldots & 0 & 0 \\
0 & d & 0 & \ldots & 0 & 0 \\
0 & 0 & d & \ldots & 0 & 0 \\
\vdots & \vdots & \vdots & \ddots & \vdots & \vdots \\
0 & 0 & 0 & \ldots & d & 0 \\
0 & 0 & 0 & \ldots & 0 & 2 d
\end{array}\right), \quad A_{n+2}=\cdots=A_m=0 ;
\]

and
\[
A_{n+1}=\left(\begin{array}{cccccc}
2 d & 0 & 0 & \ldots & 0 & 0 \\
0 & 2 d & 0 & \ldots & 0 & 0 \\
0 & 0 & 2 d & \ldots & 0 & 0 \\
\vdots & \vdots & \vdots & \ddots & \vdots & \vdots \\
0 & 0 & 0 & \ldots & 2 d & 0 \\
0 & 0 & 0 & \ldots & 0 & d
\end{array}\right), \quad A_{n+2}=\cdots=A_m=0.
\]
\end{corollary}

In \cite{10}, some consequences of Theorem \ref{T204} have also been derived, which are as follows:

\begin{corollary} {\rm \cite{10}} 
We have the following for a Riemannian manifold $N^n$ isometrically immersed in $\bar M^m$:

\noindent {\rm (I)} For $\delta_C(r ; n-1)$ with $r \in\{0, \ldots, n(n-1)\}$ we have:

\vskip.05in
{\rm (a)} $\bar M^m$ is equipped with $\alpha$ semi-symmetric metric $U$-connection
\begin{equation}
\begin{aligned}
\rho \leq & \frac{1}{\left(n^2-n\right)}\left[\delta_C(r ; n-1)\right] \\
& +\frac{(\mp \sqrt{5}+3) c_p+( \pm \sqrt{5}+3) c_q-10 \alpha^2}{10(n-1)}(n-\varepsilon) \\
& +\frac{( \pm \sqrt{5}-1) c_p+(\mp \sqrt{5}-1) c_q}{10\left(n^2-n\right)}[(2 m \varepsilon-2) \text{\rm Trace } \phi] \\
& +\frac{c_p+c_q}{5\left(n^2-n\right)}\left[(\text{\rm Trace}\, \phi)^2-\text{\rm Trace} \,\phi-n \varepsilon\right]
\end{aligned}
\end{equation}

{\rm (b)} $\bar M^m$ is equipped with $\beta$ quarter symmetric metric $U$-connection
\begin{equation}
\begin{aligned}
\rho \leq & \frac{1}{\left(n^2-n\right)}\left[\delta_C(r ; n-1)\right] \\
& +\frac{(\mp \sqrt{5}+3) c_p+( \pm \sqrt{5}+3) c_q}{10(n-1)}(n-\varepsilon) \\
& +\frac{( \pm \sqrt{5}-1) c_p+(\mp \sqrt{5}-1) c_q}{10\left(n^2-n\right)}[(2 n \varepsilon-2) \text{\rm Trace } \phi] \\
& +\frac{c_p+c_q-5 \beta^2}{5\left(n^2-n\right)}\left[(\text{\rm Trace}\, \phi)^2-\text{\rm Trace}\,\phi-n \varepsilon\right]
\end{aligned}
\end{equation}

{\rm (c)} $\bar M^m$ is equipped with semi-symmetric metric $U$-connection
\begin{equation}
\begin{aligned}
\rho \leq & \frac{1}{\left(n^2-n\right)}\left[\delta_C(r ; n-1)\right] \\
& +\frac{(\mp \sqrt{5}+3) c_p+( \pm \sqrt{5}+3) c_q-10}{10(n-1)}(n-\varepsilon) \\
& +\frac{( \pm \sqrt{5}-1) c_p+(\mp \sqrt{5}-1) c_q}{10\left(n^2-n\right)}[(2 n \varepsilon-2) \text{\rm Trace } \phi] \\
& +\frac{c_p+c_q}{5\left(n^2-n\right)}\left[(\text{\rm Trace } \phi)^2-\text{\rm Trace } \phi-n \varepsilon\right]
\end{aligned}
\end{equation}

{\rm (d)} $\bar M^m$ is equipped with quarter symmetric metric $U$-connection
\begin{equation}
\begin{aligned}
\rho \leq & \frac{1}{\left(n^2-n\right)}\left[\delta_C(r ; n-1)\right] \\
& +\frac{(\mp \sqrt{5}+3) c_p+( \pm \sqrt{5}+3) c_q}{10(n-1)}(n-\varepsilon) \\
& +\frac{( \pm \sqrt{5}-1) c_p+(\mp \sqrt{5}-1) c_q}{10\left(n^2-n\right)}[(2 n \varepsilon-2) \text{\rm Trace } \phi] \\
& +\frac{c_p+c_q-5}{5\left(n^2-n\right)}\left[(\text{\rm Trace } \phi)^2-\text{\rm Trace } \phi-n \varepsilon\right].
\end{aligned}
\end{equation}

\noindent {\rm (II)} For $\widehat{\delta}_C(r ; n-1)$ with $r>n(n-1)$, we  have:

\vskip.05in
{\rm (a)} $\bar M^m$ is equipped with $\alpha$ semi-symmetric metric $U$-connection
\begin{equation}
\begin{aligned}
\rho \leq \,& \frac{1}{\left(n^2-n\right)}\left[\widehat{\delta}_C(r ; n-1)\right] \\
& +\frac{(\mp \sqrt{5}+3) c_p+( \pm \sqrt{5}+3) c_q-10 \alpha^2}{10(n-1)}(n-\varepsilon) \\
& +\frac{( \pm \sqrt{5}-1) c_p+(\mp \sqrt{5}-1) c_q}{10\left(n^2-n\right)}[(2 n \varepsilon-2) \text{\rm Trace } \phi] \\
& +\frac{c_p+c_q}{5\left(n^2-n\right)}\left[(\text{\rm Trace } \phi)^2-\text{\rm Trace } \phi-n \varepsilon\right]
\end{aligned}
\end{equation}

{\rm (b)} $\bar M^m$ is equipped with $\beta$ quarter symmetric metric $U$-connection
\begin{equation}
\begin{aligned}
\rho \leq \,& \frac{1}{\left(n^2-n\right)}\left[\widehat{\delta}_C(r ; n-1)\right] \\
& +\frac{(\mp \sqrt{5}+3) c_p+( \pm \sqrt{5}+3) c_q}{10(n-1)}(n-\varepsilon) \\
& +\frac{( \pm \sqrt{5}-1) c_p+(\mp \sqrt{5}-1) c_q}{10\left(n^2-n\right)}[(2 n \varepsilon-2) \text{\rm Trace } \phi] \\
& +\frac{c_p+c_q-5 \beta^2}{5\left(n^2-n\right)}\left[(\text{\rm Trace } \phi)^2-\text{\rm Trace } \phi-n \varepsilon\right]
\end{aligned}
\end{equation}

{\rm (c)} $\bar M^m$ is equipped with semi-symmetric metric $U$-connection
\begin{equation}
\begin{aligned}
\rho \leq \,& \frac{1}{\left(n^2-n\right)}\left[\widehat{\delta}_C(r ; n-1)\right] \\
& +\frac{(\mp \sqrt{5}+3) c_p+( \pm \sqrt{5}+3) c_q-10}{10(n-1)}(n-\varepsilon) \\
& +\frac{( \pm \sqrt{5}-1) c_p+(\mp \sqrt{5}-1) c_q}{10\left(n^2-n\right)}[(2 n \varepsilon-2) \text{\rm Trace } \phi] \\
& +\frac{c_p+c_q}{5\left(n^2-n\right)}\left[(\text{\rm Trace } \phi)^2-\text{\rm Trace } \phi-n \varepsilon\right]
\end{aligned}
\end{equation}

{\rm (d)} $\bar M^m$ is equipped with quarter symmetric metric $U$-connection
\begin{equation}
\begin{aligned}
\rho \leq \,& \frac{1}{\left(n^2-n\right)}\left[\widehat{\delta}_C(r ; n-1)\right] \\
& +\frac{(\mp \sqrt{5}+3) c_p+( \pm \sqrt{5}+3) c_q}{10(n-1)}(n-\varepsilon) \\
& +\frac{( \pm \sqrt{5}-1) c_p+(\mp \sqrt{5}-1) c_q}{10\left(n^2-n\right)}[(2 n \varepsilon-2) \text{\rm Trace } \phi] \\
& +\frac{c_p+c_q-5}{5\left(n^2-n\right)}\left[(\text{\rm Trace}\, \phi)^2-\text{\rm Trace } \phi-n \varepsilon\right]
\end{aligned}
\end{equation}

Moreover, the relations in the above results become equalities if in some orthonormal frame $\{e_1, \ldots, e_n, e_{n+1}, \ldots, e_m\}$, the operators $A$ reduces to:
\[
A_{m+1}=\begin{pmatrix}
d & 0 & 0 & \cdots & 0 & 0 \\
0 & d & 0 & \cdots & 0 & 0 \\
0 & 0 & d & \cdots & 0 & 0 \\
\vdots & \vdots & \vdots & \ddots & \vdots & \vdots \\
0 & 0 & 0 & \cdots & d & 0 \\
0 & 0 & 0 & \cdots & 0 & \frac{d(n^2-n)}{r}
\end{pmatrix}, \quad A_{n+2}=\cdots=A_m=0.
\]
\end{corollary}

\begin{corollary}\label{C207} {\rm \cite{10}} 
When $N^n$ represents a Riemannian manifold isometrically immersed in a golden Lorentzian manifold $\bar M^m$ equipped with a g.s.m. U-connection, we have the following relations.

\noindent {\rm (I)} For $\delta_C(n-1)$ with $r \in\{0, \ldots\left(n^2-n\right)\}$, we have:

{\rm (a)} $\bar M^m$ is equipped with $\alpha$ semi-symmetric metric $U$-connection
\begin{equation}
\begin{aligned}
\rho \leq \,& \frac{1}{\left(n^2-n\right)}\left[\delta_C(n-1)\right] \\
& +\frac{(\mp \sqrt{5}+3) c_p+( \pm \sqrt{5}+3) c_q-10 \alpha^2}{10(n-1)}(n-\varepsilon) \\
& +\frac{( \pm \sqrt{5}-1) c_p+(\mp \sqrt{5}-1) c_q}{10\left(n^2-n\right)}[(2 n \varepsilon-2) \text{\rm Trace } \phi] \\
& +\frac{c_p+c_q}{5\left(n^2-n\right)}\left[(\text{\rm Trace } \phi)^2-\text{\rm Trace } \phi-n \varepsilon\right]
\end{aligned}
\end{equation}

{\rm (b)} $\bar M^m$ is equipped with $\beta$ quarter symmetric metric $U$-connection
\begin{equation}
\begin{aligned}
\rho \leq \, & \frac{1}{\left(n^2-n\right)}\left[\delta_C(n-1)\right] \\
& +\frac{(\mp \sqrt{5}+3) c_p+( \pm \sqrt{5}+3) c_q}{10(n-1)}(n-\varepsilon) \\
& +\frac{( \pm \sqrt{5}-1) c_p+(\mp \sqrt{5}-1) c_q}{10\left(n^2-n\right)}[(2 n \varepsilon-2) \text{\rm Trace } \phi] \\
& +\frac{c_p+c_q-5 \beta^2}{5\left(n^2-n\right)}\left[(\text{\rm Trace } \phi)^2-\text{\rm Trace } \phi-n \varepsilon\right]
\end{aligned}
\end{equation}

{\rm (c)} $\bar{M}^m$ is equipped with semi-symmetric metric $U$-connection
\begin{equation}
\begin{aligned}
\rho \leq & \frac{1}{\left(n^2-n\right)}\left[\delta_C(n-1)\right] \\
& +\frac{(\mp \sqrt{5}+3) c_p+( \pm \sqrt{5}+3) c_q-10}{10(n-1)}(n-\varepsilon) \\
& +\frac{( \pm \sqrt{5}-1) c_p+(\mp \sqrt{5}-1) c_q}{10\left(n^2-n\right)}[(2 n \varepsilon-2) \text{\rm Trace } \phi] \\
& +\frac{c_p+c_q}{5\left(n^2-n\right)}\left[(\text{\rm Trace } \phi)^2-\text{\rm Trace } \phi-n \varepsilon\right]
\end{aligned}
\end{equation}

{\rm (d)} $\bar M^m$ is equipped with quarter symmetric metric $U$-connection
\begin{equation}
\begin{aligned}
\rho \leq & \frac{1}{\left(n^2-n\right)}\left[\delta_C(n-1)\right] \\
& +\frac{(\mp \sqrt{5}+3) c_p+( \pm \sqrt{5}+3) c_q}{10(n-1)}(n-\varepsilon) \\
& +\frac{( \pm \sqrt{5}-1) c_p+(\mp \sqrt{5}-1) c_q}{10\left(n^2-n\right)}[(2 n \varepsilon-2) \text{\rm Trace } \phi] \\
& +\frac{c_p+c_q-5}{5\left(n^2-n\right)}\left[(\text{\rm Trace } \phi)^2-\text{\rm Trace } \phi-n \varepsilon\right]
\end{aligned}
\end{equation}

\noindent {\rm (II)} For $\widehat{\delta}_C(n-1)$ with $\left(n^2-n\right)<r$, we have:

{\rm(a)} $\bar M^m$ is equipped with $\alpha$ semi-symmetric metric $U$-connection
\begin{equation}
\begin{aligned}
\rho \leq \,& \frac{1}{\left(n^2-n\right)}\left[\widehat{\delta}_C(n-1)\right] \\
& +\frac{(\mp \sqrt{5}+3) c_p+( \pm \sqrt{5}+3) c_q-10 \alpha^2}{10(n-1)}(n-\varepsilon) \\
& +\frac{( \pm \sqrt{5}-1) c_p+(\mp \sqrt{5}-1) c_q}{10\left(n^2-n\right)}[(2 n \varepsilon-2) \text{\rm Trace } \phi] \\
& +\frac{c_p+c_q}{5\left(n^2-n\right)}\left[(\text{\rm Trace } \phi)^2-\text{\rm Trace } \phi-n \varepsilon\right]
\end{aligned}
\end{equation}

{\rm (b)} $\bar M^m$ is equipped with $\beta$ quarter symmetric metric $U$-connection
\begin{equation}
\begin{aligned}
\rho \leq \, & \frac{1}{\left(n^2-n\right)}\left[\widehat{\delta}_C(n-1)\right] \\
& +\frac{(\mp \sqrt{5}+3) c_p+( \pm \sqrt{5}+3) c_q}{10(n-1)}(n-\varepsilon) \\
& +\frac{( \pm \sqrt{5}-1) c_p+(\mp \sqrt{5}-1) c_q}{10\left(n^2-n\right)}[(2 n \varepsilon-2) \text{\rm Trace } \phi] \\
& +\frac{c_p+c_q-5 \beta^2}{5\left(n^2-n\right)}\left[(\text{\rm Trace } \phi)^2-\text{\rm Trace } \phi-n \varepsilon\right]
\end{aligned}
\end{equation}

{\rm (c)} $\bar M^m$ is equipped with semi-symmetric metric $U$-connection
\begin{equation}
\begin{aligned}
\rho \leq \, & \frac{1}{\left(n^2-n\right)}\left[\widehat{\delta}_C(n-1)\right] \\
& +\frac{(\mp \sqrt{5}+3) c_p+( \pm \sqrt{5}+3) c_q-10}{10(n-1)}(n-\varepsilon) \\
& +\frac{( \pm \sqrt{5}-1) c_p+(\mp \sqrt{5}-1) c_q}{10\left(n^2-n\right)}[(2 n \varepsilon-2) \text{\rm Trace } \phi] \\
& +\frac{c_p+c_q}{5\left(n^2-n\right)}\left[(\text{\rm Trace } \phi)^2-\text{\rm Trace } \phi-n \varepsilon\right]
\end{aligned}
\end{equation}

{\rm (d)} $\bar M^m$ is equipped with quarter symmetric metric $U$-connection
\begin{equation}
\begin{aligned}
\rho \leq \, & \frac{1}{\left(n^2-n\right)}\left[\widehat{\delta}_C(n-1)\right] \\
& +\frac{(\mp \sqrt{5}+3) c_p+( \pm \sqrt{5}+3) c_q}{10(n-1)}(n-\varepsilon) \\
& +\frac{( \pm \sqrt{5}-1) c_p+(\mp \sqrt{5}-1) c_q}{10\left(n^2-n\right)}[(2 n \varepsilon-2) \text{\rm Trace } \phi] \\
& +\frac{c_p+c_q-5}{5\left(n^2-n\right)}\left[(\text{\rm Trace } \phi)^2-\text{\rm Trace } \phi-n \varepsilon\right]
\end{aligned}
\end{equation}

Furthermore, the relations in the above results become equalities if in some orthonormal frame $\{e_1, \ldots, e_n, e_{n+1}, \ldots, e_m\}$, the operator $A$ reduces to:
\[
A_{m+1}=\begin{pmatrix}
d & 0 & 0 & \cdots & 0 & 0 \\
0 & d & 0 & \cdots & 0 & 0 \\
0 & 0 & d & \cdots & 0 & 0 \\
\vdots & \vdots & \vdots & \ddots & \vdots & \vdots \\
0 & 0 & 0 & \cdots & d & 0 \\
0 & 0 & 0 & \cdots & 0 & \frac{d(n^2-n)}{r}
\end{pmatrix}, \quad A_{n+2}=\cdots=A_m=0.
\]
\end{corollary}

\section{Further Structures on Golden Riemannian Manifolds}\label{S13}

\subsection{Integrability of golden Riemannian structure}\label{S13.1}

The significance of the golden structure on a Riemannian manifold stems from its association with pure Riemannian metrics. Given the connection between Riemannian golden structures and almost product structures, the $\varphi-$operator technique from the theory of almost product structures is applicable to golden structures. Consequently, the authors in the paper \cite{73} have established a new sufficient condition for the integrability of golden Riemannian structures, while also detailing certain characteristics of twin golden Riemannian metrics and the curvature features of locally decomposable golden Riemannian manifolds.

Consider a golden manifold denoted by $\bar M$ equipped with a golden structure $\phi$. For $\phi$ to be integrable, it is both necessary and sufficient to establish a torsion-free affine connection $\nabla$ such that the structure tensor $\phi$ remains covariantly constant under this connection. Furthermore, the integrability of $\phi$ correlates directly with the absence of the Nijenhuis tensor $\mathcal{N_\phi}$ \cite{2}. In \cite{73}, recent studies have investigated an additional potential sufficient condition for the integrability of golden structures within the framework of Riemannian manifolds.

\begin{theorem} {\rm \cite{73}}
A golden Riemannian manifold is denoted by $(\bar M, \phi,\bar g)$. If $ \varphi_\phi \bar g = 0.$, then $\phi$ is integrable.
\end{theorem}

\begin{corollary} {\rm \cite{73}}
Consider a golden Riemannian manifold $(\bar M, \phi,\bar g)$. $\varphi_\phi \bar g = 0$ is the same as $\nabla \phi =0$, where $\nabla$ is $\bar g$'s Levi-Civita connection. 
\end{corollary}

\begin{proposition} {\rm \cite{73}}
Given a golden Riemannian manifold $(\bar M, \phi,\bar g)$, let $F$ be the corresponding almost product structure. If $\varphi_F \bar g =0$, then the golden structure $\phi$ is integrable.
\end{proposition}

\begin{proposition} {\rm \cite{73}}
Let be a golden Riemannian manifold $(\bar M, \phi,\bar g)$. If and only if $\varphi_F \bar g = 0$, where $F$ is the corresponding almost product structure, then the manifold $\bar M$ is a locally decomposable golden Riemannian manifold.
\end{proposition}

\subsubsection{Twin Golden Riemannian Metrics}\label{S13.1.1}

A golden Riemannian manifold is denoted by $(\bar M, \phi,\bar g)$. The definition of the twin golden Riemannian metric is $$G (X, Y)= (\bar g \circ \phi) (X,Y)= \bar g (\phi X,Y)=\bar g (X, \phi Y)$$ for any vector field $X$ and $Y$ in $\bar M.$. It is simple to demonstrate that $G$ is a novel pure Riemannian metric:  
\begin{align*}
G(\phi X, Y) & = (\bar g \circ \phi)(\phi X, Y) = \bar g(\phi(\phi X), Y) =\bar g(\phi^2X, Y) \\
             & = \bar g(\phi X, Y) + \bar g(X, Y) = \bar g(X, \phi Y) +\bar g(X, Y) \\
             & =\bar g(X, (\phi + I)Y) =\bar g(X, \phi^2Y) \\
             & = (\bar g \circ \phi)(X, \phi Y) = G(X, \phi Y)
\end{align*}
which is called the twin metric of $\bar g.$ 

\begin{theorem} {\rm \cite{73}}
In a golden Riemannian manifold $(\bar M, \phi ,\bar g)$ the following holds: $$\varphi_\phi G =(\varphi_\phi \bar g) \circ \phi + \bar g \circ (\mathcal{N_\phi})$$
\end{theorem}

\begin{corollary}{\rm \cite{73}}
The criteria listed below are equivalent in a locally golden Riemannian manifold $(\bar M, \phi,\bar g)$:

(i) $\varphi_\phi \bar g=0$

(ii) $\varphi_\phi G =0$
\end{corollary}

\begin{theorem} {\rm \cite{73}}
Assume that the golden Riemannian manifold $(\bar M, \phi,\bar g)$ is locally decomposable. Thus, the Levi-Civita connection of the twin golden Riemannian metric $G$ and the Levi-Civita connection of the golden Riemannian metric $\bar g$ coincide.
\end{theorem}

\begin{theorem} {\rm \cite{73}}
The Riemannian curvature tensor field in a locally decomposable golden Riemannian manifold is a $\varphi-$tensor field.
\end{theorem}

\subsection{$s$-golden manifolds}\label{S13.2}

$s$-golden manifolds represent a fascinating category within almost golden Riemannian manifolds, explored by Gherici in \cite{75}.

Let $(\bar{M}, \phi, \bar{g})$ be an almost golden Riemannian manifold of dimension $n$. The tangent vector space $T_{p} \bar{M}$ splits as follows for each $p \in \bar{M}$:  $T_{p} \bar{M} = (D_{\psi^*})_{p} \oplus (D_{\psi})_{p}$, where $\psi = \frac{1 + \sqrt{5}}{2}$, $\psi^* = \frac{1 - \sqrt{5}}{2} = 1 - \psi$ and
\begin{align*}
&\left(D_\psi\right)_{\mathrm{p}}=\left\{\xi' \in T_{\mathrm{p}} \bar M: \phi_{\mathrm{p}} \xi'=\psi \xi'\right\}
\\& \left(D_{\psi^*}\right)_{\mathrm{p}}=\left\{\Phi \in T_{\mathrm{p}} \bar M: \phi_{\mathrm{p}} \Phi=\psi^* \Phi\right\} .
\end{align*}

\begin{definition} {\rm \cite{75}}
Let $\bar M$  be a differentiable manifold of dimension $(n + s)$. An almost $s$-golden structure on $\bar M$ is the data $\left(\phi,\left(\xi_\alpha, \rho_\alpha\right)_{\alpha=1}^s,\bar g\right)$ where:
$\xi_\alpha$ is global vector fields (called golden vector fields), $\rho_\beta$ is a differential 1-form on $\bar M$ such that $\rho_\beta\left(\xi_\alpha\right)=\delta_{\alpha \beta}$, where $\alpha, \beta \in$ $\{1, \ldots, s\}$, $\bar g$ is a Riemannian metric such that $\bar g\left(X, \xi_\alpha\right)=\rho_\alpha(X)$, $\phi$ is a tensor field of type $(1,1)$ satisfying
$$
\phi=\psi^* I+\sqrt{5} \sum_{\alpha=1}^s \rho_\alpha \otimes \xi_\alpha,
$$
$\forall$ vector fields $X$ on $\bar M$.
In addition, if $\phi$ is integrable, then $\left(\phi,\left(\xi_\alpha, \rho_\alpha\right)_{\alpha=1}^s,\bar g\right)$ is an $s$-golden structure and $\left(\bar M, \phi,\left(\xi_\alpha, \rho_\alpha\right)_{\alpha=1}^s,\bar g\right)$ is called an $s$-golden manifold.
\end{definition}

\begin{corollary} {\rm \cite{75}}
Any almost golden Riemannian structure $\phi$ that admits $s$ global unit eigenvectors associated with $\Phi$ is an almost s-golden structure.
\end{corollary}

\begin{proposition} {\rm \cite{75}}
Let $\left(\bar M, \phi,\left(\xi_\alpha, \rho_\alpha\right)_{\alpha=1}^s,\bar g\right)$ be an s-golden manifold and $U$ be a coordinate neighborhood on $\bar M$ and $\Phi_i$ any unit vector field on $U$ such that $\phi \Phi_i=\psi^* \Phi_i$, where $i \in\{1, \ldots, n\}$. Then, we may easily check that the set $\left\{\xi_\alpha, \Phi_i\right\}$ is a local orthonormal basis on $\bar M$.
\end{proposition}

Now,  a new sufficient condition of integrability for this class of structures is introduced as

\begin{theorem} {\rm \cite{75}}
Let $\left(\bar M, \Phi,\left(\xi_\alpha, \rho_\alpha\right)_{\alpha=1}^s,\bar g\right)$ be an almost s-golden structure. Then, $\phi$ is integrable if $\eta_\alpha$ are closed and $\left[\xi_\alpha, \xi_\beta\right]=0$ $\forall$$\alpha, \beta \in\{1, \ldots, s\}$.    
\end{theorem}

In \cite{75} author have also defined  two more special types of manifold $\mathcal{C}$-golden manifolds and $\mathcal{G}$-golden manifold.

\begin{definition} {\rm \cite{75}}
A $\mathcal{C}$-golden manifold is an almost trans-1-golden manifold $(\bar M^{n+1} , \phi ,\xi , \rho ,\bar{g})$ (or an almost trans-s-golden manifold of type $(1,0)$) which satisfies:
$$\nabla \phi=0.$$
\end{definition}

\begin{definition} {\rm \cite{75}}
A $\mathcal{G}$-golden manifold is an almost trans-1-golden manifold $(\bar M^{n+1} , \phi ,\xi , \rho ,\bar{g})$ (or an almost trans-s-golden manifold of type $(1,1)$) which satisfies:
$$(\nabla_X \phi ) Y = \sqrt{5} (\bar g (X,Y) \xi + \rho (Y) X -2 \rho (X) \rho (Y) \xi ),$$
$\forall$$X,Y$ vector fields in $\bar M.$
\end{definition}

\begin{remark} {\rm
The  geometric properties  and examples of $\mathcal{C}$-golden manifold and $\mathcal{G}$-golden manifold are also discussed in \cite{75}.}
\end{remark}

\subsection{Golden$^*$-manifold}\label{S13.3}

The $golden^*$-{\it manifold}  has been introduced in \cite{75} under the name ``$\mathcal{G}$-golden manifold" for a manifold of any dimension, but here the dimension is odd.

\begin{definition} {\rm \cite{76}}
A $golden^*$-manifold is an almost golden Riemannian almost contact manifold $(\bar M^{n+1} , \phi ,\xi , \rho ,\bar{g})$ which satisfies: $$(\nabla_X \phi ) Y = \sqrt{5} (\bar g (X,Y) \xi + \rho (Y) X -2 \rho (X) \rho (Y) \xi ),$$
$\forall$$X,Y$ vector fields in $\bar M.$
\end{definition}

\begin{lemma} {\rm \cite{76}}
Any $golden^*$-manifold is a golden Riemannian manifold.
\end{lemma}

\begin{proposition} {\rm \cite{76}}
On $golden^*$-manifold, the sectional curvature of all plane section containing $\xi$ is 1.
\end{proposition}

\begin{proposition} {\rm \cite{76}}
If the sectional curvature of any $\text{golden}^*$-manifold is a constant $c$, then $c = -1$.
\end{proposition}

\begin{proposition} {\rm \cite{76}}
On a $golden^*$-manifold,  the golden sectional curvature is 2.
\end{proposition}

\subsection{The $(\alpha , p)$-golden metrics manifolds}\label{S13.4}

Consider an even-dimensional manifold $\bar{M}$ that has a $\alpha$-structure $F_\alpha$. A Riemannian metric $\bar{g}$ is fixed such that
\begin{equation}\label{13.1}
\bar{g}\left(F_\alpha X, Y\right)=\alpha \bar{g}\left(X, F_\alpha Y\right),
\end{equation}
which is equivalent to
\begin{equation}\label{13.2}
    \bar{g}\left(F_\alpha X, F_\alpha Y\right)=\bar{g}(X, Y),
\end{equation}
for any vector fields $X, Y \in \Gamma(T \bar{M})$, where $\Gamma(T \bar{M})$ is the set of smooth sections of $T \bar{M}$ \cite{77}.

\begin{definition} {\rm \cite{77}}
The Riemannian metric $\bar{g}$ defined on an even dimensional manifold $\bar{M}$ endowed with an $\alpha$-structure $F_\alpha$ which verifies the equivalent identities \eqref{13.1} and \eqref{13.2}, is called a metric $\left(\alpha, F_\alpha\right)$-compatible.
Thus, one can obtain that the Riemannian metric $\bar{g}$ verifies the identity
\begin{equation}\label{13.3}
    \bar{g}\left(\phi_{\alpha, p} X, Y\right)-\alpha \bar{g}\left(X, \phi_{\alpha, p} Y\right)=\frac{p}{2}(1-\alpha) \bar{g}(X, Y),
\end{equation}
for any $X, Y \in \Gamma(T \bar{M})$.
Moreover, remark that $\bar{g}$ and $\left(\phi_{\alpha, p}\right)$ are related by
\begin{equation}\label{13.4}
    \bar{g}\left(\phi_{\alpha, p} X, \Phi_{\alpha, p} Y\right)=\frac{p}{2}\left(\bar{g}\left(\phi_{\alpha, p} X, Y\right)+\bar{g}\left(X, \phi_{\alpha, p} Y\right)\right)+p^2 \bar{g}(X, Y),
\end{equation}
for any $X, Y \in \Gamma(T \bar{M})$.
\end{definition}

\begin{definition} {\rm \cite{77}}
An almost $(\alpha, p)$-golden Riemannian manifold is a triple $\left(\bar{M}, \phi_{\alpha, p}, \bar{g}\right)$, where $\bar{M}$ is an even dimensional manifold, $\phi_{\alpha, p}$ is an almost $(\alpha, p)$-golden structure and $\bar{g}$ is a Riemannian metric which verifies identities \eqref{13.3} and \eqref{13.4}.
\end{definition}

\begin{proposition}\hskip-.03in  {\rm \cite{77}}
    If $\left(\bar{M}, \phi_{\alpha, p}, \bar{g}\right)$ is an almost $(\alpha, p)$-golden Riemannian manifold of dimension $2 m$, then the Trace of the $\phi_{\alpha, p}$ structure satisfies
$$\operatorname{Trace}\left(\phi_{\alpha, p}^2\right)=p \cdot \operatorname{Trace}\left(\phi_{\alpha, p}\right)+\frac{5 \alpha-1}{2} m p^2 .$$
\end{proposition}

\begin{definition} {\rm \cite{77}}
    If $\nabla$ is the Levi-Civita connection on $(\bar M,\bar g),$ then the covariant derivative $\nabla F_\alpha$ is a tensor field of type $(1,2),$ defined by $$(\nabla_X F_\alpha) Y := \nabla_X F_\alpha Y - F_\alpha \nabla_X Y,$$ for any $X , Y \in \Gamma (T \bar M).$ 
\end{definition}

\begin{remark} {\rm
   The necessary and sufficient criteria for a submanifold in an almost $(\alpha, p)$-golden Riemannian manifold to be an invariant submanifold have also been determined by Hretcanu and Crasmareanu in \cite{77}.}
\end{remark}

\noindent {\footnotesize Bang-Yen Chen }

\noindent {\footnotesize Department of Mathematics, Michigan State University, East Lansing, Michigan 48824-1027, USA.}

\noindent {\footnotesize e-mail: chenb@msu.edu}

\vskip.2in

\noindent {\footnotesize Majid Ali Choudhary}

\noindent {\footnotesize Department of Mathematics, School of Sciences, Maulana Azad National Urdu University, Hyderabad, India.}

\noindent {\footnotesize e-mail: majid.alichoudhary@gmail.com}

\vskip.2in
\noindent {\footnotesize Afshan Perween}

\noindent {\footnotesize Department of Mathematics, School of Sciences, Maulana Azad National Urdu University, Hyderabad, India.}

\noindent {\footnotesize e-mail: afshanperween99@gmail.com}

\end{document}